\documentclass{tran-l}

\usepackage{t1enc}
\usepackage[utf8]{inputenc}
\usepackage{amsthm,amsmath,amssymb}
\usepackage{graphicx}
\usepackage{enumerate}
\usepackage{hyperref}
\usepackage{bm}

\usepackage{amsfonts}
\usepackage{graphicx}
\usepackage{bm}
\usepackage{amsmath, amsthm, amssymb}
\usepackage{graphicx}
\usepackage{hyperref}
 \usepackage{relsize}

\usepackage{amsmath, amssymb, amscd, verbatim, xspace,amsthm}
\usepackage{latexsym, epsfig, color}

\newcommand{\mathbbm}[1]{\text{\usefont{U}{bbm}{m}{n}#1}}

\newcommand\eqd{\mathrel{\overset{\makebox[0pt]{\mbox{\normalfont\tiny d}}}{=}}}

\DeclareMathOperator{\rankk}{Rank}

\DeclareMathOperator{\Symm}{Symm}

\DeclareMathOperator{\Rep}{Re}
\DeclareMathOperator{\Aff}{MinC}
\DeclareMathOperator{\RowS}{RowSpace}
\DeclareMathOperator{\Cosets}{Cos}
\DeclareMathOperator{\cok}{cok}

\DeclareMathOperator{\Sur}{Sur}
\DeclareMathOperator{\Hom}{Hom}
\DeclareMathOperator{\Aut}{Aut}
\DeclareMathOperator{\tors}{tors}

\DeclareMathOperator{\corank}{co-rank}

\DeclareMathOperator{\SBP}{ symmetric, bilinear, perfect }

\DeclareMathOperator{\rank2}{Rank_2}

\newcommand{\Z}{\ensuremath{{\mathbb{Z}}}\xspace}
\renewcommand{\P}{\ensuremath{{\mathbb{P}}}}

\newcommand{\ra}{\rightarrow}

\newcommand\f{\mathfrak{f}}
\newcommand\tensor{\otimes}
\newcommand\isom{\simeq}

\newtheorem{theorem}{Theorem}[section]
\newtheorem{lemma}[theorem]{Lemma}
\newtheorem{corr}[theorem]{Corollary}

\newtheorem{prop}[theorem]{Proposition}

\theoremstyle{definition}
\newtheorem{definition}[theorem]{Definition}

\theoremstyle{remark}
\newtheorem{remark}[theorem]{Remark}

\numberwithin{equation}{section}

\def\rver{1}

\begin{document}

\title{The distribution of sandpile groups of random regular graphs}


\author{Andr\'as M\'esz\'aros}
\address{Central European University, Budapest, and  \newline 
\indent Alfr\'ed R\'enyi Institute of Mathematics, Budapest}
\email{Meszaros\_Andras@phd.ceu.edu}


\date{}

\dedicatory{}

\begin{abstract}
We study the distribution of the sandpile group of random \break d-regular
graphs. For the directed model, we prove that it follows the Cohen-Lenstra heuristics, that
is, the limiting  probability that the $p$-Sylow subgroup of the sandpile group is
a given $p$-group $P$, is proportional to $|\Aut(P)|^{-1}$. For finitely many
primes, these events get independent in the limit. Similar results hold for undirected random regular graphs, where for odd primes the limiting distributions are the ones given by Clancy, Leake and Payne.

This answers an open question of Frieze and Vu whether the adjacency matrix of a random regular graph is  invertible with high probability. Note that for directed graphs this was recently proved by Huang. It also gives an alternate proof of a theorem of Backhausz and Szegedy.

\end{abstract}

\maketitle

\section{Introduction}

We start by defining our random graph models. Let  $d\ge3$. The graph of a permutation $\pi$ consists of the directed edges $i\pi(i)$. The \textit{random directed graph $D_n$} is defined by taking the union of the graphs of $d$ independent uniform random permutations of $\{1,2,\dots,n\}$.    Thus, the adjacency matrix $A_n$ of $D_n$ is just obtained as $A_n=P_1+P_2+...+P_d$, where $P_1,P_2,\dots,P_d$ are independent uniform random $n\times n$ permutation matrices. 

For the undirected model, assume that $n$ is even. The \textit{random $d$-regular graph $H_n$} is obtained by taking the union of $d$ independent uniform random perfect matchings. The adjacency matrix of $H_n$ is denoted by $C_n$. 

The reduced Laplacian $\Delta_n$ of $D_n$ is obtained from $A_n-dI$ by deleting its last row and last column. The subgroup of $\mathbb{Z}^{n-1}$ generated by the rows of $\Delta_n$ is denoted by $\RowS(\Delta_n)$. The  group $\Gamma_n=\mathbb{Z}^{n-1}/\RowS(\Delta_n)$ is called the \textit{sandpile group} of $D_n$. If $D_n$ is strongly connected (which happens with high probability as $n\to \infty$), then  $\Gamma_n$ is a finite abelian group of order $|\det \Delta_n|$. Note that from the Matrix-Tree Theorem, $|\det \Delta_n|$ is the number of spanning trees in $D_n$ oriented towards the vertex $n$. For general directed graphs the sandpile group may depend on the choice of deleted row and column, but not in our case, because $D_n$ is Eulerian. The sandpile group of $H_n$ is defined the same way. Assuming that $H_n$ is connected, the order of the sandpile group is equal to the number of spanning trees in $H_n$.

Our main results are the following.

\begin{theorem}\label{CohenlenstraD}
Let $p_1,p_2,\dots,p_s$ be  distinct primes. Let $\Gamma_n$ be the sandpile group of $D_n$. Let $\Gamma_{n,i}$ be the $p_i$-Sylow subgroup of $\Gamma_n$. For $i=1,2,\dots, s$, let $G_i$ be a finite abelian $p_i$-group. Then 
\begin{equation}\label{Cohenlenstraeq}
\lim_{n\to\infty} \mathbb{P}\left(\bigoplus_{i=1}^s \Gamma_{n,i}\isom \bigoplus_{i=1}^s G_i\right)=\prod_{i=1}^s \left(|\Aut(G_i)|^{-1} \prod_{j=1}^\infty (1-p_i^{-j})\right).
\end{equation}
\end{theorem}

\begin{theorem}\label{CohenlenstraUD}
Let $\Gamma_n$ be the sandpile group of $H_n$. Again let $\Gamma_{n,i}$ be the $p_i$-Sylow subgroup of $\Gamma_n$, and for $i=1,2,\dots, s$, let $G_i$ be a finite abelian $p_i$-group. Assuming that $d$ is odd, we have
\begin{multline}\label{Cohenlenstraeq2}
\lim_{n\to\infty} \mathbb{P}\left(\bigoplus_{i=1}^s \Gamma_{n,i}\isom \bigoplus_{i=1}^s G_i\right)\\=\prod_{i=1}^s \left(\frac{|\{\phi:G_i \times G_i \to \mathbb{C}^* \SBP \}|}{|G_i||\Aut(G_i)|} \prod_{j=0}^\infty (1-p_i^{-2j-1})\right).
\end{multline} 
Assume that $d$ is even and $p_1=2$. Then the $2$-Sylow subgroup  of $\Gamma_n$ has odd rank\footnote{The rank of a group is the minimum number of generators.}. Furthermore, if we assume that $G_1$ has odd rank, then  
\begin{multline*}
\lim_{n\to\infty} \mathbb{P}\left(\bigoplus_{i=1}^s \Gamma_{n,i}\isom \bigoplus_{i=1}^s G_i\right)=\\ 2^{\rankk(G_1)}\prod_{i=1}^s \left(\frac{|\{\phi:G_i \times G_i \to \mathbb{C}^* \SBP \}|}{|G_i||\Aut(G_i)|} \prod_{j=0}^\infty (1-p_i^{-2j-1})\right).
\end{multline*}
\end{theorem}

The distribution appearing in  \eqref{Cohenlenstraeq} is the one that appears in the Cohen-Lenstra heuristics. 
It was introduced by Cohen and Lenstra \cite{cohlens}  in a conjecture on the distribution of class groups of quadratic number fields.  The distribution appearing in  \eqref{Cohenlenstraeq2} is a modified version of the distribution from the Cohen-Lenstra heuristics that was introduced by Clancy et al  \cite{clp13,clp14}. 

A recent breakthrough paper of Wood \cite{wood} shows that the sandpile group of
dense Erd\H{o}s-R\'enyi random graphs satisfies the latter heuristic. 
That is, Theorem \ref{CohenlenstraUD} says that in terms of the sandpile group, random
3-regular graphs exhibit the same level of randomness as dense
Erd\H{o}s-R\'enyi graphs.
 The conceptual explanation is that the random matrices coming from
both models mix the space extremely well, as we will see in Theorem \ref{FoFo22} for our model. 

We can gain information about the sandpile group by counting the surjective homomorphisms from it to a fixed finite abelian group $V$. For a random abelian group $\Gamma$ and a fixed finite abelian group $V$, we call the expectation $\mathbb{E} |\Sur(\Gamma,V)|$  the \textit{surjective $V$-moment} of $\Gamma$. Our next theorems determine the limits of the surjective moments of the sandpile groups for our random graph models. The convergence of these moments then implies Theorem \ref{CohenlenstraD} and Theorem \ref{CohenlenstraUD},  using the
work of Wood \cite{wood}.

\begin{theorem}\label{momentumokD}
Let $\Gamma_n$ be the sandpile group of $D_n$. For any finite abelian group $V$, we have
\[\lim_{n\to\infty} \mathbb{E}|\Sur(\Gamma_n,V)|=1.\]
\end{theorem}

Recall that the exterior power $\wedge^2 V$  is defined to be the quotient of $V\otimes V$ by the subgroup generated by elements of the form $v\otimes v$.

\begin{theorem}\label{momentumokUD}
Let $\Gamma_n$ be the sandpile group of $H_n$. Let $V$ be a finite abelian group.  If $d$ is odd, then
\[\lim_{n\to\infty} \mathbb{E}|\Sur(\Gamma_n,V)|=|\wedge^2 V|,\]
if $d$ is even, then
\[\lim_{n\to\infty} \mathbb{E}|\Sur(\Gamma_n,V)|=2^{\rank2(V)}|\wedge^2 V|,\]
where ${\rank2(V)}$ is the rank of the $2$-Sylow subgroup of $V$.
\end{theorem} 

These theorems are proved by using the fact that, when they are acting on $V^n$, the adjacency matrices $A_n$ and $C_n$ both exhibit strong mixing properties, described as follows: 
For $q=(q_1,q_2,\dots,q_n)\in V^n$,  the minimal coset in $V$ containing $q_1,q_2,\dots,q_n$ is denoted by $\Aff_q$. Note that $\Aff_q$ is the coset $q_n+V_0$ where $V_0$ is the subgroup of $V$ generated by $q_1-q_n,q_2-q_n,\dots,q_{n-1}-q_n$. The sum of the components of $q$ is denoted by $s(q)=\sum_{i=1}^n q_i$, and we define
\[R(q,d)=\{r\in (d\cdot \Aff_q)^n\quad|\quad s(r)=ds(q)\}.\footnote{By definition $d\cdot\Aff_q=\{g_1+g_2+\dots+g_d|g_1,g_2,\dots,g_d\in \Aff_q\}$.}\]
It is straightforward to check that $A_n q\in R(q,d)$. Let $U_{q,d}$ be a uniform random element of $R(q,d)$. Given two random variables $X$ and $Y$ taking values of the finite set $\mathcal{R}$, we define $d_{\infty}(X,Y)=\max_{r\in \mathcal{R}} |\mathbb{P}(X=r)-\mathbb{P}(Y=r)|$. We prove that the distribution of $A_n q$ is close to that of $U_{q,d}$ in the following sense.

\begin{theorem}\label{FoFo}
For $d\ge 3$, we have
\[\lim_{n\to\infty} \sum_{q\in V^n} d_\infty(A_n q,U_{q,d})=0.\]
\end{theorem} 

We have a similar theorem for $C_n$. 
For $q,w\in V^n$, we define \[<q\otimes w>=\sum_{i=1}^n q_i\otimes w_i.\] 

Furthermore, let $I_2=I_2(V)$ be the subgroup of $V\otimes V$ generated by the set $\{a\otimes b+b\otimes a|\quad a,b\in V\}$.  Let $\rank2(V)$ be the rank  
of the $2$-Sylow of $V$, and let $I=I(V)$ be the subgroup of $V\otimes V$ generated by all elements of the form $a\otimes a$ for $a\in V$. Note that $I_2$ is a subgroup of $I$ of index $2^{\rank2(V)}$.    Since the random matrix $C_n$ is symmetric and the diagonal entries are all equal to $0$, for any $q\in V^n$, we have $<q\otimes C_n q>\in I_2$.  Let us define $R^S(q,d)$ as
\[R^S(q,d)=\{r\in (d\cdot \Aff_q)^n\quad|\quad s(r)=ds(q)\text{ and }<q\otimes r>\in I_2\}.\]
It is clear from what is written above that $C_n q\in R^S(q,d)$. Similarly as before, let $U_{q,d}^S$ be a uniform random element of $R^S(q,d)$. Then, we have

\begin{theorem}\label{FoFo22}
For $d\ge 3$, we have
\[\lim_{n\to\infty} \sum_{q\in V^n} d_\infty(C_n q,U_{q,d}^S)=0.\]
\end{theorem} 

Note that the limits in Theorems \ref{momentumokD}, \ref{momentumokUD}, \ref{FoFo} and \ref{FoFo22} are uniform in $d$. See Section \ref{secuniform} for further discussion. However, until Section \ref{secuniform}, we never claim any uniformity over the choice of $V$ and $d$.

Recently, Huang \cite{Huang} considered a slightly different random $d$-regular directed graph model on $n$ vertices, the  configuration model introduced by Bollob\'as \cite{bollob}. Let $F_n$ be the adjacency matrix of this random graph. Huang proves that for a prime $p$ such that $\text{gcd}(p,d)=1$, we have

\[\mathbb{E} |\{0\neq x\in \mathbb{F}_p^n|\quad F_nx=0\}|=1+o(1),\]
as $n$ goes to infinity, where $F_n$ is  considered as a matrix over $\mathbb{F}_p$. Then he combines this with Markov's inequality to obtain that
\[\mathbb{P}(F_n\text{ is singular in }\mathbb{F}_p)\le \frac{1+o(1)}{p-1}.\]

Consequently, as a random matrix in $\mathbb{R}$,
\[\mathbb{P}(F_n\text{ is singular in }\mathbb{R})=o(1).\]
This solves an open problem of Frieze \cite{Frieze} and Vu \cite{Vu} for random regular bipartite graphs.

Using Theorem \ref{FoFo22}, we can answer this question in its original form.
\begin{theorem}\label{inverti}
For the adjacency matrix $C_n$ of $H_n$, we have
\[\mathbb{P}(C_n\text{ is singular in }\mathbb{R})=o(1).\]
\end{theorem}
Indeed, from Theorem \ref{FoFo22} with the choice of $V=\mathbb{F}_p$, it is straightforward to prove that for an odd prime $p$ such that $\text{gcd}(p,d)=1$, we have

\[\mathbb{E} |\{0\neq x\in \mathbb{F}_p^n|\quad C_nx=0\}|=1+o(1).\]
Therefore, the statement follows as above.

There are contiguity results \cite{cont1,cont2} which allow us to pass from one random $d$-regular graph model to another. In particular, Theorem \ref{inverti} also true for  uniform random $d$-regular graphs with even number of vertices. See also the work of Nguyen and Wood \cite{ngwood}. After the first version of this paper appeared online, Huang \cite{Huang2} also extended his results to the undirected configuration model, giving credit to this paper.

Theorem \ref{CohenlenstraUD} describes the local behavior of the sandpile group $\Gamma_n$ of $H_n$. Now we try to gain some global information on these groups. The next statement gives the asymptotic order of $\Gamma_n$. This was first proved by McKay \cite{mckay}, but it also follows from the more general theorem of Lyons \cite{lyons}. Let us choose $H_2,H_4,\dots$ independently. The torsion part of $\Gamma_n$ is denoted by $\tors(\Gamma_n)$.

\begin{theorem}[McKay, Lyons]\label{ThmLyons}
With probability $1$, we have
\[\lim_{n\to\infty} \frac{\log |\tors(\Gamma_n)|}{n}=\log\frac{(d-1)^{d-1}}{[d(d-2)]^{d/2-1}}.\]
\end{theorem}

Theorem \ref{momentumokUD} leads to the following statement on the rank of $\Gamma_n$.

\begin{theorem}\label{rankThm}
With probability $1$, we have
\[\lim_{n\to\infty}\frac{\rankk(\Gamma_n)}{n}=0.\]
\end{theorem}
Observe that $\rankk(\tors(\Gamma_n))=\max_{p\text{ is a prime}} \rankk_p(\tors(\Gamma_n))$, 
where \break $\rankk_p(\tors(\Gamma_n))$ 
is the rank of the $p$-Sylow subgroup of $\tors(\Gamma_n)$. Thus, this theorem  suggests that many primes should contribute to reach the growth described in Theorem \ref{ThmLyons}, but we do not have a definite result in this direction.

A conjecture of Ab\'ert and Szegedy \cite{Abert} states that if $G_1,G_2,\dots$ is a Benjamini-Schramm convergent sequence of finite graphs, then for any prime $p$ the limit
\[\lim_{n\to\infty}\frac{\corank_p G_n}{|V(G_n)|}\]
exists, here $\corank_p G_n=\dim\ker \text{Adj}(G_n)$, where $\text{Adj}(G_n)$ is the adjacency matrix of $G_n$ considered as a matrix over the finite field $\mathbb{F}_p$. One of the most common examples of a Benjamini-Scramm convergent sequence is the sequence of random d-regular graphs $H_n$. This means that if we choose $H_n$ independently, then with probability $1$, the sequence converges. Following along the lines of the proof of Theorem \ref{rankThm}, one can prove that
\[\lim_{n\to\infty}\frac{\max_{p\text{ is a prime}} \corank_p(H_n)}{n}=0\]
with probability $1$, which settles this special case of the conjecture, and we even get a uniform convergence in $p$. Note that this has been proved by  Backhausz and  Szegedy \cite{back} using a different method.

Theorem \ref{CohenlenstraD} follows from Theorem \ref{momentumokD} using the results of Wood \cite{wood} on the \textit{moment problem}. The general question is the following. Given a random finite abelian $p$-group $X$, is it true that the surjective $V$-moments of $X$ uniquely determine the distribution of $X$? Note that we can restrict our attention to the surjective $V$-moments, where $V$ is a $p$-group, because any other moment is $0$. Furthermore, is it true that if $X_1,X_2,\dots$ is a sequence of random abelian $p$-groups such that the surjective $V$-moments of $X_n$ converge to those of $X$, then the distribution of $X_n$ converge weakly to the distribution of $X$?  Ellenberg, Venkatesh and  Westerland \cite{ellenb} proved that the answer is affirmative for both questions in the special case when each surjective moment of $X$ is $1$. In this case $X$ has the distribution from the Cohen-Lenstra heuristic. Later, it was proved by Wood \cite{wood} that the answer is yes for both questions if the moments do not grow too fast, namely, if $\mathbb{E}|\Sur(X,V)|\le |\wedge^2 V|$ for any finite abelian $p$-group $V$. The proof generalizes the ideas of Heath-Brown~\cite{HeathBrown}. In \cite{wood} this is stated only in the special case, when the limiting surjective $V$-moments of $X$ are exactly $|\wedge^2 V|$, but in a later paper of Wood \cite{wood2} it is stated in its full generality above. In fact, Wood proved this theorem in a slightly more general setting. Instead of abelian $p$-groups, one can consider groups which are direct sums of finite abelian $p_i$-groups for a fixed finite set of primes. See Section~\ref{mixingboldist} for details. Note that for even $d$, the moments of the sandpile groups of $H_n$ are larger than the bounds above. But using the extra information that the $2$-Sylow subgroups have odd rank in this case, we can modify the arguments of Wood to obtain the convergence of probabilities. See Section \ref{parosparos}.

Now we discuss the Cohen-Lenstra heuristic in terms of random matrices over the $p$-adic integers. Let $\mathbb{Z}_p$ be the ring of $p$-adic integers. Given an $n\times m$ matrix $M$ over $\mathbb{Z}_p$ we define $\RowS(M)=\{xM|x\in \mathbb{Z}_p^n\}$. The \textit{cokernel} of $M$ is defined as $\cok(M)=\mathbb{Z}_p^m/\RowS(M)$. Freidman and Washington \cite{frwa} proved that if $M_n$ is an $n\times n$ random matrix over $\mathbb{Z}_p$, with respect to the Haar-measure, then $\cok(M_n)$ asymptotically follows the distribution from the Cohen-Lenstra heuristic, that is, for any finite abelian $p$-group $G$, we have
\[\lim_{n\to\infty} \mathbb{P}(\cok(M_n)\isom G)=|\Aut(G)|^{-1} \prod_{j=1}^\infty (1-p^{-j}).\]   
In fact this is true even in a more general setting. It is enough to assume that the entries of  $M_n$ are independent and they are not degenerate in a certain sense. This was proved by Wood \cite{wood2}. Her paper also contains similar results for non-square matrices. 

Bhargava, Kane, Lenstra, Poonen and Rains \cite{bkl} proved that the 
cokernels of Haar-uniform skew-symmetric random matrices over $\mathbb{Z}_p$ are asymptotically distributed according to  
Delaunay's heuristics. The following somewhat analogous result was obtained by Clancy,  Leake, Kaplan, Payne and Wood \cite{clp14}. Let $M_n$ be a Haar-uniform symmetric random matrix over $\mathbb{Z}_p$. Then, for any finite abelian $p$-group $G$, we have
\begin{multline}\label{szimformula}
\lim_{n\to\infty} \mathbb{P}(\cok(M_n)\isom G)\\= \frac{|\{\phi:G \times G \to \mathbb{C}^* \SBP \}|}{|G||\Aut(G)|} \prod_{j=0}^\infty (1-p^{-2j-1}).
\end{multline}
This is exactly the distribution appearing in Theorem \ref{CohenlenstraUD}.  
Note that this is not the original formula given in \cite{clp14}, but it can be easily deduced from it, see \cite{wood}. Here, a map $\phi: G \times G \to \mathbb{C}^*$ is called a symmetric, bilinear, perfect pairing if  \break (i) $\phi(x,y)=\phi(y,x)$, (ii)  $\phi(x,y+z)=\phi(x,y)\phi(x,z)$, and (iii) for $\phi_x(y)=\phi(x,y)$, we have $\phi_x\equiv 1$ if and only if $x=0$. We can give a more explicit formula for the limiting probability above by using the following fact from \cite{wood}. If $G=\bigoplus_{i} \mathbb{Z}/p^{\lambda_i}\mathbb{Z}$ with $\lambda_1\geq\lambda_2\geq\cdots$ and  
 $\mu$ is the transpose of the partition $\lambda$, then
\begin{multline}\label{explicite}
\frac{|\{\phi:G \times G \to \mathbb{C}^* \SBP \}|}{|G||\Aut(G)|}\\=
p^{-\sum_i \frac{\mu_i(\mu_i+1)}{2}} {\prod_{i=1}^{\lambda_1} \prod_{j=1}^{\lfloor \frac{\mu_i-\mu_{i+1}}{2} \rfloor} (1-p^{-2j})^{-1}}.
\end{multline}

Now we give a brief summary of results on distribution of sandpile groups. We already defined the Laplacian and the sandpile group of a $d$-regular graph, now we give the general definitions. We start by directed graphs. Let $D$ be a strongly connected directed graph on the $n$ element vertex set $V$. The Laplacian $\Delta$ of $D$ is an $n\times n$ matrix, where the rows and the columns are both indexed by $V$, and for $i,j\in V$, we have
\[\Delta_{ij}=
\begin{cases}
d(i,j)&\text{for }i\neq j,\\
d(i,i)-d_{\text{out}}(i)&\text{for }i=j.
\end{cases}\]
 Here $d(i,j)$ is the multiplicity of the directed edge $ij$, $d_{\text{out}}(i)$ is the out-degree of $i$, that is, $d_{\text{out}}(i)=\sum_{j\in V} d(i,j)$.  For $s\in V$, the reduced Laplacian $\Delta_s$ is obtained from $\Delta$ by deleting the row and column corresponding to $s$. The group $\Gamma_s=\mathbb{Z}^{n-1}/\RowS(\Delta_s)$ is called the \textit{sandpile group at vertex $s$}. The order of $\Gamma_s$ is the number of spanning trees in $D$ oriented towards $s$. Let us define\break $\mathbb{Z}_0^n=\{x\in \mathbb{Z}^n|\sum_{i=1}^n x_i=0\}$. Note that every row of $\Delta$ is in $\mathbb{Z}_0^n$. Thus the following definition makes sense. The group $\Gamma=\mathbb{Z}_0^n/\RowS(\Delta)$ is called the \textit{total sandpile group}. If $D$ is Eulerian, then all of these definitions 
of sandpile groups coincide, so it is justified to speak about the sandpile group of $D$. In fact, the converse of the above statement about Eulerian graphs is also true, see Farrel and Levine \cite{farlev}.

For an undirected graph $G$, let $D$ be the directed graph obtained from $G$ by replacing each edge $\{i,j\}$ of $G$ by the directed edges $ij$ and $ji$. Then $D$ is Eulerian. The sandpile group of $G$ is defined as the sandpile group of $D$. See \cite{jarai,SPS1,SPS2,SPS3} for more information on sandpile groups.

We already mentioned the result of Wood \cite{wood} on Erd\H{o}s-R\'enyi random graphs. Here we give more details. For $0\le \varrho\le 1$, the Erd\H{o}s-R\'enyi random graph $G(n,\varrho)$ is a graph on the vertex set $\{1,2,\dots,n\}$, such that for each pair of vertices, there is an edge connecting them with probability $\varrho$ independently. Let $p_1,p_2,\dots,p_s$ be distinct primes. Fix $0<\varrho<1$. Let $\Gamma_n$ be the sandpile group of $G(n,\varrho)$. Let $\Gamma_{n,i}$ be the $p_i$-Sylow subgroup of $\Gamma_n$, and for $i=1,2,\dots, s$, let $G_i$ be a finite abelian $p_i$-group. Then
\begin{multline*}
\lim_{n\to\infty} \mathbb{P}\left(\bigoplus_{i=1}^s \Gamma_{n,i}\isom \bigoplus_{i=1}^s G_i\right)\\=\prod_{i=1}^s \left(\frac{|\{\phi:G_i \times G_i \to \mathbb{C}^* \SBP \}|}{|G_i||\Aut(G_i)|} \prod_{j=0}^\infty (1-p_i^{-2j-1})\right).
\end{multline*}  
See Equation \eqref{explicite} for an even more explicit formula.

Koplewitz \cite{kopl} proved the analogous result for directed graphs. For $0\le \varrho\le 1$, the random directed graph $D(n,\varrho)$ is a graph on the vertex set $\{1,2,\dots,n\}$, such that for each ordered pair of vertices, there is a directed edge connecting them with probability $\varrho$ independently. Let $p_1,p_2,\dots,p_s$ be distinct primes. Fix $0<\varrho<1$. Let $\Gamma_n$ be the total sandpile group of $D(n,\varrho)$. Let $\Gamma_{n,i}$ be the $p_i$-Sylow subgroup of $\Gamma_n$, and for $i=1,2,\dots, s$, let $G_i$ be a finite abelian $p_i$-group. Then
\begin{equation*}
\lim_{n\to\infty} \mathbb{P}\left(\bigoplus_{i=1}^s \Gamma_{n,i}\isom \bigoplus_{i=1}^s G_i\right)=\prod_{i=1}^s  \frac{\prod_{j=2}^\infty (1-p_i^{-j})}{|G||\Aut(G)|}.
\end{equation*}
Note that, unlike what we would expect knowing the undirected case, this distribution is not the same as the one given in Theorem \ref{CohenlenstraD}  for the random directed $d$-regular graph $D_n$. A quick explanation is that $D_n$ is Eulerian, while $D(n,\varrho)$ is not. Indeed, the total sandpile group is defined as $\mathbb{Z}_0^n\isom \mathbb{Z}^{n-1}$ factored out by $n$ relations, so for a general directed graph, we expect that it behaves like the cokernel of a random $n\times (n-1)$ matrix. However, for an Eulerian graph these $n$ relations are linearly dependent, because their sum is zero, so we expect that the total sandpile group behaves like the cokernel of a random $(n-1)\times (n-1)$ matrix. The results above indeed support these intuitions.

\vspace{30pt}


\textbf{Acknowledgements} The author is grateful to Mikl\'os Ab\'ert for the useful discussions throughout the writing of this paper, to Melanie Wood for her comments and the proof of Lemma \ref{explicitnu}, and to Van Vu for pointing out relevant references. The exceptionally detailed and thorough reports of the anonymous referees were of great help in improving the presentation of the paper. The author was partially supported by the ERC Consolidator Grant 648017 and the Hungarian National
Research, Development and Innovation Office, NKFIH grant K109684. 

\vspace{30pt}

\textbf{The structure of the paper}

Section \ref{preli} contains the basic definitions that we need, including the notion of typical vectors.  
In Section \ref{mixing}, we investigate the distribution of $A_nq$, where $q$ is a typical vector. 
The results in this section allow us to handle the contribution of the typical vectors to the sum $\sum_{q\in V^n} d_\infty(A_n^{(d)}q,U_{q,d})$ in Theorem \ref{FoFo}, but we still need to control the contribution of the non-typical vectors. This is done in  Section~\ref{Secmom}. 
The connection between the mixing property of the adjacency matrix and the sandpile group is explained in Section \ref{mixingboldist}. In Section \ref{secuniform}, we prove that several results hold uniformly in $d$. Most of the paper deals with the directed random graph model, the necessary modifications for the undirected model are given in Section \ref{modificat} and Section~\ref{parosparos}. In Section \ref{sublin}, we prove Theorem \ref{rankThm}. 
At many points of the paper we need to estimate the probabilities of certain non-typical events, the proofs of these lemmas are collected in Section \ref{Bounding}.

\section{Preliminaries}\label{preli}
In most of the paper we will consider the directed model, and then later give the modifications of the arguments that are needed to be done for the undirected model.

Consider a vector $q=(q_1,q_2,...,q_n)\in V^n$. For a permutation $\pi$ of the set $\{1,2,\dots, n\}$, the vector $q_{\pi}=(q_{\pi(1)},q_{\pi(2)},\dots,q_{\pi(n)})$ is called a permutation of $q$. We write $q_1\sim q_2$ if $q_1$ and $q_2$ are permutations of each other. The relation $\sim$ is an equivalence relation, the equivalence class of $q$, i.e., the set of permutations of $q$ is denoted by $S(q)$. A random permutation of $q$ is defined as the random variable $q_{\pi}$, where $\pi$ is chosen uniformly from the set of all permutations, or equivalently, as a uniform random element of $S(q)$.  


Note that for $q\in V^n$, the equivalence class $S(q)$ can be described by $|V|$ non-negative integers summing up to $n$.  Namely,  for $c\in V$, we define

\[m_q(c)=|\{i\quad|\quad q_i=c\}|,\]
so $m_q$ can be considered as a vector in $\mathbb{R}^V$. 

Fix $\frac{1}{2}<\alpha<\beta<\gamma<\frac{2}{3}$. We keep these choices fixed throughout the whole paper. All the (explicit or implicit) constants  are allowed to depend on the choice of $\alpha,\beta$ and $\gamma$. However, since we view $\alpha,\beta$ and $\gamma$ as fixed, we will never emphasize this.

Note that if we choose a uniform random element $q$ of $V^n$, then the expectation of $m_q(c)$ is $\frac{n}{|V|}$ for any $c\in V$. This makes the following definition  quite natural.

\begin{definition}
A vector $q\in V^n$ is called $\alpha$-typical if $\left\|m_q-\frac{n}{|V|}\mathbbm{1}\right\|_{\infty}<n^{\alpha}$. Here $\mathbbm{1}$ is the all $1$ vector and $\|.\|_{\infty }$ is the maximum norm. 
\end{definition}     

Similarly, we can can define $\beta$-typical vectors. Note that, since $\alpha>\frac{1}2$, a uniform element of $V^n$ will be $\alpha$-typical with probability $1-o(1)$.

We write $A_n^{(d)}$ in place of $A_n$ to emphasize the value of $d$.

One of the key steps towards  Theorem \ref{FoFo} is the following theorem.
\begin{theorem}\label{mixingThm2}
For any fixed finite abelian group $V$ and $d\ge 3$, 
we have
\[\lim_{n\to\infty} |V|^n \sup_{q\in V^n\quad \alpha-\text{typical}} d_\infty(A_n^{(d)}q,U_{q,d})=0.\]
\end{theorem}
This will be an easy consequence of the following theorem. 

\begin{theorem}\label{mixingThm}
For any fixed finite abelian group $V$ and   $h\ge 2$,  we have
\[\lim_{n\to\infty} \sup_{\substack{q\in V^n\quad \alpha-\text{typical}\\r\in R(q,h)\quad \beta-\text{typical}}} \left|\mathbb{P}(A_n^{(h)}q=r) |V|^{n-1}-1\right|=0.\]
\end{theorem}

In the proofs we often need to consider $h$-tuples $Q=(q^{(1)},q^{(2)},\dots,q^{(h)})$ where each $q^{(i)}$ is a permutation of a fixed $q\in V^n$. Such $h$-tuples will be called \break $(q,h)$-tuples. Let $\mathcal{Q}_{q,h}$ be the set of $(q,h)$-tuples. A random $(q,h)$-tuple is a tuple \break $\bar{Q}=(\bar{q}^{(1)},\bar{q}^{(2)},\dots,\bar{q}^{(h)})$, where $\bar{q}^{(1)},\bar{q}^{(2)},...,\bar{q}^{(h)}$ are independent random permutations of $q$.

Whenever we use the symbols $Q$ and $\bar{Q}$, they stand for a $(q,h)$-tuple, and a random $(q,h)$-tuple respectively, even if this is not mentioned explicitly. The value of $q$ should be clear from the context.

 Sometimes, it will be convenient to view a $(q,h)$-tuple $Q$ as a vector $Q=(Q_1,Q_2,\dots, Q_n)$ in $\left (V^h\right)^n$, where $Q_i=(q^{(1)}_i,q^{(2)}_i,\dots,q^{(h)}_i)$. The vector $m_q$ was used to extract the important information from a vector $q\in V^n$, we do the same for $(q,h)$-tuples, that is,  for $t\in V^h$, we define
\[m_Q(t)=|\{i\quad|\quad Q_i=t\}|.\]

For a subset $S$ of $V^h$, the sum $\sum_{t\in S} m_Q(t)$ is denoted by $m_Q(S)$. Instead of~$S$, we usually just write the property that defines the subset $S$. For example, $m_Q(\tau_1=c)$ stands for  $m_Q(\{\tau\in V^h|\quad \tau_1=c\})$.


\begin{definition}
A $(q,h)$-tuple $Q$ or $m_Q$ itself will be called $\gamma$-typical if \[\left\|m_Q-\frac{n}{|V|^h}\mathbbm{1}\right\|_{\infty}<n^\gamma.\]
\end{definition}
The sum $\Sigma(Q)$ of a $(q,h)$-tuple $Q$ is defined as $\Sigma(Q)=\sum_{i=1}^h q^{(i)}$. 

Note that for a random $(q,h)$-tuple $\bar{Q}$, the distribution of $\Sigma(\bar{Q})$ is the same as that of $A_n^{(h)}q$.

Later in the paper  we will give asymptotic formulas that will be true uniformly in the following sense. 
\begin{definition}\label{unidef}
Let $X_1,X_2,...$ and $Y_1,Y_2,...$ be two sequences of finite sets, \break $P_n\subset X_n\times Y_n$, $f:\cup_{n=1}^{\infty} X_n\to \mathbb{R}$ and $g:\cup_{n=1}^{\infty} Y_n\to \mathbb{R}$. 

The term $f(x_n)\sim g(y_n)$ \textit{uniformly} for $(x_n,y_n)\in P_n$  means that
\[\lim_{n\to\infty} \sup_{(x_n,y_n)\in P_n}\left|\frac{f(x_n)}{g(y_n)}-1\right|=0.\]   
\end{definition}

The statement of Theorem \ref{mixingThm} then can be reformulated as \[\mathbb{P}(\Sigma(\bar{Q})=r)\sim \frac{1}{|V|^{n-1}}\] uniformly for any $\alpha$-typical $q\in V^n$ and $\beta$-typical $r\in R(q,h)$.

 \section{Behavior of typical vectors}\label{mixing}

In this section and the next section, we keep $V$ and $h$ fixed. 
All the (explicit or implicit) constants  are allowed to depend on $V$ and $h$. Moreover, whenever we claim the convergence of any quantity, it is meant that the convergence is only true for fixed $V$ and $h$. We never claim any uniformity over the choice of $V$ and $h$. Note that we deal with the question of uniformity in $d$ in Section \ref{secuniform} separately.  

We assume that $h\ge 2$ throughout this section.

\subsection{The proof of Theorem \ref{mixingThm}} 


We express the event $\Sigma(\bar{Q})=r$ as the disjoint union of smaller events, which can be handled more easily. Let 
\[\mathcal{M}(q,r)=\{m_Q\quad|\quad Q\in\mathcal{Q}_{q,h},\Sigma(Q)=r\}.\footnote{Here we omitted from the notation the dependence on $h$, later we will do this several times without mentioning it.}\]

Then the event $\Sigma(\bar{Q})=r$ can be written as the disjoint union of the events $(\Sigma(\bar{Q})=r)\wedge (m_{\bar{Q}}=m)$ where $m$ runs through $\mathcal{M}(q,r)$, so
\[\mathbb{P}(\Sigma(\bar{Q})=r)=\sum_{m\in\mathcal{M}(q,r)} \mathbb{P}((\Sigma(\bar{Q})=r)\wedge (m_{\bar{Q}}=m)).\]







Observe that $\mathcal{M}(q,r)$ consists of the non-negative integral points of a certain affine subspace $A(q,r)$ of $\mathbb{R}^{V^h}$. This affine subspace $A(q,r)$ is determined by linear equations expressing that whenever $\Sigma(Q)=r$ for a $(q,h)$-tuple  \break $Q=(q^{(1)},q^{(2)},\dots,q^{(h)})$, we have $m_{q^{(i)}}=m_q$ for every $i=1,2,\dots,h$ and \break $m_{\Sigma(Q)}=m_r$, as the following lemma shows.

For $t=(t_1,t_2,\dots,t_h)\in V^h$, we define $t_\Sigma$ as $t_\Sigma=\sum_{i=1}^h t_i$.

\begin{lemma}\label{euq12lemma}
Consider  $q,r\in V^n$. 
If $m\in \mathcal{M}(q,r)$, then $m$ is a non-negative integral vector satisfying the following linear equations:


\begin{align}\label{equ1}
m(\tau_i=c)&=m_q(c)  & \forall i\in\{1,2,\dots,h\}, c\in V,\\
\label{equ2}
m(\tau_\Sigma=c)&=m_r(c) & \forall c\in V. 
\end{align}

Now assume that $m$ is a nonnegative integral vector satisfying the equations above, then

\begin{align}\label{epr}
\mathbb{P}((\Sigma(\bar{Q})=r)\wedge (m_{\bar{Q}}=m))&=\frac{\prod_{c\in V} m_r(c)!}{\prod_{t\in V^h} m(t)!}\Bigg/\left(\frac{n!}{\prod_{c\in V} m_q(c)!}\right)^h\\&=
\frac{\prod_{c\in V} m(\tau_\Sigma=c)!}{\prod_{t\in V^h} m(t)!}\Bigg/\left(\frac{n!}{\prod_{c\in V} m_q(c)!}\right)^h.\nonumber
\end{align}

In particular, $\mathbb{P}((\Sigma(\bar{Q})=r)\wedge (m_{\bar{Q}}=m))>0$ so $m\in 
\mathcal{M}(q,r)$. Thus, $\mathcal{M}(q,r)$ is the set of non-negative 
integral points of the affine subspace $A(q,r)$ given by the linear equations 
above.    
\end{lemma}
\begin{proof}
We only give the proof of Equation \eqref{epr}, since all the other statements of the lemma are straightforward to prove. For $c\in V$, let \[I_c=\{i\in \{1,2,\dots,n\}|\quad r_i=c\},\] and let $W_c=\{t\in V^h|t_\Sigma=c\}$. Let $Q=(Q_1,Q_2,\dots,Q_n)\in \left(V^h\right)^n$.  Assume that $m$ is a nonnegative integral vector satisfying  Equation \eqref{equ1} and Equation \eqref{equ2} above.  Observe that $Q\in\mathcal{Q}_{q,h}$, $m_Q=m$ and $\Sigma(Q)=r$ if and only if for every $c\in V$, the sets
\[\left(\{i\in \{1,2,\dots,n\}\quad|\quad Q_i=t\}\right)_{t\in W_c}\]
give us a partition of $I_c$, such that for every $t\in W_c$, the size of the corresponding part is $m(t)$. 

Note that for any $c\in V$, we have
\[\frac{|I_c|!}{\prod_{t\in W_c} m(t)!}=\frac{m_r(c)!}{\prod_{t\in W_c} m(t)!}\]
such partitions of $I_c$.

Clearly, the total number $(q,h)$-tuples is
\[\left(\frac{n!}{\prod_{c\in V} m_q(c)!}\right)^h.\]
Putting everything together the statement follows. 
\end{proof}
The left hand sides of Equation \eqref{equ1} and Equation \eqref{equ2} in Lemma \ref{euq12lemma} do not depend on $q$ or $r$, therefore the affine subspaces $A(q,r)$ are all parallel for any choice of $q$ and $r$. Hence, for every $q,r_1,r_2\in V^n$, there is a translation that moves $A(q,r_1)$ to $A(q,r_2)$. There are many such translations, and we will use the one given in
the next lemma.

\begin{lemma}\label{eltolas}
For any $r_1,r_2\in V^n$, we define the vector  $v=v_{r_1,r_2}\in \mathbb{R}^{V^h}$  by
\[v(t)=\frac{m_{r_2}(t_\Sigma)-m_{r_1}(t_\Sigma)}{|V|^{h-1}}\]
for every $t\in V^h$. Then, for any $q\in V^h$, we have

\[A(q,r_1)+v_{r_1,r_2}=A(q,r_2).\]
\end{lemma}
\begin{proof}
It is enough to prove that $A(q,r_1)+v_{r_1,r_2}\subset A(q,r_2)$ or equivalently if $m$ satisfies Equation \eqref{equ1} and Equation \eqref{equ2} in Lemma \ref{euq12lemma} above for $r=r_1$, then $m'=m+v_{r_1,r_2}$ satisfies  Equation \eqref{equ1} and Equation \eqref{equ2}  for $r=r_2$. Observe that for any $i=1,2,\dots,h$ and $c,s\in V$, we have
\[|\{t\in V^h|\quad t_i=c,t_\Sigma=s\}|=|V|^{h-2}.\]
(Here we need to use that $h\ge 2$.) So we have
\begin{align*}
\sum_{\substack{t\in V^h\\ t_i=c}} m'(t)&=\sum_{\substack{t\in V^h\\ t_i=c}} m(t)+\sum_{\substack{t\in V^h\\ t_i=c}} v_{r_1,r_2}(t)\\&=
m_q(c)+\sum_{s\in V} |\{t\in V^h|\quad t_i=c,t_\Sigma=s\}|\frac{m_{r_2}(s)-m_{r_1}(s)}{|V|^{h-1}}\\&=m_q(c)+\frac{1}{|V|}\left(\sum_{s\in V} m_{r_2}(s)-\sum_{s\in V} m_{r_1}(s)\right)\\&=m_q(c)+\frac{1}{|V|}(n-n)=m_q(c),
\end{align*}
that is, Equation \eqref{equ1} is satisfied. Furthermore, for any $c\in V$, we have
\begin{align*}
\sum_{\substack{t\in V^h\\ t_\Sigma=c}} m'(t)&=\sum_{\substack{t\in V^h\\ t_\Sigma=c}} m(t)+\sum_{\substack{t\in V^h\\ t_\Sigma=c}} v_{r_1,r_2}(t)\\&= m_{r_1}(c)+|V|^{h-1} \frac{m_{r_2}(c)-m_{r_1}(c)}{|V|^{h-1}}=m_{r_2}(c),
\end{align*}
that is, Equation \eqref{equ2} is satisfied.
\end{proof}

Whenever $A(q,r)$ contains integral points, the integral points of $A(q,r)$ are placed densely, in the sense that there is a $D$, depending only on $h$ and $V$, such that for any point $x\in A(q,r)$, there is an integral point $y\in A(q,r)$ with $\|x-y\|_\infty<D$. Actually, this is a general fact as the following lemma shows.

\begin{lemma}
Let $A$ be an affine subspace of $\mathbb{R}^k$ which is given by a set of equations with rational coefficients. Assume that $A$ contains an integral point $p$. Then there is a $D$ such that for any point $x\in A$, there is an integral point $y\in A$ with \break $\|x-y\|_\infty<D$. For parallel subspaces, we can choose the same $D$.
\end{lemma}
\begin{proof}
Observe that we can write $A$ as $A=p+A_0$, where $A_0$ is a linear subspace generated by a set of rational vectors $\{a_1,a_2,\dots,a_\ell\}$. Multiplying these vectors with an  appropriate scalar, we may assume that they are all integral vectors. Let 
\[D=\sum_{i=1}^\ell \|a_i\|_\infty.\]
Note that $x-p\in A_0$, so $x-p=\sum_{i=1}^\ell \alpha_i a_i$ for some constants $\alpha_i$. Then \[y=p+\sum_{i=1}^\ell \lfloor\alpha_i\rfloor a_i\]
is an integral vector such that $\|x-y\|_\infty<D$.   
\end{proof}
For $c\in V$, let $w_c\in \mathbb{R}^{V^h}$ be such that $w_c(t)=1$ if $t_\Sigma=c$ and $w_c(t)=0$ otherwise. For $i=1,2,\dots, h$ and $c\in V$, let $u_{i,c}\in \mathbb{R}^{V^h}$ be such that $u_{i,c}(t)=1$ if $t_i=c$ and $u_{i,c}(t)=0$ otherwise.

\begin{lemma}\label{vanegesz}
If $r\in R(q,h)$, then $A(q,r)$ contains an integral point. 
\end{lemma}
\begin{proof}
We need to show that the system of linear equations given by Equation~\eqref{equ1} and Equation \eqref{equ2} admits an integral solution. Using the integral analogue of Farkas' lemma  \cite[Corollary 4.1a.]{schr}, we obtain that there exists an integral solution if and only if for every choice of rational numbers $0\le \gamma(i,c)<1$ \break ($i=1,2,\dots,h,\quad c\in V)$ and $0\le \delta(c)<1$ ($c\in V)$ such that 
\begin{equation}\label{cond1}
\sum_{i=1}^h \sum_{c\in V} \gamma(i,c) u_{i,c}+\sum_{c\in V} \delta(c) w_{c}\text{ is an integral vector}   
\end{equation}
the number $\sum_{i=1}^h \sum_{c\in V} \gamma(i,c) m_q(c)+\sum_{c\in V} \delta(c) m_r(c)$ is an integer. We project the rational numbers $\gamma(i,c)$ and $\delta(c)$ to the group $S^1=\mathbb{Q}/\mathbb{Z}$. From now on we work in the group $S^1$. The condition given in \eqref{cond1} translates as follows. For every $t\in V^h$,
\begin{equation}\label{equt1}
\sum_{i=1}^h \gamma(i,t_i)+\delta(t_{\Sigma})=0
\end{equation}
in the group $S^1$. We define $\gamma'(i,c)=\gamma(i,c)-\gamma(i,0)$ and $\delta'(c)=\delta(c)+\sum_{i=1}^h \gamma(i,0)$. Clearly $\gamma'(i,0)=0$. Moreover, from Equation \eqref{equt1} with $t=0$, we get that \break$\delta'(0)=0$. Equation \eqref{equt1} can be rewritten as
\[\sum_{i=1}^h \gamma'(i,t_i)+\delta'(t_{\Sigma})=0.\]
For every $i$ and $c$, if $t\in V^h$ is such that $t_i=c$ and $t_j=0$ for $i\neq j$, then we obtain that $\gamma'(i,c)=-\delta'(c)$. Therefore, Equation \eqref{equt1} can be once again rewritten as
\[\sum_{i=1}^h \delta'(t_i)=\delta'(t_{\Sigma})=\delta'\left(\sum_{i=1}^h t_i\right),\]
which means that $\delta'$ is a group homomorphism between $V$ and $\mathbb{Q}/\mathbb{Z}$. Thus, we get that
\begin{align*}
\sum_{i=1}^h &\sum_{c\in V} \gamma(i,c) m_q(c)+\sum_{c\in V} \delta(c) m_r(c)\\&=
\sum_{i=1}^h \sum_{c\in V} \left(\gamma'(i,c)+\gamma(i,0)\right) m_q(c)+\sum_{c\in V} \left(\delta'(c)-\sum_{i=1}^h \gamma(i,0)\right) m_r(c)\\&=
\sum_{i=1}^h \sum_{c\in V} -\delta'(c) m_q(c)+\sum_{c\in V} \delta'(c) m_r(c)\\&=-h \sum_{c\in V} \delta'(c) m_q(c)+\sum_{c\in V} \delta'(c) m_r(c)\\&=
-h\sum_{i=1}^n \delta'(q_i)+\sum_{i=1}^n \delta'(r_i)=\delta'\left(-h\cdot s(q)+s(r)\right)=\delta'(0)=0
\end{align*} 
using that $r\in R(q,h)$. That is, $\sum_{i=1}^h \sum_{c\in V} \gamma(i,c) m_q(c)+\sum_{c\in V} \delta(c) m_r(c)$ is indeed an integer.
\end{proof}

Suppose that $r_1,r_2\in R(q,h)$. Let $v=v_{r_1,r_2}$. Then there is an integral point $m_1$ in $A(q,r_1)$. Since $m_1+v\in A(q,r_2)$, there is an integral point $m_2$ in $A(q,r_2)$ such that  $\|m_1+v-m_2\|_\infty<D$. Set $\hat{v}=\hat{v}_{r_1,r_2}=m_2-m_1$, then $\|\hat{v}-v\|_{\infty}<D$ and the map $m\mapsto m+\hat{v}$ gives a bijection between the integral points of $A(q,r_1)$ and the integral points of $A(q,r_2)$.

For each $\alpha$-typical $q\in V^n$, fix an arbitrary $\beta$-typical $r_0=r_0(q)\in R(q,h)$, that is, let $r_0$ be any $\beta$-typical $r_0\in V^n$ such that $s(r_0)=h\cdot s(q)$. Set \[\mathcal{M}^*(q,r_0)=\left\{m\in \mathcal{M}(q,r_0)\quad\Big|\quad \left\|m-\frac{n}{|V|^h}\mathbbm{1}\right\|_\infty<2n^{\gamma}\right\}.\]

For any other $\beta$-typical $r\in R(q,h)$, we define
\[\mathcal{M}^*(q,r)=\{m+\hat{v}_{r_0,r}\quad|\quad m\in \mathcal{M}^*(q,r_0)\}\subset \mathcal{M}(q,r).\]

Observe that for large enough $n$, if both $r_0$ and $r$ are $\beta$-typical, then \[\|\hat{v}_{r_0,r}\|_{\infty}<D+\frac{2n^{\beta}}{|V|^{h-1}}<n^{\gamma}.\] Thus, using that the map $m\mapsto m+\hat{v}_{r_0,r}$ is a bijection between the integral points of $A(q,r_0)$ and the integral points of $A(q,r)$, we obtain that if $n$ is large enough, then for every $\alpha$-typical $q\in V^n$ and  $\beta$-typical $r\in R(q,h)$, we have

\begin{equation}\label{tartalmazza}
\left\{m\in \mathcal{M}(q,r)\quad\Big|\quad \left\|m-\frac{n}{|V|^h}\mathbbm{1}\right\|_\infty<n^{\gamma}\right\}\subset \mathcal{M}^*(q,r).
\end{equation}

Here the set on the left is just the set of the $\gamma$-typical elements of $\mathcal{M}(q,r)$.

The crucial point of our argument is the next lemma. 

\begin{lemma}\label{lemma20}
For an $\alpha$-typical $q\in V^n$, a $\beta$-typical $r\in R(q,h)$, $r_0=r_0(q)$ and $m\in \mathcal{M}^*(q,r_0)$, we have that
\[\mathbb{P}((\Sigma(\bar{Q})=r_0)\wedge (m_{\bar{Q}}=m))\sim \mathbb{P}((\Sigma(\bar{Q})=r)\wedge (m_{\bar{Q}}=m+\hat{v}_{r_0,r}))\]
uniformly in the sense of Definition \ref{unidef}.
\end{lemma} 

\begin{remark}
For clarity, we write out the definition  of the uniform convergence above. That is, Lemma \ref{lemma20} is equivalent with the statement that for any fixed $V$ and $h$, we have
\[\lim_{n\to\infty} \sup_{\substack{q\in V^n\quad\alpha\text{-typical}\\
m\in \mathcal{M}^*(q,r_0(q))\\r\in R(q,h)\quad\beta\text{-typical}}}\left|\frac{\mathbb{P}((\Sigma(\bar{Q})=r_0(q))\wedge (m_{\bar{Q}}=m))}{ \mathbb{P}((\Sigma(\bar{Q})=r)\wedge (m_{\bar{Q}}=m+\hat{v}_{r_0(q),r}))}-1\right|=0.\]
\end{remark}

To prove Lemma \ref{lemma20}, we need a few lemmas. 

The following approximation will be useful for Lemma  \ref{ujapprox}.
\begin{lemma}\label{factorial}
Fix $K(n)$ such that $K(n)=o\left(n^{\frac{2}{3}}\right)$. Then for  $|k|<K(n)$, we have
\[(n+k)!\sim \sqrt{2\pi n}\left(\frac{n}{e}\right)^n \exp\left(k\log n+\frac{k^2}{2n}\right)\]
uniformly. In other words, we have
\[\lim_{n\to\infty} \sup_{|k|<K(n)} \left|\frac{\sqrt{2\pi n}\left(\frac{n}{e}\right)^n\exp\left(k\log n+\frac{k^2}{2n}\right)}{(n+k)!}-1\right|=0.\]
\end{lemma}
\begin{proof}
Using Taylor's theorem with the Lagrange form of the remainder \break \cite[Theorem 5.15]{rudin} for the function $f(x)=x\log x$, we get that
\[\left|(n+k)\log(n+k)-\left(n\log n+(\log n+1)k+\frac{k^2}{2n}\right)\right|=\left|\frac{f^{(3)}(c)}{6}k^3\right|=\frac{|k|^3}{6c^2}\]
for some $c\in (n,n+k)$. This implies that
\[\lim_{n\to\infty} \sup_{|k|<K(n)} \left|(n+k)\log(n+k)-\left(n\log n+(\log n+1)k+\frac{k^2}{2n}\right)\right|= 0.\]
\begin{samepage}
It is also clear that
\[\frac{\sqrt{n+k}}{\sqrt{n}}\sim 1\]
uniformly for $|k|\le K(n)$. 
\end{samepage}

Recall that  Stirling's formula \cite[(8.22)]{rudin} states that
\[n!\sim \sqrt{2\pi n}\exp(n \log n-n).\]

If we put everything together, then we get that
\begin{samepage}
\begin{align*}
(n+k)!&\sim \sqrt{2\pi (n+k)}\exp\left((n+k)\log (n+k)-(n+k)\right)\\
&\sim \sqrt{2\pi n} \exp\left(\left(n\log n+(\log n+1)k+\frac{k^2}{2n}\right)-(n+k)\right)\\
&=\sqrt{2\pi n}\left(\frac{n}{e}\right)^n \exp\left(k\log n+\frac{k^2}{2n}\right)
\end{align*}
uniformly for $|k|\le K(n)$.
\end{samepage}
\end{proof}

Note that in the lemma above, we do not need to assume that $n$ is an integer, as long as $n+k$ is an integer. 

In the next lemma, we use the notation $a(n)=\sqrt{2\pi n}(\frac{n}{e})^n$.

\begin{lemma}\label{ujapprox}
For $q,r\in V^n$ and $m\in \mathcal{M}(q,r)$ such that $\left\|m-\frac{n}{|V|^h}\mathbbm{1}\right\|_\infty<3n^\gamma$, we have

\[\mathbb{P}((\Sigma(\bar{Q})=r)\wedge (m_{\bar{Q}}=m))\sim f(q) \exp\left(\frac{1}{2n}B\left(m-\frac{n}{|V|^h}\mathbbm{1},m-\frac{n}{|V|^h}\mathbbm{1}\right)\right)\]
uniformly, where \[f(q)=\left (\frac{n!}{\prod_{c\in  V} m_q(c)!}\right)^{-h} \frac{\left(a\left(\frac{n}{|V|}\right)\right)^{|V|}}{\left(a\left(\frac{n}{|V|^h}\right)\right)^{|V|^h}},\] and $B:\mathbb{R}^{V^h}\times \mathbb{R}^{V^h}\to \mathbb{R}$ is a bilinear form defined as
\[B(x,y)=|V|\sum_{c\in V}\left(\sum_{\substack{t\in V^h\\t_\Sigma=c}} x(t)\right)\left( \sum_{\substack{t\in V^h\\t_\Sigma=c}} y(t) \right)-|V|^h\sum_{t\in V^h} x(t)y(t).\]
Note that $f(q)$ does not depend on $r$ and $m$.
\end{lemma}
\begin{proof}
Recall that $\gamma<\frac{2}{3}$, so for any $t\in V^h$, Lemma \ref{factorial} can be applied to expand $m(t)!$ at the point $\frac{n}{|V|^h}$. Thus, we obtain the approximation

\[
m(t)!\sim a\left(\frac{n}{|V|^h}\right)\cdot \exp\left(\left(m(t)-\frac{n}{|V|^h}\right)\log \frac{n}{|V|^h}+\frac{|V|^h\left(m(t)-\frac{n}{|V|^h}\right)^2}{2n}\right). 
\]
Similarly, for every $c\in V$, by expanding  $m(\tau_\Sigma=c) !$ at the point $\frac{n}{|V|}$, we obtain the approximation 
\begin{multline*}
 m(\tau_\Sigma=c) !\sim \\ a\left(\frac{n}{|V|}\right) \cdot\exp\left(\left( \sum_{\substack{t\in V^h\\t_\Sigma=c}} m(t)-\frac{n}{|V|}\right)\log\frac{n}{|V|}+\frac{|V|\left(\sum_{\substack{t\in V^h\\t_\Sigma=c}}\left(m(t)-\frac{n}{|V|^h}\right)\right)^2}{2n}\right).
 \end{multline*}

Substituting these approximations in Equation \eqref{epr}, we obtain the statement.
\end{proof}


We made all the necessary preparations to prove Lemma \ref{lemma20}. 
\begin{proof}(Lemma \ref{lemma20})
It is easy to check that $w_c$  is in the radical of the bilinear form $B$, that is, $B(.,w_c)=B(w_c,.)=0$. ($w_c$ was defined before Lemma \ref{vanegesz}.) Since $v_{r_0,r}\in \text{Span}_{c\in V} w_c$, we get that $v_{r_0,r}$ is also in the radical. Observe that if $n$ is large enough, then  $\|\hat{v}_{r_0,r}\|_{\infty}<D+\frac{2n^{\beta}}{|V|^{h-1}}<n^{\gamma}$, so both $m$ and $m+\hat{v}_{r_0,r}$ satisfies the conditions of Lemma \ref{ujapprox}. It is also clear that $B(x,y)=O(\|x\|_\infty \|y\|_\infty)$. Thus,

\begin{align*}
&\frac{1}{2n} B \left(m+\hat{v}_{r_0,r}-\frac{n}{|V|^h}\mathbbm{1},m+\hat{v}_{r_0,r}-\frac{n}{|V|^h}\mathbbm{1}\right)\\&=
\frac{1}{2n}B\left(m+(\hat{v}_{r_0,r}-v_{r_0,r})+v_{r_0,r}-\frac{n}{|V|^h}\mathbbm{1},m+(\hat{v}_{r_0,r}-v_{r_0,r})+v_{r_0,r}-\frac{n}{|V|^h}\mathbbm{1}\right)\\&= \frac{1}{2n}\Bigg(B\left(m-\frac{n}{|V|^h}\mathbbm{1},m-\frac{n}{|V|^h}\mathbbm{1}\right)+2B\left(m-\frac{n}{|V|^h}\mathbbm{1},\hat{v}_{r_0,r}-v_{r_0,r}\right)\\&\qquad\qquad\qquad+B\left(\hat{v}_{r_0,r}-v_{r_0,r},\hat{v}_{r_0,r}-v_{r_0,r})\right)\Bigg)\\&= \frac{1}{2n}\left(B(m-\frac{n}{|V|^h}\mathbbm{1},m-\frac{n}{|V|^h}\mathbbm{1})+O(4Dn^\gamma+D^2)\right)\\&= \frac{1}{2n}B\left(m-\frac{n}{|V|^h}\mathbbm{1},m-\frac{n}{|V|^h}\mathbbm{1}\right)+O(n^{\gamma-1}).
\end{align*}
 Then, the statement follows from Lemma \ref{ujapprox}.
\end{proof}

From Lemma \ref{lemma20}, it follows immediately that for an $\alpha$-typical $q$ and  $\beta$-typical $r_1,r_2\in R(q,h)$, we have
\[\sum_{m\in \mathcal{M}^*(q,r_1)} \mathbb{P}((\Sigma(\bar{Q})=r_1)\wedge (m_{\bar{Q}}=m))\sim \sum_{m\in \mathcal{M}^*(q,r_2)} \mathbb{P}((\Sigma(\bar{Q})=r_2)\wedge (m_{\bar{Q}}=m))\]
uniformly, or equivalently
\begin{equation}\label{kovetk}
\mathbb{P}((\Sigma(\bar{Q})=r_1)\wedge (m_{\bar{Q}}\in \mathcal{M}^*(q,r_1)))\sim \mathbb{P}((\Sigma(\bar{Q})=r_2)\wedge (m_{\bar{Q}}\in \mathcal{M}^*(q,r_2)))
\end{equation} 
uniformly.

The content of the next lemma can be summarized as "only the typical events matter".
\newpage

\begin{lemma}\label{OTEM}
We have
 \begin{enumerate}[(i)]
 \item\label{OTEM1} A uniformly chosen element of $V^n$ is $\beta$-typical with probability $1-o(1)$.

 \item\label{OTEM2} There is a $C_1$ such that for any $\alpha$-typical $q\in V^n$, we have
 \[\mathbb{P}(\bar{Q}\text{ is not }\gamma-\text{typical})\le C_1\exp(-n^{2\gamma-1}/C_1).\]

In particular, for an $\alpha$-typical $q\in V^n$, we have $\mathbb{P}(\bar{Q}\text{ is }\gamma-\text{typical})\sim 1$ uniformly in the sense of Definition \ref{unidef}.

 \item\label{OTEM3} There is a $C_2$ such that for any $\alpha$-typical $q\in V^n$, we have
 \[\mathbb{P}(\Sigma(\bar{Q})\text{ is not }\beta-\text{typical})\le C_2\exp(-n^{2\beta-1}/C_2).\]

 In particular, for an $\alpha$-typical $q\in V^n$, we have $\mathbb{P}(\Sigma(\bar{Q})\text{ is }\beta-\text{typical})\sim 1$ uniformly in the sense of Definition \ref{unidef}.

 \item\label{OTEM4} The following holds
 \[\lim_{n\to\infty} \sup_{\substack{q\in V^n\quad \alpha-\text{typical}\\r\in R(q,h)\quad \beta-\text{typical}}} \mathbb{P}\left((\Sigma(\bar{Q})=r)\wedge(\bar{Q}\text{ is not }\gamma-\text{typical})\right)|V|^{n-1}=0.\] 
 \end{enumerate}
\end{lemma}
\begin{proof}

Part \eqref{OTEM1} can be proved using standard concentration results. We omit the details. To prove the other statements of Lemma \ref{OTEM}, we need the following result.

\begin{lemma}\label{interm}
Fix $K(n)$ such that $n^{\alpha}=o(K(n))$. There is a $C$ such that for any $\alpha$-typical  $q\in V^n$ and a random $(q,h)$-tuple $\bar{Q}$, we have
\[\mathbb{P}\left(\left\|m_{\bar{Q}}-\frac{n}{|V|^h}\mathbbm{1}\right\|_{\infty}\ge K(n)\right)\le C\exp\left(-\frac{K(n)^2}{Cn}\right).\]
\end{lemma}
\begin{proof}
Observe that for any $\alpha$-typical $q\in V^n$ and $t\in V^h$, we have
\[\left|n\prod_{i=1}^h \frac{m_q(t_i)}n-\frac{n}{|V|^h}\right|=O(n^{\alpha})=o(K(n)),\]
where the hidden constant does not depend on $q$ or $t$. Thus,  for an $\alpha$-typical  $q\in V^n$ and a $(q,h)$-tuple $Q$, if we have 
\[\left|m_{Q}(t)-\frac{n}{|V|^h}\right|\ge K(n)\]
for some $t\in V^h$, then
\[\left|m_{Q}(t)-n\prod_{i=1}^h \frac{m_q(t_i)}n\right|\ge (1-o(1))K(n).\]
The lemma follows from Lemma \ref{itbazuma} and the union bound.
\end{proof}

With the choice of $K(n)=n^\gamma$ Lemma \ref{interm} implies part \eqref{OTEM2}.

To prove part \eqref{OTEM3}, choose $K(n)=|V|^{-(h-1)} n^\beta$, and observe the following. For  $(q,h)$-tuple $Q$, if we have 
\[\left\|m_{{Q}}-\frac{n}{|V|^h}\mathbbm{1}\right\|_{\infty}<K(n),\]
then $\Sigma(Q)$ is $\beta$-typical.

To prove part \eqref{OTEM4}, we need the following lemma.


\begin{lemma}\label{tech2}
There is a $C_3>0$ such that for every $\beta$-typical $r\in V^n$, if we consider the number of permutations of $r$, i. e., the cardinality of the set $S(r)=\{r'\text{ is a permutation of }r\}$, then we have
\[|S(r)|\ge |V|^n \exp\left(-C_3 n^{2\beta-1}\right).\]  
\end{lemma}
\begin{proof}
This can be proved using Lemma \ref{factorial}.
\end{proof}

Part \eqref{OTEM4} follows from the next lemma.
\begin{lemma}\label{tech3}
We will use the constants $C_1$ and $C_3$ provided by Lemma \ref{tech2} and part \eqref{OTEM2}. For every $\alpha$-typical $q\in V^n$, $\beta$-typical $r\in V^n$ and a random $(q,h)$-tuple $\bar{Q}$, we have
\[\mathbb{P}(\Sigma(\bar{Q})=r\text{ and } \bar{Q} \text{ is not }\gamma\text{-typical})<\frac{C_1 \exp\left(-n^{2\gamma-1}/C_1+C_3 n^{2\beta-1}\right)}{|V|^n}.\]
Here the numerator $C_1 \exp\left(-n^{2\gamma-1}/C_1+C_3 n^{2\beta-1}\right)$ on the right hand side goes to $0$ as $n$ goes to infinity.
\end{lemma}
\begin{proof}
For every $r'\in S(r)$, consider the event that $\Sigma(\bar{Q})=r'$ and $\bar{Q}$ is not \break $\gamma$-typical. These events are disjoint, and by symmetry, they have the same probability. Moreover, they are all contained by the event that $Q$ is not $\gamma$-typical. Thus,
\[\mathbb{P}(\Sigma(\bar{Q})=r\text{ and } \bar{Q} \text{ is not }\gamma\text{-typical})\le\frac{\mathbb{P}(\bar{Q}\text{ is not }\gamma\text{-typical})}{|S(r)|}.\]
The statement then follows from part \eqref{OTEM2} and Lemma \ref{tech2}. 
\end{proof}
This concludes the proof of Lemma \ref{OTEM}.
\end{proof}

Fix an $\alpha$-typical $q\in V^n$. For every $\beta$-typical $r\in R(q,h)$, consider the events \break $(\Sigma(\bar{Q})=r)\wedge (m_{\bar{Q}}\in \mathcal{M}^*(q,r))$. These events are pairwise disjoint. Moreover, from  \eqref{tartalmazza} above, we see that their union contains the event  $(\Sigma(\bar{Q})\text{ is }\beta-\text{typical})\wedge (\bar{Q}\text{ is }\gamma-\text{typical})$ for large enough $n$. So for large enough $n$, we have
\begin{multline}\mathbb{P}((\Sigma(\bar{Q})\text{ is }\beta-\text{typical})\wedge (\bar{Q}\text{ is }\gamma-\text{typical}))\\ \le \sum_{r\in R(q,h)\quad\beta-\text{typical}} \mathbb{P}((\Sigma(\bar{Q})=r)\wedge (m_{\bar{Q}}\in \mathcal{M}^*(q,r)))\le 1.\end{multline}
From part \eqref{OTEM2} and \eqref{OTEM3} of Lemma \ref{OTEM}, we get that \[\mathbb{P}((\Sigma(\bar{Q})\text{ is }\beta-\text{typical})\wedge (\bar{Q}\text{ is }\gamma-\text{typical}))\sim 1\] uniformly for all $\alpha$-typical $q\in V^n$.  
Thus  
\[\sum_{r\in R(q,h)\quad\beta-\text{typical}} \mathbb{P}((\Sigma(\bar{Q})=r)\wedge( m_{\bar{Q}}\in \mathcal{M}^*(q,r)))\sim 1\]
uniformly 
for every $\alpha$-typical $q\in V^n$. Combining this with Equation \eqref{kovetk}, we obtain that 
\begin{multline*}
\mathbb{P}((\Sigma(\bar{Q})=r)\wedge (m_{\bar{Q}}\in \mathcal{M}^*(q,r)))\sim \\
|\{r\in R(q,h)|\quad r\text{ is }\beta\text{-typical}\}|^{-1}\sim |R(q,h)|^{-1}=|V|^{-(n-1)}
\end{multline*}
uniformly for all $\alpha$-typical $q\in V^n$ and $\beta$-typical $r\in R(q,h)$. Here in the second line, we used part \eqref{OTEM1} of Lemma \ref{OTEM}. Finally, using part \eqref{OTEM4} of Lemma \ref{OTEM} and \eqref{tartalmazza}, we get Theorem \ref{mixingThm}.

\subsection{The proof of Theorem \ref{mixingThm2}}

We start by a simple lemma.
\begin{lemma}\label{lemmamixingthm2}
For $q,r\in V^n$, and $h\ge 2$, we have $\mathbb{P}(A_n^{(h)}q=r)\le |S(q)|^{-1}$.
\end{lemma}
\begin{proof}
Let $q'$ be a uniform random permutation of $q$ independent from $A_n^{(h-1)}$. Observe that $A_n^{(h)}q$ has the same distribution as $A_n^{(h-1)}q+q'$.    
The statement of the lemma follows from the facts that
\[\mathbb{P}(A_n^{(h-1)}q+q'=r|\quad r-A_n^{(h-1)}q\sim q)=|S(q)|^{-1}\]
and 
\[\mathbb{P}(A_n^{(h-1)}q+q'=r|\quad r-A_n^{(h-1)}q\not\sim q)=0.\]
\end{proof}

Now we prove Theorem \ref{mixingThm2} from Theorem \ref{mixingThm}.

\begin{proof}
Let $q\in V^n$ be $\alpha$-typical, and let $r\in R(q,d)$. Let $q'$ be a uniform random permutation of $q$ independent from $A_n^{(d-1)}$. Observe that $A_n^{(d)}q$ has the same distribution as $A_n^{(d-1)}q+q'$. Now, we have \[\mathbb{P}(A_n^{(d)}q=r)=\mathbb{E}\mathbb{P}(A_n^{(d-1)}q=r-q'),\] where the expectation is over the  random choice of $q'$. 

Observe that
\begin{itemize}
\item $\mathbb{P}(A_n^{(d-1)}q=r-q')\sim |V|^{-(n-1)}$ uniformly, if $r-q'$ is $\beta$-typical.

\item  $0\le \mathbb{P}(A_n^{(d-1)}q=r-q')\le |S(q)|^{-1}$ otherwise.
\end{itemize}

Indeed, the first statement follows from Theorem \ref{mixingThm} and the fact that \break $r-q'\in R(q,d-1)$. The second statement follows from Lemma \ref{lemmamixingthm2}. 

Moreover, combining Lemma \ref{bazuma} with the union bound, we get the following statement. There is a $c>0$ such that
\[\mathbb{P}(r-q'\text{ is not }\beta-\text{typical})\le \exp(-cn^{2\beta-1}).\]

From the law of total probability, we have
\begin{multline*}\mathbb{P}(A_n^{(d)}q=r)=\mathbb{P}(A_n^{(d-1)}q=r-q'|r-q'\text{ is }\beta-\text{typical})\mathbb{P}(r-q'\text{ is }\beta-\text{typical})\\+ \mathbb{P}(A_n^{(d-1)}q=r-q'|r-q'\text{ is not }\beta-\text{typical})\mathbb{P}(r-q'\text{ is not }\beta-\text{typical}).
\end{multline*}
Inserting the inequalities  above into this, we obtain that 
\begin{multline*}(1+o(1))|V|^{-(n-1)}(1-\exp(-cn^{2\beta-1}))\\\le \mathbb{P}(A_n^{(d)}q=r)\le(1+o(1)) |V|^{-(n-1)}+\frac{\exp(-cn^{2\beta-1})}{|S(q)|}.
\end{multline*}
Since there is $c'$ such that $|S(q)|\ge |V|^n\exp(-c'n^{2\alpha-1})$ for every $\alpha$-typical $q\in V^n$, we get that ${\exp(-cn^{2\beta-1})}/{|S(q)|}=o(|V|^{-n})$. The theorem follows.  
 
\end{proof}

\section{Only the typical vectors matter}\label{Secmom}





The aim of this section to prove Theorem \ref{FoFo}. Let $\Cosets(V)$ be the set of all cosets in $V$. 
Given a function $f(n)$, and a subset $W$ of $V$, a vector $q\in V^n$ will be called $(W,f(n))$-typical if for every $c\in W$, we have $\left|m_q(c)-\frac{n}{|W|}\right|<n^\alpha$ and $\sum_{c\not\in W}m_q(c)\le f(n)$. In the previous section, we used the term $\alpha$-typical for \break $(V,0)$-typical vectors. 

We start by a simple corollary of Theorem \ref{mixingThm2}.
\begin{lemma}\label{tipikuskulonbseg}
We have
\[\lim_{n\to\infty}\sum_{W\in\Cosets(V)} \sum_{\substack{q\text{ is }\\(W,0)-\text{typical}}} d_{\infty}(A_nq,U_{q,d})=0.\]
\end{lemma}
\begin{proof}
If $W$ is a subgroup of $V$, then from Theorem \ref{mixingThm2}, we know that $d_\infty(A_n q,U_{q,d})$ is $o(|W|^{-n})$ uniformly for all $(W,0)$-typical $q$. On the other hand, the number of $(W,0)$-typical vectors is at most $|W|^n$. Thus, \[\lim_{n\to\infty} \sum_{q\text{ is }(W,0)-\text{typical}} d_{\infty}(A_n q,U_{q,d})=0.\]

Consider a coset $W\in\Cosets(V)$ such that $W$ is not a subgroup of $V$. Let $t\in W$, then $W_0=W-t$ is a subgroup of $V$. For $q=(q_1,q_2,\dots,q_n)\in W^n$, we define $q'=(q_1-t,q_2-t,\dots,q_n-t)$. Note that $q\mapsto q'$ is a bijection between $W^n$ and $W_0^n$, and it is also a bijection between $(W,0)$-typical and $(W_0,0)$-typical vectors.  Using this, it is easy to see that $d_\infty(A_nq,U_{q,d})=d_\infty(A_n q',U_{q',d})$, which implies that \[\lim_{n\to\infty}\sum_{q\text{ is }(W,0)-\text{typical}} d_{\infty}(A_n q,U_{q,d})=\lim_{n\to\infty}\sum_{q'\text{ is }(W_0,0)-\text{typical}} d_{\infty}(A_n q',U_{q',d})=0,\]
using the already established case. Since  $\Cosets(V)$ is finite, the statement follows.
 \end{proof}

For $q\in V^n$, choose $r_q$ such that \[\mathbb{P}(A_n q=r_q)=\max_{r\in V^n} \mathbb{P}(A_n q=r).\] 
For $W\in\Cosets(V)$, we define $I(W^n)=\{q\in W^n\quad|\quad \Aff_q=W\}$. Note that $V^n=\cup_{W\in\Cosets(V)} I(W^n)$, where this is a disjoint union. 

Then
\begin{align}
\limsup_{n\to\infty} & \sum_{q\in V^n} d_\infty(A_n q,U_{q,d})\nonumber\\
&=\limsup_{n\to\infty} \sum_{W\in\Cosets(V)} \sum_{q\in I(W^n)} d_\infty(A_n q,U_{q,d})\nonumber\\
&=\limsup_{n\to\infty} \sum_{W\in\Cosets(V)} \sum_{\substack{q\text{ is }\\(W,0)-\text{typical}}} d_\infty(A_n q,U_{q,d})\nonumber\\ 
&\qquad+\limsup_{n\to\infty} \sum_{W\in\Cosets(V)}\quad  \sum_{\substack{q\in I(W^n)\text{ is }\\\text{not }(W,0)-\text{typical}}} d_\infty(A_n q,U_{q,d}).\label{felbecsles} 
\end{align}
Using Lemma \ref{tipikuskulonbseg}, we have
\[\limsup_{n\to\infty} \sum_{W\in\Cosets(V)} \sum_{\substack{q\text{ is }\\(W,0)-\text{typical}}} d_\infty(A_n q,U_{q,d})=0.\]
For $q\in I(W^n)$, we have \[d_\infty(A_n q,U_{q,d})\le |W|^{-(n-1)}+\mathbb{P}(A_n q=r_q)\] from the triangle inequality. Moreover, \[|\{q\in I(W^n)\quad|\quad q\text{ is not }(W,0)-\text{typical}\}|=o(|W|^n)\]
from standard concentration results.

Inserting these into  Equation \eqref{felbecsles}, we obtain that 
\begin{align*}
\limsup_{n\to\infty} & \sum_{q\in V^n} d_\infty(A_n q,U_{q,d})\\
&\le\limsup_{n\to\infty} \sum_{W\in\Cosets(V)} \quad\sum_{\substack{q\in I(W^n)\text{ is}\\\text{not }(W,0)-\text{typical}}} \left(|W|^{-(n-1)}+\mathbb{P}(A_n q=r_q)\right)\\
&=\limsup_{n\to\infty} \sum_{W\in\Cosets(V)} |\{q\in I(W^n)\quad|\quad q\text{ is not }(W,0)-\text{typical}\}||W|^{-(n-1)}\\
&\qquad+ \limsup_{n\to\infty} \sum_{W\in\Cosets(V)} \quad \sum_{\substack{q\in I(W^n)\text{ is}\\\text{not }(W,0)-\text{typical}}} \mathbb{P}(A_n q=r_q)\\
&= \limsup_{n\to\infty} \sum_{W\in\Cosets(V)}\quad \sum_{\substack{q\in I(W^n)\text{ is}\\ \text{not }(W,0)-\text{typical}}} \mathbb{P}(A_n q=r_q).
\end{align*}

Thus, in order to prove Theorem \ref{FoFo}, it is enough to prove that \[\limsup_{n\to\infty} \sum_{W\in\Cosets(V)} \quad\sum_{\substack{q\in I(W^n)\text{ is}\\ \text{not }(W,0)-\text{typical}}} \mathbb{P}(A_n q=r_q)=0.\]

We establish this in three steps, namely, we prove that

\begin{align}
\label{eegyenlo1}
\limsup_{n\to\infty}  \quad\sum_{\substack{q\in V^n\text{ is not}\\ (W,n^{\alpha})-\text{typical for any }W\in\Cosets(V)}} \mathbb{P}(A_n q=r_q)&=0,\\
\label{eegyenlo2}
\limsup_{n\to\infty} \sum_{W\in\Cosets(V)} \quad\sum_{\substack{q\text{ is  }(W,n^{\alpha})-\text{typical,}\\ \text{but not } (W,C\log n)-\text{typical}}} \mathbb{P}(A_n q=r_q)&=0,\\
\label{eegyenlo3}
\limsup_{n\to\infty} \sum_{W\in\Cosets(V)} \quad\sum_{\substack{q\text{ is  }(W,C\log n)-\text{typical,}\\ \text{but not } (W,0)-\text{typical}}} \mathbb{P}(A_n q=r_q)&=0,
\end{align}
where $C$ is a constant to be chosen later.

Equations \eqref{eegyenlo1}, \eqref{eegyenlo2} and \eqref{eegyenlo3} are proved in Subsections \ref{egyenlo1}, \ref{egyenlo2} and \ref{egyenlo3} respectively.

\subsection {Proof of Equation \eqref{eegyenlo1}}\label{egyenlo1}

The following terminology will be useful for us. With every $(q,d-1)$-tuple $Q=(Q_1,Q_2,\dots,Q_n)$ we associate the random variables $Z\in V$ and $X^Q=(X^Q_1,X^Q_2,\dots,X^Q_{d-1})\in V^{d-1}$, such that $Z=r_q(i)$ and $X^Q=Q_i$, where $i$ is a uniform random element of the set $\{1,2,\dots,n\}$. 
Each $X^Q_j$ has the same distribution as $q_i$ where $i$ is chosen uniformly from $\{1,2\dots,n\}$.  The random variable $X^Q_\Sigma\in V$ is defined as $X^Q_\Sigma=\sum_{i=1}^{d-1} X^Q_i$. 
These two sets of $(q, d - 1)$-tuples are equal:
\[\{Q\quad|\quad r_q-\Sigma(Q)\sim q\}=\{Q\quad|\quad Z-X_\Sigma^Q\eqd X^Q_1\}.\]
Here $\eqd $ means that the two random variables have the same distribution. Thus, \[\mathbb{P}\left(r_q -A_n^{(d-1)}q\sim q\right)=\mathbb{P}_{\bar{Q}}\left(Z-X^{\bar{Q}}_\Sigma \eqd X^{\bar{Q}}_1\right),\] 
where the subscript in the notation $\mathbb{P}_{\bar{Q}}$ indicates that the probability is over the random choice of $\bar{Q}$. 

We call the random variables $Z,X_1,X_2,...,X_{d-1}\in V$ $\varepsilon$-independent, if for every $z,x_1,x_2,...,x_{d-1}\in V$, we have
\begin{multline*}
|\mathbb{P}(Z=z,X_1=x_1,...,X_{d-1}=x_{d-1})-\mathbb{P}(Z=z)\mathbb{P}(X_1=x_1)\cdots \mathbb{P}(X_{d-1}=x_{d-1})|\\<\varepsilon.
\end{multline*}

Fix $\frac{1}{2}<\eta<\alpha$. The next lemma follows from Lemma \ref{itbazuma} and the union bound. 
\begin{lemma}\label{tech26}
For any $q\in V^n$, we have 
\begin{multline*}\mathbb{P}_{\bar{Q}}(Z,X^{\bar{Q}}_1,X^{\bar{Q}}_2,\dots,X^{\bar{Q}}_{d-1}\text{ are not }n^{\eta-1}\text{-independent})\\\le |V|^d 2(d-1)\exp\left(-\frac{2 n^{2\eta-1}}{(d-1)^2}\right).\qed
\end{multline*}
\end{lemma}


The crucial step in the proof of Equation \eqref{eegyenlo1} is the following lemma, which is proved in the next subsection.

\begin{lemma}\label{lenyeges}
Let  $d\ge 3$. There is $C$ and $\varepsilon_0>0$ (which may depend on $d$ and $V$), such that the following holds. Assume that $Z,X_1,X_2,...,X_{d-1}$ are $\varepsilon$-independent $V$-valued random variables, for some $0<\varepsilon<\varepsilon_0$. Let $X_{\Sigma}=X_1+X_2+\dots +X_{d-1}$. Assume that $X_1,X_2,\dots,X_{d-1}$ and $Z-X_\Sigma$ have the same distribution $\pi$. Then there is a coset $W$ in $V$ such that $d_{\infty}(\pi,\pi_W)<C\varepsilon$.
\end{lemma}

Here $\pi_W$ is the uniform distribution on $W$. For two distribution $\pi$ and $\mu$ on the same finite set $\mathcal{R}$, their distance $d_{\infty}(\pi,\mu)$ is defined as \[d_{\infty}(\pi,\mu)=\max_{r\in\mathcal{R}} |\pi(r)-\mu(r)|.\]  

Combining the last lemma with Lemma \ref{tech26}, we get the following lemma.
\begin{lemma}\label{lemma13nak} Assume that $n$ is large enough. Let $q\in V^n$. If 
\[\mathbb{P}\left(r_q -A_n^{(d-1)}q\sim q\right)=\mathbb{P}_{\bar{Q}}\left(Z-X^{\bar{Q}}_\Sigma \eqd X^{\bar{Q}}_1\right)> |V|^d 2(d-1)\exp\left(-\frac{2 n^{2\eta-1}}{(d-1)^2}\right),\]
then $q$ is $(W,n^{\alpha})$-typical for some coset $W$ in $V$. In other words, if $q$ is not \break $(W,n^{\alpha})$-typical for any coset $W$, then  
\[\mathbb{P}\left(r_q -A_n^{(d-1)}q\sim q\right)=\mathbb{P}_{\bar{Q}}\left(Z-X^{\bar{Q}}_\Sigma \eqd X^{\bar{Q}}_1\right)\le |V|^d 2(d-1)\exp\left(-\frac{2 n^{2\eta-1}}{(d-1)^2}\right).\]
\end{lemma}
\begin{proof}

Combining our assumptions on $q$ with Lemma \ref{tech26}, we have
\[\mathbb{P}_{\bar{Q}}\left(Z,X^{\bar{Q}}_1,X^{\bar{Q}}_2,\dots,X^{\bar{Q}}_{d-1}\text{ are }n^{\eta-1}\text{-independent and }Z-X^{\bar{Q}}_\Sigma \eqd X^{\bar{Q}}_1\right)>0.\]
 So there exist $n^{\eta-1}$-independent random variables $Z,X_1,X_2,\dots,X_{d-1}$, such that $X_1,X_2,\dots,X_{d-1}$ and $Z-X_\Sigma=Z-\sum_{i=1}^{d-1} X_i$ all have the same distribution as $q_i$ where $i$ is chosen uniformly from $\{1,2,...,n\}$. Let us call this distribution $\pi$. For large enough $n$, we have $n^{\eta-1}<\varepsilon_0$, so Lemma \ref{lenyeges} can be applied to give us that  there is a coset $W$ in $V$ such that $d_{\infty}(\pi,\pi_W)<Cn^{\eta-1}$. Since $n^\alpha>C |V|n^\eta$, this implies that $q$ is $(W,n^{\alpha})$-typical. 
\end{proof}

Now we made all the necessary preparations to prove Equation \eqref{eegyenlo1}.

Due to symmetry if $q_1\sim q_2$, then $\mathbb{P}(A_n^{(d)} q_1=r_{q_1})=\mathbb{P}(A_n^{(d)} q_2=r_{q_2})$. Let $q^{(d)}$ be a uniform random permutation of $q$ independent from $A_n^{(d-1)}$.  

We have
\begin{align*}
\sum_{q'\sim q} \mathbb{P}(A_n^{(d)} q'=r_{q'})  &=|S(q)| \mathbb{P}(A_n^{(d)}q=r_q)\\
&=|S(q)| \mathbb{P}(A_n^{(d-1)}q+q^{(d)}=r_q)\\
&=|S(q)|\sum_{q'\sim q} \mathbb{P}(A_n^{(d-1)}q=r_q-q')\mathbb{P}(q^{(d)}=q')\\
&=\sum_{q'\sim q} \mathbb{P}(A_n^{(d-1)}q=r_q-q')=\mathbb{P}(r_q-A_n^{(d-1)}q\sim q).
\end{align*}
Let $T_n\subset V^n$ be such that it contains exactly one element of each equivalence class. Then, assuming that $n$ is large enough, we have

\begin{align*}
\sum_{\substack{q\in V^n\text{ is not}\\ (W,n^{\alpha})-\text{typical for any }W\in\Cosets(V)}} &\mathbb{P}(A^{(d)}_n q=r_q)\\&=\sum_{\substack{q\in T_n\text{ is not}\\ (W,n^{\alpha})-\text{typical for any }W\in\Cosets(V)}}\mathbb{P}(r_q-A_n^{(d-1)}q\sim q)  \\&\le |T_n| |V|^d 2(d-1)\exp\left(-\frac{2 n^{2\eta-1}}{(d-1)^2}\right).
\end{align*} 

In the last step, we used Lemma \ref{lemma13nak}. Equation \eqref{eegyenlo1} follows from the fact that $|T_n|=o\left(n^{|V|+1}\right)=o\left(\exp\left(\frac{2 n^{2\eta-1}}{(d-1)^2}\right)\right)$.


\subsection{The proof of Lemma \ref{lenyeges}}

Although we will not use the following lemma directly, we include it and its proof, because it contains many ideas, that will occur later, in a much clearer form. 
\begin{lemma}\label{szemlelet}
Let $Z,X_1,X_2,...,X_{d-1}$ be  independent $V$-valued random variables. Let $X_{\Sigma}=X_1+X_2+\dots+X_{d-1}$. Assume that $X_1,X_2,\dots,X_{d-1}$ and $Z-X_\Sigma$ have the same distribution $\pi$. Then $\pi=\pi_W$ for some coset $W$ in $V$.
\end{lemma}

\begin{proof}
We use discrete Fourier transform, that is, for $\varrho\in \hat{V}=\Hom(V,\mathbb{C}^*)$, we define
\[\hat{\pi}(\varrho)=\sum_{v\in V} \pi(v)\varrho(v)\]
and 
\[\hat{\mu}(\varrho)=\sum_{v\in V} \mathbb{P}(Z=v)\varrho(v).\]
The assumptions of the lemma imply that 
\[\hat{\mu}(\varrho)\left(\overline{\hat{\pi}(\varrho)}\right)^{d-1}=\hat{\pi}(\varrho)\]
for every $\varrho\in \hat{V}$. In particular $|\hat{\mu}(\varrho)|\cdot\left|\hat{\pi}(\varrho)\right|^{d-1}=|\hat{\pi}(\varrho)|$ for every $\varrho\in \hat{V}$. Since $|\hat{\mu}(\varrho)|, |\hat{\pi}(\varrho)|\le 1$, this is only possible if $|\hat{\pi}(\varrho)|\in \{0,1\}$ for every $\varrho\in \hat{V}$. Let us define $\hat{V}_1=\{\varrho\in \hat{V}|\quad |\hat{\pi}(\varrho)|=1\}$. Note that $\hat{V}_1$ always contains the trivial character. Then for every $\varrho\in \hat{V}_1$, the character $\varrho$ is constant on the support of $\pi$. Or in other words, the support of $\pi$ is contained in $W_\varrho=\varrho^{-1}(\hat{\pi}(\varrho))$, which is a coset of $\ker \varrho$. Therefore, the support of $\pi$ is contained in the coset $W=\cap_{\varrho\in \hat{V}_1}W_\varrho$.  Now we prove that $\hat{\pi}(\varrho)=\hat{\pi}_W(\varrho)$ for every $\varrho\in \hat{V}$, which implies that $\pi=\pi_W$. This is clear for $\varrho\in \hat{V}_1$, so assume that $\varrho\not\in \hat{V_1}$, that is, $\hat{\pi}(\varrho)=0$. This implies that $\varrho$ is not constant on $W$. So there are $w_1,w_2\in W$ such that $\varrho(w_1)\neq\varrho(w_2)$. For $w=w_1-w_2$, we have $\varrho(w)\neq 1$ and $W=w+W$. Thus
\begin{align}\label{zerohatpi}
\hat{\pi}_W(\varrho)&=\frac{1}{|W|}\sum_{v\in W} \varrho(v)=\frac{1}{|W|}\sum_{v\in W} \varrho(w+v)\\&=\frac{1}{|W|}\varrho(w)\sum_{v\in W} \varrho(v)=\varrho(w)\hat{\pi}_W(\varrho).\nonumber
\end{align}
Since $\varrho(w)\neq 1$, this means that $\hat{\pi}_W(\varrho)=0$.
\end{proof}

Now we turn to the proof of Lemma \ref{lenyeges}.
\begin{proof}
Using the notations of the proof of Lemma \ref{szemlelet}, the conditions of the lemma imply that
\[\left|\hat{\pi}(\varrho)-\hat{\mu}(\varrho)\left(\overline{\hat{\pi}(\varrho)}\right)^{d-1}\right|\le |V|^d \varepsilon\]
for every $\varrho\in \hat{V}$. Using the fact that $|\hat{\mu}(\varrho)|\le 1$, we obtain

\[\left|\hat{\pi}(\varrho)-\hat{\mu}(\varrho)\left(\overline{\hat{\pi}(\varrho)}\right)^{d-1}\right|\ge |\hat{\pi}(\varrho)|-\left|\hat{\mu}(\varrho)\right|\cdot\left|\hat{\pi}(\varrho)\right|^{d-1} \ge |\hat{\pi}(\varrho)|-\left|\hat{\pi}(\varrho)\right|^{d-1} ,\]  
which gives us $|\hat{\pi}(\varrho)|-\left|\hat{\pi}(\varrho)\right|^{d-1}\le |V|^d \varepsilon$ for every $\varrho\in \hat{V}$. 


Consider the $[0,1]\to [0,1]$ function $x\mapsto x-x^{d-1}$, this function only vanishes at $0$ and $1$. Moreover, the derivative of this function does not vanish at $0$ and $1$. This implies that there is an $\varepsilon_1>0$ and a $C_1>0$ such that for every $0<\varepsilon<\varepsilon_1$ the following holds. For $x\in [0,1]$, if we have $x-x^{d-1} \le |V|^{d} \varepsilon$, then either $x<C_1\varepsilon$ or $x>1-C_1\varepsilon$. In the rest of the proof, we assume that $\varepsilon<\varepsilon_1$. Then for every $\varrho\in \hat{V}$, we have either $|\hat{\pi}(\varrho)|<C_1\varepsilon$ or $|\hat{\pi}(\varrho)|>1-C_1\varepsilon$. 

Let $\hat{V}_1=\{\varrho\in \hat{V}| 1-C_1\varepsilon<|\hat{\pi}(\varrho)|\}$. Take any $\varrho\in \hat{V}_1$. Set \[z=\frac{\overline{\hat{\pi}(\varrho)}}{|\hat{\pi}(\varrho)|}.\] 

Choose $\xi_0=\xi_0(\varrho)$ in the range $R(\varrho)$ of the character $\rho$, such that \break $\Rep z\xi_0=\max_{\xi\in R(\varrho)} \Rep z \xi$. An elementary geometric argument gives that for \break $\xi_0	\neq \xi\in R(\varrho)$, we have $\Rep z\xi\le 1-\delta$, where $\delta=1-\cos\frac{\pi}{|V|}>0$.\footnote{Here $\pi=3.14\dots$ is the well-known constant.}  Clearly $\Rep z\xi_0 \le 1$. Then we have
\[|\hat{\pi}(\varrho)|=z\hat{\pi}(\varrho)=\Rep z \hat{\pi}(\varrho)=\sum_{\xi\in R(\varrho)} \pi(\varrho^{-1}(\xi)) \Rep z\xi\le 1-\left(1-\pi(\varrho^{-1}(\xi_0))\right)\delta.\]
Thus, $|\hat{\pi}(\varrho)|>1-C_1\varepsilon$ implies that  for the coset $W_\varrho=\varrho^{-1}(\xi_0)$, we have \break $\pi(W_\varrho)>1-C_1\delta^{-1}\varepsilon$. So the coset $W=\cap_{\varrho\in \hat{V}_1} W_\varrho$ satisfies  $\pi(W)>1-C_1\delta^{-1}|V|\varepsilon$. 

Consider a $\varrho\in \hat{V}_1$. Let $\xi_0=\xi_0(\varrho)$ be like above. Note that $\varrho(v)=\xi_0$ for any $v\in W_\varrho$. In particular, we have $\hat{\pi}_W(\varrho)=\xi_0$. Thus,
\begin{align*}
|\hat{\pi}_W(\varrho)-\hat{\pi}(\varrho)|&=\left|\xi_0-\left(\pi(W_\varrho)\xi_0-\sum_{v\in V\backslash W_\varrho}\pi(v)\varrho(v)\right)\right|\\&=\left|(1-\pi(W_\varrho))\xi_0-\sum_{v\in V\backslash W_\varrho}\pi(v)\varrho(v)\right|\\
&\le 1-\pi(W_\varrho)+\sum_{v\in V\backslash W_\varrho}\pi(v)=  2 (1-\pi(W_\varrho))\le 2C_1\delta^{-1}\varepsilon.
\end{align*}
Now take $\varrho\in \hat{V}\backslash \hat{V}_1$. 
We know that $|\hat{\pi}(\varrho)|<C_1\varepsilon$. We claim that $\varrho$ is not constant on $W$. To show this, assume that $\varrho$ is constant on $W$, then 
\[|\hat{\pi}(\varrho)|\ge \pi(W)-\pi(V\backslash W)\ge 1-2C_1\delta^{-1}|V|\varepsilon>C_1\varepsilon\]
provided that $\varepsilon$ is small enough, which gives us a contradiction. Using that $\varrho$ is not constant on $W$, Equation \eqref{zerohatpi} gives us $\hat{\pi}_W(\varrho)=0$. Thus,
\[|\hat{\pi}(\varrho)-\hat{\pi}_W(\varrho)|=|\hat{\pi}(\varrho)|\le C_1\varepsilon.\]

This gives us that $|\hat{\pi}(\varrho)-\hat{\pi}_W(\varrho)| \le 2C_1\delta^{-1}\varepsilon$ for any $\varrho\in \hat{V}$. Since the map \break $\pi\mapsto \hat{\pi}$   is an invertible linear map,  there is a constant $L=L_V$ such that \break $d_\infty(\pi,\pi_W)\le L\max_{\varrho\in \hat{V}} |\hat{\pi}(\varrho)-\hat{\pi}_W(\varrho)|$. This gives the statement.
\end{proof}

\subsection {Proof of Equation \eqref{eegyenlo2}}\label{egyenlo2}
We start by the following lemma.

\begin{lemma}\label{ee2lemma}
There is a $C$ such that if $W\in\Cosets(V)$ and $q\in V^n$ is $(W,n^{\alpha})$-typical, but not $(W,C\log n)$-typical, then for a random $(q,d-1)$-tuple $\bar{Q}$, we have \[\mathbb{P}(r_q-\Sigma(\bar{Q})\sim q)\le n^{-(|V|+1)}.\]
\end{lemma}
\begin{proof}
Let $E=\sum_{c\not\in W} m_q(c)$. Since $q$ is $(W,n^\alpha)$-typical, we have $E\le n^{\alpha}$. Assume that $r=\sum_{i=1}^d q^{(i)}$, where $q^{(i)}\sim q$. Note that \[\{j|\quad r_j\not\in dW\}\subset \cup_{i=1}^d \{j|\quad q^{(i)}(j)\not\in W\},\] so $\sum_{c\not\in dW} m_r(c)\le d E$. In particular, this is true for $r_q$, that is,  \[\sum_{c\not\in dW} m_{r_q}(c)\le d E.\]

Let \[H_0=\{j\quad|\quad r_q(j)\not\in dW\}.\] For $i=1,2,...,d-1$, we define the random subset $H_i$ of $\{1,2,...,n\}$ using the random $(q,d-1)$-tuple $\bar {Q}=(\bar {q}^{(1)},\bar {q}^{(2)},\dots,\bar {q}^{(d-1)})$ as 
\[H_i=\{j\quad|\quad \bar{q}^{(i)}(j)\not\in W\},\] and let the random subset $H^*\subset\{1,2,\dots,n\}$ be defined as 
\[H^*=\{j\quad| r_q(j)-\Sigma(\bar{Q})(j)\not\in W\}.\] Then $0\le |H_0|\le dE$ and  $|H_1|=|H_2|=...=|H_{d-1}|=E$. Let 
\[B=\{j\quad|\quad j\text{ is contained in exactly one of the sets $H_0,H_1,H_2,...,H_{d-1}$}\}.\]

 Then $B\subset H^*$, therefore we have
\[\mathbb{P}(r_q-\Sigma(\bar{Q})\sim q)\le \mathbb{P}(|H^*|=E)\le \mathbb{P}(|B|\le E).\]
We will need the following inequality
\[|B|\ge \sum_{i=0}^{d-1} |H_i|-2\sum_{0\le i<j\le d-1} |H_i\cap H_j|\ge(d-1)E-2\sum_{0\le i<j\le d-1} |H_i\cap H_j|.\]

The proof of this is straightforward, or see \cite[Chapter IV, 5.(c)]{feller}. Thus, if $|B|\le E$, then \[2\sum_{0\le i<j\le d-1} |H_i\cap H_j|\ge (d-2) E.\] So $|H_i\cap H_j|\ge \frac{(d-2)E}{d(d-1)}$ for some  $i<j$. Therefore,
\begin{equation}\label{szita}
\mathbb{P}(r_q-\Sigma(Q)\sim q)\le \mathbb{P}(|B|\le E)\le \sum_{0\le i<j\le d-1} \mathbb{P}\left(|H_i\cap H_j|\ge \frac{(d-2)E}{d(d-1)}\right).
\end{equation}
\begin{lemma}\label{most}
There is a constant $C$ such that, for all $a,b$ and $E$ satisfying  \break $C\log n< E< n^{\alpha}$ and $a,b\le dE$, if $A$ and $B$ are two random subset of $\{1,2,...,n\}$ of size $a$ and $b$ respectively chosen independently and uniformly, then
\[\mathbb{P}\left(|A\cap B|\ge \frac{(d-2)E}{d(d-1)}\right)<n^{-(|V|+1)}\Big/{{d}\choose{2}}.\]
\end{lemma}
\begin{proof}
We may assume that $n$ is large enough, because we can always increase $C$ to handle the small values of $n$. Let $\delta=\frac{(d-2)}{d(d-1)}$. For large enough $n$, we have $\frac{ab}{n}\le \frac{\delta}{2}E$. Using Lemma \ref{bazuma}, we obtain that
\begin{align*}
\mathbb{P}\left(|A\cap B|\ge \frac{(d-2)E}{d(d-1)}\right)&=\mathbb{P}\left(|A\cap B|\ge \delta E\right)\\
&\le \mathbb{P}\left(\left||A\cap B|-\frac{ab}{n}\right|\ge \frac{\delta}{2} E\right)\le 2\exp\left(-\frac{\delta^2E^2}{2a}\right)\\&\le 2\exp\left(-\frac{\delta^2E}{2d}\right)\le 2\exp\left(-\frac{\delta^2C\log n}{2d}\right)\\&=2n^{-\frac{\delta^2C}{2d}}<n^{-(|V|+1)}\Big/{{d}\choose{2}}  
\end{align*}
for large enough $C$.
\end{proof}
Combining this lemma with Inequality \eqref{szita}, we get the statement of Lemma~\ref{ee2lemma}. 
\end{proof}
Then Equation \eqref{eegyenlo2} follows, because
\begin{align*}
\limsup_{n\to\infty} \sum_{W\in\Cosets(V)}& \quad\sum_{\substack{q\text{ is  }(W,n^{\alpha})-\text{typical,}\\ \text{but not } (W,C\log n)-\text{typical}}} \mathbb{P}(A_n^{(d)} q=r_q)\\&=
\limsup_{n\to\infty} \sum_{W\in\Cosets(V)} \quad\sum_{\substack{q\in T_n\text{ is  }(W,n^{\alpha})-\text{typical,}\\ \text{but not } (W,C\log n)-\text{typical}}} \mathbb{P}(r_q-A_n^{(d-1)}q\sim q)\\&\le \limsup_{n\to\infty} |\Cosets(V)|\cdot|T_n| n^{-(|V|+1)}=0.
\end{align*}

\subsection {Proof of Equation \eqref{eegyenlo3}}\label{egyenlo3}

Since there are only finitely many cosets in $V$, it is enough to prove that for any coset $W\in \Cosets (V)$, we have
\[\lim_{n\to\infty} \sum_{q\in D_W^n} |S(q)|\mathbb{P}(\Sigma(\bar{Q})=r_q)=0,\]
where \[D_W^n=\{q\in T_n\quad|\quad q\text{ is }(W,C\log n)-\text{typical, but not }(W,0)\text{-typical}\},\] and $\bar{Q}$ is a random $(q,d)$-tuple. (Recall that $S(q)$ is the set of permutations of $q$.) 


Given a $q\in V^n$, a $(q,d)$-tuple $Q$ or $m_Q$ itself will be called $W$-decent if for any $u\in W^d$, we have \[\frac{1+m_{\Sigma(Q)}(u_\Sigma)}{1+m_Q(u)}\le \log^2 n,\] 
and it will be called $W$-half-decent if $(1+m_{\Sigma(Q)}(u_\Sigma))/(1+m_Q(u))\le \log^4 n$. 
Or even more generally, a non-negative integral vector $m$ indexed by $V^d$ will be called $W$-half-decent if for every $u\in W^d$, we have \[\frac{1+ m(\tau_\Sigma=u_\Sigma)}{1+m(u)}\le \log^4 n,\] where $n=\sum_{t\in V^d} m(t)$.

\begin{lemma}\label{csakdecent}
For any coset $W\in\Cosets(V)$, we have
\begin{multline*}
\limsup_{n\to\infty} \sum_{q\in D_W^n} |S(q)|\mathbb{P}(\Sigma(\bar{Q})=r_q)\\=\limsup_{n\to\infty} \sum_{q\in D_W^n} |S(q)|\mathbb{P}(\Sigma(\bar{Q})=r_q\text{ and }\bar{Q}\text{ is }W-\text{decent}).
\end{multline*}
\end{lemma}
\begin{proof}

It is enough to show that if $n$ is large enough, then \[|S(q)|\mathbb{P}(\Sigma(\bar{Q})=r_q\text{ and }\bar{Q}\text{ is not  }W-\text{decent})\le n^{-(|V|+1)}\] for every  $q\in D_W^n$.  
Indeed, once we establish this, it follows that
\begin{multline*}\limsup_{n\to\infty} \sum_{q\in D_W^n} |S(q)|\mathbb{P}(\Sigma(\bar{Q})=r_q\text{ and }\bar{Q}\text{ is not }W-\text{decent})\\\le \limsup_{n\to\infty}|T_n|n^{-(|V|+1)}=0,
\end{multline*}
which gives the statement.

Just for this proof $(q,h)$-tuples and random $(q,h)$-tuples will be denoted by $Q^h$ and $\bar{Q}^h$, because it will be important to emphasize the value of $h$.
 Given any $(q,d-1)$-tuple $Q^{d-1}=(q^{(1)},q^{(2)},\dots,q^{(d-1)})$ such that $r_q-\Sigma(Q^{d-1})\sim q$ the tuple $(q^{(1)},q^{(2)},\dots,q^{(d-1)},r_q-\Sigma(Q^{d-1}))$ will be a $(q,d)$-tuple and it is denoted by $\text{Ext}(Q^{d-1})$. 
It is also clear that $\Sigma(\text{Ext}(Q^{d-1}))=r_q$, and for any $(q,d)$-tuple $Q^d$ such that $\Sigma(Q^d)=r_q$ there is a unique $(q,d-1)$-tuple $Q^{d-1}$ such that $r_q-\Sigma(Q^{d-1})\sim q$ and $Q^d=\text{Ext}(Q^{d-1})$.  Also note that $\mathbb{P}(\bar{Q}^{d-1}=Q^{d-1})=|S(q)|\mathbb{P}(\bar{Q}^d=Q^d)$. 

Therefore, for any $q\in D_W^n$, we have
\begin{multline*}
|S(q)|\mathbb{P}(\Sigma(\bar{Q}^{d})=r_q\text{ and }\bar{Q}\text{ is not }W-\text{decent})\\=\mathbb{P}(r_q-\Sigma(\bar{Q}^{d-1})\sim q\text{ and }\text{Ext}(\bar{Q}^{d-1})\text{ is not }W-\text{decent}).
\end{multline*}
The event on the right-hand side is contained in the even that
\[
\text{there are }t\in W^{d-1}\text{ and }c\in dW\text{, such that }\qquad\qquad\qquad\qquad\]
\begin{equation}\label{ineqdec}
\frac{1+m_{r_q}(c)}{1+|\{i|\quad r_q(i)=c\text{ and }\bar{Q}^{d-1}(i)=t\}|}> \log^2 n.   
\end{equation}
This event has probability at most $n^{-(|V|+1)}$ for every $(W,C\log n)$-typical vector $q\in V^n$, if $n$ is large enough. Indeed, for a $c\in dW$ such that  $m_{r_q}(c)< \log^2 n$, Inequality \eqref{ineqdec}  can  not be true. On the other hand, if $m_{r_q}(c)\ge \log^2 n$, then with high probability \[|\{i|\quad r_q(i)=c\text{ and }\bar{Q}^{d-1}(i)=t\}|>\frac{1}{2}\frac{m_{r_q}(c)}{|W|^{d-1}}>\frac{1+m_{r_q}(c)}{\log^2 n}\] for any $t\in W^{d-1}$, as it follows from Lemma \ref{itbazuma}.
\end{proof}

As before, we define 
\[\mathcal{M}(q,r)=\{m_Q\quad|\quad Q\in\mathcal{Q}_{q,d},\Sigma(Q)=r\}.\]
Let
\[\mathcal{M}^\sharp(q,r)=\{m\in \mathcal{M}(q,r)|\quad m\text{ is }W-\text{decent}\}.\]
From the previous lemma, we need to prove that
\[\lim_{n\to\infty}\sum_{q\in D_W^n}\quad\sum_{m\in \mathcal{M}^\sharp(q,r_q)} |S(q)|\mathbb{P}((\Sigma(\bar{Q})=r_q)\wedge (m_{\bar{Q}}=m ))=0.\]

Let \[\mathcal{M}=\{m_Q\quad|\quad Q\text{ is a }(q,d)\text{-tuple for some }n\ge 0\text{ and }q\in V^n\}.\]

The set $\mathcal{M}$ is the set of non-negative integral points of the linear subspace of $\mathbb{R}^{V^{d}}$ consisting of the vectors $m$ satisfying the following linear equations:

\[ m(\tau_i=c)= m(\tau_1=c)\]
for every $c\in V$ and $i=1,2,\dots,d$.


In other words, $\mathcal{M}$ consists of the integral points of a rational polyhedral cone. From  \cite[Theorem 16.4]{schr}, we know that this cone is generated by an integral Hilbert basis, i. e., we have the following lemma. 
\begin{lemma}\label{Hilbert}
There are finitely many vectors $m_1,m_2,...,m_{\ell}\in \mathcal{M}$, such that
\[
\pushQED{\qed}\mathcal{M}=\{c_1 m_1+c_2m_2+\dots+ c_\ell m_\ell \quad |\quad c_1,c_2,\dots, c_{\ell}\text{ are non-negative integers}\}.\qedhere
\popQED
\]
\end{lemma}

We may assume that the indices in the lemma above are chosen such that there is an $h$ such that the supports of $m_1,m_2,\dots,m_h$ are contained in $W^{d}$, and the supports of $m_{h+1},m_{h+2},...,m_\ell$ are not contained in $W^{d}$. 
\begin{definition}
Given a vector $m\in \mathcal{M}$, write $m$ as $m=\sum_{i=1}^\ell c_im_i$, where $c_1,c_2,...,c_\ell$ are non-negative integers, and let $\Delta(m)=\sum_{i=h+1}^\ell c_im_i$. (If the decomposition of $m$ is not unique just pick and fix a decomposition.)
\end{definition}

With the notation  $\|m\|_{W^C}=m(\tau\not\in W^d)$, 
we have $\|m\|_{W^C}=\|\Delta(m)\|_{W^C}$ and $\|m-\Delta(m)\|_{W^C}=0$.

For any non-negative integral vector $m\in \mathbb{R}^{V^{d}}$, we define

\begin{equation}\label{Emdef}
E(m)=\frac{\prod_{c\in V}  m(\tau_\Sigma=c)!}{\prod_{t\in V^d} m(t)!}
\left(\prod_{i=1}^d \frac{\prod_{c\in V} m(\tau_i=c)!}{m(V^d)!}\right)^{\frac{d-1}{d}}.
\end{equation}

\begin{lemma}
For every $q,r\in V^n$ and $m\in\mathcal{M}(q,r)$, we have
\[|S(q)|\mathbb{P}((\Sigma(\bar{Q})=r)\wedge(m_{\bar{Q}}=m))=\frac{\prod_{c\in V} m_r(c)!}{\prod_{t\in V^{d}} m(t)!}\Big/\left(\frac{n!}{\prod_{c\in V} m_q(c)!}\right)^{d-1}=E(m).\]
\end{lemma}
\begin{proof} 
The first equality is a consequence of Lemma \ref{euq12lemma}. To prove the second equality, note that since $m\in \mathcal{M}(q,r)$, for any $c\in V$ and $i\in \{1,2,\dots, d\}$, we have $m_q(c)=m(\tau_i=c)$. By taking factorials, we get that $m_q(c)!=m(\tau_i=c)!$. Multiplying all these equations, we get that
\[\prod_{i=1}^d \prod_{c\in V} m(\tau_i=c)!=\left(\prod_{c\in V} m_q(c)!\right)^d,\]
that is,
\[\left(\prod_{i=1}^d \prod_{c\in V} m(\tau_i=c)!\right)^{\frac{d-1}d}=\left(\prod_{c\in V} m_q(c)!\right)^{d-1}.\]
\end{proof} 

Of course there are many other equivalent ways to express the quantity \break $|S(q)|\mathbb{P}((\Sigma(\bar{Q})=r)\wedge(m_{\bar{Q}}=m))$ and each of them suggests a way to extend the formula to all non-negative integral vectors, but the formula given in Equation~\eqref{Emdef} will be useful for us later.



\begin{lemma}\label{segedlem}
Consider a non-negative integral $W$-half-decent vector $m_0\in\mathbb{R}^{V^{d}}$, such that $\|m_0\|_{W^C}=O(\log n)$, where $n=\sum_{t\in V^d} m(t)$. For $u\in V^{d}$, let $\chi_u\in\mathbb{R}^{V^{d}}$ be such that $\chi_u(u)=1$ and $\chi_u(t)=0$ for every $t\neq u\in V^{d}$.
\begin{itemize}
\item If $u\in W^{d}$, then $E(m_0+\chi_u)/E(m_0)=O( \log^4 n);$
\item If $u\not\in W^{d}$, then $E(m_0+\chi_u)/E(m_0)=O( n^{-(d-2)/d}\log^{2} n).$
\end{itemize}

\end{lemma}
\begin{proof}

Let \[g=\frac{1+ m_0(\tau_\Sigma=u_\Sigma)}{1+m_0(u)}\qquad\text{ and }\qquad f_i=\frac{1+ m_0(\tau_i=u_i) }{n+1}.\]

Note that\[E(m_0+\chi_u)/E(m_0)=g\cdot\left(\prod_{i=1}^d f_i \right)^{\frac{d-1}{d}}.\]

If $u\in W^d$, then since $m_0$ is $W$-half-decent, we have $g\le\log^4 n$, and clearly $f_i\le 1$, so the statement follows.

If $u\not\in W^d$, we consider the following two cases:

\begin{enumerate}
\item If $u_\Sigma\not\in dW$, then 
\[ g\le 1+m_0(\tau_\Sigma=u_\Sigma)\le 1+ \|m_0\|_{W^c}=O(\log n),\]
 and there is an $i$ such that $u_i\not\in W$. This imply that $f_i=O\left(\frac{\log n}{n}\right)$. So 
\[E(m_0+\chi_u)/E(m_0)=O\left(\log n \left(\frac{\log n}{n}\right)^{\frac{d-1}{d}}\right)=O\left( n^{-\frac{d-2}{d}}\log^2 n\right).\] 
\item If $u_\Sigma\in W^d$, then there are at least two indices $i$ such that $u_i\not\in W$, for such an index $i$,  we have $f_i=O\left(\frac{\log n}{n}\right)$, clearly $g=O(n)$, so \[E(m_0+\chi_u)/E(m_0)=O\left(n\left(\frac{\log n}{n}\right)^{\frac{2(d-1)}{d}}\right)=O\left( n^{-\frac{d-2}{d}}\log^2 n\right).\]\qedhere
\end{enumerate}
\end{proof}

The next lemma follows easily from the previous one.
\begin{lemma}\label{DeltaMbecsElo}
There is a $D$, such that for any $i\in \{h+1,h+2,\dots,\ell\}$ and any  non-negative integral $W$-half-decent vector $m_0\in\mathbb{R}^{V^{d}}$, such that $\|m_0\|_{W^C}=O(\log n)$,  we have
\[\pushQED{\qed}E(m_0+m_i)/E(m_0)=O\left(\left(n^{-(d-2)/d}\log^D n\right)^{\|m_i\|_{W^C}}\right).\qedhere
\popQED\]
\end{lemma}

\begin{lemma}\label{WHDLemma}
Assume that $n$ is large enough.  Let  $q\in V^n$ be $(W,C\log n)$-typical, and let $m\in\mathcal{M}^\sharp(q,r_q)$. If $m_0$ is an integral vector indexed by $V^d$ such that \break $(m-\Delta(m))(t)\le m_0(t)\le m(t)$ for every $t\in V^d$, then $m$ is $W$-half-decent.

\end{lemma}
\begin{proof}
Let $L=\max_{i=h+1}^\ell \|m_i\|_\infty$. Note that $m(t)-m'(t)\le L\|m\|_{W^C}\le LC\log n$ for every $t\in V^d$.  Let $n_0=\sum_{t\in V^d} m_0(t)$. Then
\[n_0\ge n-L\cdot|V|^d\cdot \|m\|_{W^C}\ge n-L|V|^d C\log n.\]
If $n$ is large enough, then $LC\log^3 n\le  \frac{1}2 \log ^4 n_0$. We need to prove that 
\[\frac{1+ m_0(\tau_\Sigma=u_\Sigma)}{1+m_0(u)}\le \log^4 n_0,\]
for every $u\in W^d$. If $1+ m_0(\tau_\Sigma=u_\Sigma)\le \log^4 n_0$, then it is clear. Thus, assume that  $1+ m_0(\tau_\Sigma=u_\Sigma)> \log^4 n_0$. Then,
\begin{align*}
1+ m_0(\tau_\Sigma=u_\Sigma)&\le 1+m(\tau_\Sigma=u_\Sigma)\\&\le (1+m(u))\log^2 n\\
&\le (1+m_0(u)+LC\log n)\log^2 n\\&\le (1+m_0(u))\log^2 n +\frac{1}{2}\log^4 n_0\\&\le (1+m_0(u))\log^2 n +\frac{1}{2}\left(1+ m_0(\tau_\Sigma=u_\Sigma)\right).
\end{align*}
Therefore, if $n$ is large enough, then we have
\[\frac{1+ m_0(\tau_\Sigma=u_\Sigma)}{1+m_0(u)}\le 2\log^2 n\le \log^4 n_0.\]
\end{proof}

The following estimate will be crucial later. 

\begin{lemma}\label{DeltaMbecs}
There is a $K$ such that for any $(W,C\log n)$-typical $q\in V^n$ and $m\in\mathcal{M}^\sharp(q,r_q)$, we have
\[E(m)\le \left(Kn^{-(d-2)/d}\log^D n\right)^{\|\Delta(m)\|_{W^C}}E(m-\Delta(m)).\] 
\end{lemma}
\begin{proof}
We may assume that $n$ is large enough, because we can increase $K$ to handle the small values of $n$. Then the statement follows from repeated application of  Lemma \ref{DeltaMbecsElo}.  Observe that $m-\Delta(m)$ and all other $m_0$ we need to apply that lemma is $W$-half-decent by Lemma \ref{WHDLemma}.
\end{proof}

Now we made all the necessary preparations to prove Equation \eqref{eegyenlo3}. With our new notations, we have to prove that
\[\lim_{n\to\infty} \sum_{q\in D_n^W} \sum_{m\in \mathcal{M}^\sharp (q,r_q)} E(m)=0.\]
We prove it by induction on $|V|$. The statement is clear if $W=V$, because in that case $D_n^W$ is empty. So we may assume that $| W|<| V|$. 

\begin{lemma}
There is a finite $B=B_W$ such that for every $n$, we have that \[\sum_{q\in W^n\cap T_n} |S(q)|\mathbb{P}(A_n^{(d)}q=r_q) <B.\]
\end{lemma}
\begin{proof}
First consider the case when the coset $W$ is  a subgroup.  Then from the induction hypothesis, we can use Theorem \ref{FoFo} to get that that 

\[
\sum_{q\in W^n} \mathbb{P}(A_n^{(d)}q=r_q)=
\sum_{q\in W^n} \mathbb{P}(U_{q,d}=r_q)+o(1).
\]

Recall that for $W_0\in \Cosets(W)$, we defined $I(W_0^n)$ as \[I(W_0^n)=\{q\in W_0^n\quad|\quad\Aff_q=W_0\}.\] Now, we have
\begin{align*}
\sum_{q\in W^n} \mathbb{P}(U_{q,d}=r_q)&=  \sum_{W_0\in\Cosets(W)}\sum_{q\in I(W_0^n)}\mathbb{P}(U_{q,d}=r_q)\\&=\sum_{W_0\in\Cosets(W)} |I(W_0^n)|\cdot|W_0|^{-(n-1)}\le \sum_{W_0\in\Cosets(W)} |W_0|.
\end{align*}
Thus,
\begin{align*}
\sum_{q\in W^n\cap T_n} |S(q)|\mathbb{P}(A_n^{(d)}q=r_q) &= \sum_{q\in W^n} \mathbb{P}(A_n^{(d)}q=r_q)\\&=
\sum_{q\in W^n} \mathbb{P}(U_{q,d}=r_q)+o(1)\le \sum_{W_0\in \Cosets(W)}|W_0|+o(1).
\end{align*}
This proves the lemma when $W$ is a subgroup of $V$. If the coset $W$ is not a subgroup, then  we need to use the bijection given in the proof of Lemma \ref{tipikuskulonbseg}.
\end{proof}

We need a few notations, let
\[\mathcal{M}^\Delta_n=\cup_{q\in D_n^W} \{\Delta(m)\quad|\quad m\in\mathcal{M}^\sharp(q,r_q)\}.\]
For $m_{\Delta}\in \mathcal{M}^\Delta_n$ let 
\[\Delta^{-1}_n(m_{\Delta})=\cup_{q\in D_n^W} \{m\in \mathcal{M}^\sharp(q,r_q)\quad|\quad \Delta(m)=m_\Delta\}.\]

Using Lemma \ref{DeltaMbecs}, we obtain that
\begin{multline}\label{folytkov}
\sum_{q\in D_n^W} \sum_{m\in \mathcal{M}^\sharp(q,r_q)} E(m)=\sum_{m_\Delta\in\mathcal{M}^{\Delta}_n}\quad \sum_{m\in \Delta^{-1}_n(m_\Delta)} E(m)\le\\ \sum_{m_\Delta\in\mathcal{M}^{\Delta}_n} \left(Kn^{-(d-2)/d}\log^{D} n\right)^{\|m_\Delta\|_{W^C}}\sum_{m\in \Delta^{-1}_n(m_\Delta)} E(m-m_\Delta).
\end{multline}

Fix a vector  $m_{\Delta}\in \mathcal{M}^\Delta_n$. 
Set $n'=n-\sum_{t\in V^{d}}m_\Delta(t)$. 
Let $X$ be the set of $q\in D_n^W$, such that $\mathcal{M}^\sharp(q,r_q)\cap \Delta^{-1}_n(m_{\Delta})$ is non-empty. 

For each $q\in X$, there is a unique $q'\in W^{n'}\cap T_{n'}$ such that for every 
$c\in V$, we have $m_{q'}(c)=m_q(c)-m_{\Delta}(\tau_1=c)$, 
and a unique $w_q\in W^{n'}\cap T_{n'} $ such that for every $c\in V$, we have $m_{w_q}(c)=m_{r_q}(c)- m_{\Delta}(\tau_\Sigma=c)$. 

Note that for any $m\in \mathcal{M}^\sharp(q,r_q)\cap \Delta^{-1}_n(m_{\Delta})$, we have 
$m-m_\Delta\in \mathcal{M}(q',w_q)$. Moreover, \[E(m-m_\Delta)=|S(q')|\mathbb{P}((\Sigma(\bar{Q})=w_q)\wedge (m_{\bar{Q}}=m-m_\Delta)),\] where $\bar{Q}$ is a random $(q',d)$
-tuple. The map $m\mapsto m-m_\Delta$ is injective, so it follows that \[\sum_{m\in \mathcal{M}^\sharp(q,r_q)\cap \Delta^{-1}_n(m_{\Delta})} E(m-m_\Delta)\le |S(q')|\mathbb{P}(A_{n'}^{(d)}q'=w_q).\]  Also note that that the map $q\mapsto q'$ is injective. Therefore,

\begin{align*}
\sum_{m\in \Delta^{-1}_n(m_\Delta)} E(m-m_\Delta)&=\sum_{q\in X}\quad \sum_{m\in\mathcal{M}^\sharp(q,r_q)\cap \Delta^{-1}_n(m_{\Delta})} E(m-m_\Delta)\\&\le
\sum_{q\in X} |S(q')| \mathbb{P}(A_{n'}^{(d)}q'=w_q)\\&\le \sum_{q'\in W^{n'}\cap T_{n'}} |S(q')| \mathbb{P}(A_{n'}^{(d)}q'=r_{q'})<B .
\end{align*}

 Thus, continuing Inequality \eqref{folytkov}, we have
 \[\sum_{q\in D_n^W} \sum_{m\in \mathcal{M}^\sharp(q,r_q)} E(m)\le B \sum_{m_\Delta\in\mathcal{M}^{\Delta}_n} \left(Kn^{-(d-2)/d}\log^{D} n\right)^{\|m_\Delta\|_{W^C}}.\] 

There is an $F$ such that $|\mathcal{M}_n^\Delta|\le n^F$. We choose a constant $G$ such that for a large enough $n$, we have $\left(Kn^{-(d-2)/d}\log^{d-1} n\right)^{\|m_\Delta\|_{W^C}}<n^{-(F+1)}$, whenever \break $\|m_\Delta\|_{W^C}\ge G$. Let $H$ be the cardinality of the set
\[
\{m\quad|\quad m=\sum_{i=h+1}^{\ell} c_i m_i,\quad c_{h+1},c_{h+2},\dots,c_\ell\text{ non-negative integers, }\|m\|_{W^c}< G\}.
\]
Note that $H\le G^{\ell-h}$. 
Finally observe that $\|m_\Delta\|_{W^C}\ge 1$ for all $m_\Delta\in\mathcal{M}_n^\Delta$. So for large enough $n$
\begin{align*}
B &\sum_{m_\Delta\in\mathcal{M}^{\Delta}_n} \left(Kn^{-(d-2)/d}\log^D n\right)^{\|m_\Delta\|_{W^C}}\\
&=
B \sum_{\substack{m_\Delta\in\mathcal{M}^{\Delta}_n\\ \|m_\Delta\|_{W^C}\ge G}} \left(Kn^{-(d-2)/d}\log^{D} n\right)^{\|m_\Delta\|_{W^C}}\\
&\qquad +B \sum_{\substack{m_\Delta\in\mathcal{M}^{\Delta}_n\\ \|m_\Delta\|_{W^C}< G}} \left(Kn^{-(d-2)/d}\log^D n\right)^{\|m_\Delta\|_{W^C}}\\ &\le Bn^Fn^{-(F+1)}+BHK n^{-(d-2)/d}\log^{D} n=o(1).
\end{align*}

Thus, we have proved Equation \eqref{eegyenlo3}.

\section{The connection between the mixing property of the adjacency matrix and the sandpile group}\label{mixingboldist}

The random $(n-1)\times (n-1)$ matrix $A_n'$ is obtained from $A_n$ by deleting its last row and last column. For $q\in V^{n-1}$, the subgroup generated by $q_1,q_2,\dots,q_{n-1}$ is denoted by $G_q$. Let $U_q$ be a uniform random element of $G_q^{n-1}$. The next corollary of Theorem \ref{FoFo} states that the distribution of $A_n' q$ is close to that of $U_q$.

\begin{corr}\label{correduced}
We have
\[\lim_{n\to\infty} \sum_{q\in V^{n-1}} d_\infty(A_n'q,U_{q})=0.\]
\end{corr} 
\begin{proof}
For $q\in V^{n-1}$ and $r\in G_q^{n-1}$, we define $\bar{q}=(q_1,q_2,\dots,q_{n-1},0)\in V^n$ and \break $\bar{r}=(r_1,r_2,\dots,r_{n-1},d\cdot s(q)-s(r))\in G_q^n$. 

Note that $s(\bar{r})=d\cdot s(q)=d\cdot s(\bar{q})$ 
and $\Aff_{\bar{q}}=G_q$, 
so $\bar{r}\in R(\bar{q},d)$. 
Moreover, $A_n'q=r$ if and only if $A_n\bar{q}=\bar{r}$, so $\mathbb{P}(A_n'q=r)=\mathbb{P}(A_n\bar{q}=\bar{r})$. From these observations, it follows easily that $d_\infty(A_n' q,U_q)=d_\infty(A_n \bar{q},U_{\bar{q},d})$. The rest of the proof follows from Theorem \ref{FoFo}.
\end{proof}

Recall that the reduced Laplacian $\Delta_n$ of $D_n$ was defined as $\Delta_n=A_n'-dI$. The next well-known proposition connects $\Hom(\Gamma_n,V)$ and $\Sur(\Gamma_n,V)$ with the kernel of $\Delta_n$ when $\Delta_n$ acts on $V^{n-1}$.

\begin{prop}\label{prop1}
For  any finite abelian group $V$, we have 
\[|\Hom(\Gamma_n,V)|=|\{q\in V^{n-1}\quad|\quad \Delta_n q=0\}|\]
and
\[|\Sur(\Gamma_n,V)|=|\{q\in V^{n-1}\quad|\quad \Delta_n q=0,\quad G_q=V\}|.\]
\end{prop}   
\begin{proof}
There is an obvious bijection between $\Hom(\Gamma_n,V)$ and \[\{\varphi\in \Hom(\mathbb{Z}^{n-1},V)|\quad \RowS(\Delta_n)\subset \ker \varphi\}.\] Moreover, any $\varphi\in \Hom(\mathbb{Z}^{n-1},V)$ is uniquely determined by the vector \break $q=(\varphi(e_1),\varphi(e_2),\dots,\varphi(e_{n-1}))\in V^{n-1}$, where $e_1,e_2,\dots,e_{n-1}$ is the standard generating set of $\mathbb{Z}^{n-1}$. Furthermore, $\RowS(\Delta_n)\subset\ker \varphi$ if and only if $\Delta_n q=0$, so the first statement follows. The second one can be proved similarly.
\end{proof}

Combining  Proposition \ref{prop1} with with Corollary \ref{correduced}, we obtain
\begin{align*}
\lim_{n\to\infty}\mathbb{E}|\Sur(\Gamma_n,V)|&=\lim_{n\to\infty}\sum_{\substack{q\in V^{n-1}\\ G_q=V}}\mathbb{P}(\Delta_n q=0)=\lim_{n\to\infty}\sum_{\substack{q\in V^{n-1}\\ G_q=V}}\mathbb{P}(A_n' q=dq)\\&= \lim_{n\to\infty}\sum_{\substack{q\in V^{n-1}\\ G_q=V}}\mathbb{P}(U_q=dq)\\&=\lim_{n\to\infty}|\{q\in V^{n-1}|\quad G_q=V \}|\cdot |V|^{-(n-1)}=1.
\end{align*}

This proves  Theorem \ref{momentumokD}. 

To obtain Theorem \ref{CohenlenstraD} from this theorem, we need to use the results of Wood on the moment problem.

\begin{theorem}(Wood \cite[Theorem 3.1]{wood2} or  \cite[Theorem 8.3]{wood})\label{woodmoment}
Let $X_n$ and $Y_n$ be sequences of random finitely generated abelian groups.
Let $a$ be a positive integer and $A$ be the set of (isomorphism classes of) abelian groups with exponent dividing~$a$.
Suppose that for every $G\in A$, we have a number $M_G\leq |\wedge^2 G|$ such that
$$
\lim_{n\to \infty} \mathbb{E}| \Sur(X_n, G)| = M_G, 
$$
and
$$
\lim_{n\to \infty} \mathbb{E}| \Sur(Y_n, G)| = M_G. 
$$

Then for every $H\in A$,
 the limits
\[
\lim_{n\to\infty} \mathbb{P}(X_n\otimes \mathbb{Z}/a\mathbb{Z} \isom H)\quad\text{ and }\quad\lim_{n\to\infty} \mathbb{P}(Y_n\otimes \mathbb{Z}/a\mathbb{Z} \isom H)
\]
exist, and they are equal.
\end{theorem}
This has the following consequence.
\begin{theorem}
Let $p_1,p_2,\dots,p_s$ be  distinct primes. Let $X_n$ and $Y_n$ be sequences of random finitely generated abelian groups. Assume that for any finite abelian group $G$, we have a number $M_G\leq |\wedge^2 G|$ such that
$$
\lim_{n\to \infty} \mathbb{E}| \Sur(X_n, G)| = M_G, 
$$
and
$$
\lim_{n\to \infty} \mathbb{E}| \Sur(Y_n, G)| = M_G. 
$$ 
Let $X_{n,i}$ (resp. $Y_{n,i}$) be the $p_i$-Sylow subgroup  of $X_n$ (resp. $Y_n$). For $i=1,2,\dots,s$, let $G_i$ be a finite abelian $p_i$-group. Then   
the limits
\[\lim_{n\to\infty} \mathbb{P}\left(\bigoplus_{i=1}^s X_{n,i}\isom \bigoplus_{i=1}^s G_i\right)\quad\text{ and }\quad\lim_{n\to\infty} \mathbb{P}\left(\bigoplus_{i=1}^s Y_{n,i}\isom \bigoplus_{i=1}^s G_i\right)\]
exist, and they are equal.
\end{theorem}
\begin{proof}
Let $a_0$ be the exponent of the group $\bigoplus_{i=1}^s G_i$. Let $a=a_0\cdot\prod_{i=1}^s p_i$. Observe that 
$\bigoplus_{i=1}^s X_{n,i}\isom \bigoplus_{i=1}^s G_i$ if and only if $X_n\otimes \mathbb{Z}/a\mathbb{Z}\isom \bigoplus_{i=1}^s G_i$. Thus, the previous theorem gives the statement. 
\end{proof}
The next theorem gives two special cases which are of particular interest for us.
\begin{theorem}\label{momentproblemspec}
Let $p_1,p_2,\dots,p_s$ be  distinct primes. Let $\Gamma_n$ be sequence of random finitely generated abelian groups. Let $\Gamma_{n,i}$ be the $p_i$-Sylow subgroup of $\Gamma_n$. 
\begin{enumerate}
\item
Assume that for any finite abelian group $V$, we have
$$
\lim_{n\to \infty} \mathbb{E}| \Sur(\Gamma_n, V)| = 1. 
$$
 For $i=1,2,\dots,s$, let $G_i$ be a finite abelian $p_i$-group. Then   
\[
\lim_{n\to\infty} \mathbb{P}\left(\bigoplus_{i=1}^s \Gamma_{n,i}\isom \bigoplus_{i=1}^s G_i\right)=\prod_{i=1}^s \left(|\Aut(G_i)|^{-1} \prod_{j=1}^\infty (1-p_i^{-j})\right).
\]

\item
Assume that for any finite abelian group $V$, we have
$$
\lim_{n\to \infty} \mathbb{E}| \Sur(\Gamma_n, V)| = |\wedge^2 V|. 
$$
 For $i=1,2,\dots,s$, let $G_i$ be a finite abelian $p_i$-group. Then   
\begin{multline*}
\lim_{n\to\infty} \mathbb{P}\left(\bigoplus_{i=1}^s \Gamma_{n,i}\isom \bigoplus_{i=1}^s G_i\right)=\\\prod_{i=1}^s \left(\frac{|\{\phi:G_i \times G_i \to \mathbb{C}^* \SBP \}|}{|G_i||\Aut(G_i)|} \prod_{j=0}^\infty (1-p_i^{-2j-1})\right).
\end{multline*}
\end{enumerate}
\end{theorem}
\begin{proof}
The first part follows from the previous theorem and \cite[Lemma 3.2]{wood2} with the choice of $u=0$. Or alternatively, we can use the results of \cite[Section 8]{ellenb}. The second part follows from the previous theorem and \cite[Theorem 2 and Theorem 11]{clp14}. See also the proof of Corollary 9.2 in \cite{wood}. 
\end{proof}

Combining the first statement of the previous theorem with Theorem \ref{momentumokD}, we obtain Theorem \ref{CohenlenstraD}. The proofs of the corresponding statements about the sandpile group of $H_n$ are postponed to Section \ref{modificat} and \ref{parosparos}.

\section{A version of
Theorem \ref{FoFo} with uniform convergence}\label{secuniform}

We sate our results for the directed random graph model, but the arguments can be repeated for the undirected model as well. We write $A_n^{(d)}$ in place of $A_n$ to emphasize the dependence on $d$. We start by a simple lemma.
\begin{lemma}
For a fixed $n$ and $q\in V^n$, we have \[d_\infty(A_n^{(d)}q,U_{q,d})\le d_\infty(A_n^{(d-1)}q,U_{q,d-1}).\]
\end{lemma} 
\begin{proof}
Take any $r\in R(q,d)$. Observe that for $q'\sim q$, we have $r-q'\in R(q,d-1)$. Let $q'$ be a  uniform random element of $S(q)$ independent from $A_n^{(d-1)}$, then
\begin{align*}
|\mathbb{P}(A_n^{(d)}q=r)-\mathbb{P}(U_{q,d}=r)|&=|\mathbb{E}\mathbb{P}(A_n^{(d-1)}q=r-q')-|R(q,d)|^{-1}|\\&\le \mathbb{E}|\mathbb{P}(A_n^{(d-1)}q=r-q')-|R(q,d-1)|^{-1}|\\&\le d_\infty(A_n^{(d-1)}q,U_{q,d-1}).
\end{align*}
Note that here the expectations are over the random choice of $q'$. Since this is true for any $r\in R(q,d)$, the statement follows.
\end{proof}

Using this we can deduce the following uniform version of Theorem \ref{FoFo}.
\begin{corr}\label{FoFouni22}
 We have
 \[\pushQED{\qed}\lim_{n\to\infty} \sup_{d\ge 3} \sum_{q\in V^n} d_\infty(A_n^{(d)}q,U_{q,d})=0.\qedhere
\popQED\]
 \end{corr}

This also implies a uniform version of Corollary \ref{correduced}. Therefore, the limits in Theorem \ref{momentumokD} are uniform in $d$. Consequently, Theorem \ref{CohenlenstraD} remains true if we allow $d$ to vary with $n$.

\section{Sum of matching matrices: Modifications of the  proofs}\label{modificat}

A fixed point free permutation of order $2$ is called a matching permutation. The permutation matrix of a matching permutation is called matching matrix.  Then  $C_n=M_1+M_2+\dots+M_d$, where $M_1,M_2,\dots,M_d$ are independent uniform random $n\times n$  matching matrices.

Consider a vector $q=(q_1,q_2,...,q_n)\in V^n$. For a matching permutation $\pi$ of the set $\{1,2,\dots, n\}$ the vector $q_{\pi}=(q_{\pi(1)},q_{\pi(2)},\dots,q_{\pi(n)})$ is called a matching permutation of $q$.  A random matching permutation of $q$ is defined as the random variable $q_{\pi}$, where $\pi$ is chosen uniformly from the set of all matching permutations.

A $(q,1,h)$-tuple is a $1+h$-tuple $Q=(q^{(0)},q^{(1)},\dots,q^{(h)})$, where $q^{(0)}=q$ and $q^{(1)},q^{(2)},\dots,q^{(h)}$ are  matching permutations of $q$. A random $(q,1,h)$-tuple is a tuple  $\bar{Q}=(\bar{q}^{(0)},\bar{q}^{(1)},\dots\bar{q}^{(h)})$, where $\bar{q}^{(0)}=q$ and $\bar{q}^{(1)},\bar{q}^{(2)},\dots,\bar{q}^{(h)}$ are independent   random matching permutations of $q$. Similarly as before, a $(q,1,h)$-tuple can be viewed as a vector $Q=(Q_1,Q_2,\dots,Q_n)$ in $(V^{1+h})^n$. For $t\in V^{1+h}$, we define 
\[m_Q(t)=|\{i\quad|\quad Q_i=t\}|.\]
In this section the components of a vector $t\in V^{1+h}$ are indexed from $0$ to $h$,\break that is,  $t=(t_0,t_1,\dots,t_h)$. For $t\in V^{1+h}$, we define $t_\Sigma=\sum_{i=1}^n t_i$. 
The sum $\Sigma(Q)$ of a $(q,1,h)$-tuple $Q$ is defined as $\Sigma(Q)=\sum_{i=1}^h q^{(i)}$. Note that the sums above do not include $t_0$ and $q^{(0)}$. 

We define 
\[\mathcal{M}^S(q,r)=\{m_Q|\quad Q\text{ is a }(q,1,h)\text{-tuple such that }\Sigma(Q)=r\}.\]

A $(q,1,h)$-tuple $Q$ is $\gamma$-typical if $\left\|m_Q-\frac{n}{|V|^{1+h}}\mathbbm{1}\right\|_\infty<n^{\gamma}$. 

For two vectors $q,r\in V^n$ and $a,b\in V$, we define 
\[m_{q,r}(a,b)=|\{i|\quad q_i=a\text{ and }r_i=b\}|.\]

The vector $r$ is called $(q,\beta)$-typical if \[\left\|m_{q,r}-\frac{n}{|V|^{2}}\mathbbm{1}\right\|_\infty<n^{\beta}.\] 

With these notations, we have the following analogue of Theorem \ref{mixingThm}. 

\begin{theorem}\label{mixingThmk}
For any fixed finite abelian group $V$ and $h\ge 2$,  we have
\[\lim_{n\to\infty} \sup_{\substack{q\in V^n\quad \alpha-\text{typical}\\r\in R^S(q,h)\quad (q,\beta)-\text{typical}}} \left|\mathbb{P}(C_n^{(h)}q=r)\Big/\left(\frac{2^{\rank2(V)}|\wedge^2 V|}{|V|^{n-1}}\right)-1\right|=0.\]
\end{theorem}
\begin{proof}
The proof is analogous with the proof of Theorem \ref{mixingThm}. We need to replace the notion of $(q,h)$-tuple with the notion of $(q,1,h)$-tuple, the notion of $\beta$-typical vector with the notion of $(q,\beta)$-typical vector. Moreover, some of the statements should be slightly changed. Now we list the modified statements. 

We start by determining the size of $R^S(q,h)$.
\begin{lemma}
Let $q\in V^n$ such that $\Aff_q=V$, then  \[|R^S(q,h)|=\frac{|V|^{n-1}}{2^{\rank2(V)}|\wedge^2 V|} .\]
\end{lemma}
\begin{proof} 
We define the homomorphism $\varphi:V^n\to (V\otimes V)\times V$ by setting \[\varphi(r)=(<q\otimes r>,s(r))\] for every $r\in V^n$. We claim that it is surjective. First, take any $a,b\in V$. The condition $\Aff_q=V$  implies that $q_1-q_n,q_2-q_n,\dots,q_{n-1}-q_n$ generate $V$. In particular, there are integers $c_1,c_2,\dots, c_{n-1}$ such that $a=\sum_{i=1}^{n-1} c_i(q_1-q_n)$. Let us define $r=(c_1b,c_2b,\dots,c_{n-1}b,-\sum_{i=1}^{n-1} c_i b)\in V^n$. Then
\[<q\otimes r>=\sum_{i=1}^{n-1}q_i\otimes c_i b+q_n\otimes \left(-\sum_{i=1}^{n-1} c_i b\right)=\left(\sum_{i=1}^{n-1} c_i(q_i-q_n)\right)\otimes b=a\otimes b,\]
and $s(r)=0$, that is, $\varphi(r)=(a\otimes b,0)$. Thus, $V\otimes V\times\{0\}$ is contained in the range of $\varphi$. 

Now take any $(x, v)\in (V\otimes V)\times V$. Clearly, we can pick an $r_1\in V^n$ such that $s(r_1)=v$. Then from the previous paragraph, there is an $r_2$ such that \break $\varphi(r_2)=(x-<q\otimes r_1>, 0)$. Then $\varphi(r_1+r_2)=(x,v)$. This proves that $\varphi$ is indeed surjective.
Since $R^S(q,h)=\varphi^{-1}(I_2\times\{h\cdot s(q)\})$, we have

\[|R^S(q,h)|=\frac{|I_2|}{|(V\otimes V)|\cdot| V|}|V|^n=\frac{|V|^{n-1}}{2^{\rank2(V)}|\wedge^2 V|}.\]     
\end{proof}

\begin{lemma}[The analogue of Lemma \ref{euq12lemma}]\label{matching3}
Consider  $q,r\in V^n$. 
Let $m\in \mathcal{M}^S (q,r)$. Then $m$ is a nonnegative integral vector with the following properties.



\begin{align}
\label{equ1k}
&m(\tau_0=a\text{ and }\tau_i=b)=m(\tau_0=b\text{ and }\tau_i=a)&\forall i\in\{1,2,\dots,h\},\quad a,b\in V,\\
\label{equ2k}
&m(\tau_0=a\text{ and }\tau_\Sigma=b)=m_{q,r}(a,b)&\forall a,b\in V.\\
\intertext{Moreover,}
\label{equ3k}
&m(\tau_0=c\text{ and }\tau_i=c)\text{ is even}&\forall i\in\{1,2,\dots,h\},\quad c\in V.
\end{align}

Now assume that $m$ is a nonnegative integral vector satisfying the conditions above. Then

\begin{multline}\label{eprk}
\mathbb{P}(\Sigma(\bar{Q})=r\text{ and }m_{\bar{Q}}=m)=\\
\left(\frac{n!}{2^{n/2} (n/2)!}\right)^{-h} \frac{\prod_{a,b\in V} m(\tau_0=a,\tau_\Sigma=b)!}{\prod_{t\in V^{1+h}} m(t)!} \times\\ \prod_{i=1}^h\left(\left(\prod_{a\in V} \frac{m(\tau_i=a,\tau_0=a)!}{2^{m(\tau_i=a,\tau_0=a)/2} (m(\tau_i=a,\tau_0=a)/2)!}\right)\left(\prod_{a\neq b\in V} \sqrt{m(\tau_0=a,\tau_i=b)!}\right)\right).
\end{multline}

In particular, $\mathbb{P}((\Sigma(\bar{Q})=r)\wedge (m_{\bar{Q}}=m))>0$ so $m\in 
\mathcal{M}^S(q,r)$. Let $A^S(q,r)$ be the affine subspace  given by the linear equations 
\eqref{equ1k} and \eqref{equ2k}   
above. Then  $\mathcal{M}^S(q,r)$ is the set of non-negative 
integral points of the affine subspace $A^S(q,r)$ satisfying  the parity constraints in \eqref{equ3k} above.    
\end{lemma}
\begin{proof}
We only give the proof of Equation \eqref{eprk}, since all the other statements of the lemma are straightforward to prove. The number of $(q,1,h)$-tuples $Q$ such that $\Sigma(Q)=r$ and $m_Q=m$ is
\[\frac{\prod_{a,b\in V} m(\tau_0=a,\tau_\Sigma=b)!}{\prod_{t\in V^{1+d}} m(t)!}.\]

Fix any $(q,1,h)$-tuple $Q=(q^{(0)},q^{(1)},\dots,q^{(h)})$ such that $\Sigma(Q)=r$ and $m_Q=m$. Now, we calculate the probability that $\mathbb{P}(\bar{Q}=Q)$ for a random $(q,1,h)$-tuple $\bar{Q}$. For $i\in \{1,2,\dots,h\}$ and $a,b\in V$, we define 

\[I_{i,a,b}=\{j\in\{1,2,\dots,n\}\quad|\quad q^{(i)}_j=a\text{ and }q^{(0)}_j=b \}.\]

First, for $i=1,2,\dots, h$, we determine the number of matching permutations $\pi$ such that $q_\pi=q^{(i)}$. In other words, we are interested in the number of perfect matchings $M$ on the set $\{1,2,\dots,n\}$ such that
\begin{enumerate}
\item For every $a\in V$, the restriction of $M$ to the set $I_{i,a,a}$  
is a perfect matching.
\item For every unordered pair $\{a,b\}\subset V$, where $a\neq b$, the restriction of $M$ gives a perfect matching between the disjoint set $I_{i,a,b}$ and $I_{i,b,a}$.
\end{enumerate} 
Since $|I_{i,a,a}|=m(\tau_i=a,\tau_0=a)$, we have
\[\frac{m(\tau_i=a,\tau_0=a)!}{2^{m(\tau_i=a,\tau_0=a)/2} (m(\tau_i=a,\tau_0=a)/2)!}\]
perfect matchings on the set $I_{i,a,a}$.  

For every unordered pair $\{a,b\}\subset V$, where $a\neq b$, let \[n_{i,\{a,b\}}=m(\tau_i=a,\tau_0=b)=m(\tau_i=b,\tau_0=a)\] be the common size of $I_{i,a,b}$ and $I_{i,b,a}$. Then there are \[n_{i,\{a,b\}}!=\sqrt{m(\tau_i=a,\tau_0=b)!}\cdot \sqrt{m(\tau_i=a,\tau_0=b)!}\]
perfect matchings between $I_{i,a,b}$ and $I_{i,b,a}$. We choose to express $n_{i,\{a,b\}}!$ as above, because this way we get a symmetric expression. 

Since the total number perfect matchings is $\frac{n!}{2^{n/2} (n/2)!}$, we obtain that for a uniform random matching matrix $M$, we have
\begin{multline*}\mathbb{P}(Mq=q^{(i)})=\left(\frac{n!}{2^{n/2} (n/2)!}\right)^{-1}  \\\times \left(\prod_{a\in V} \frac{m(\tau_i=a,\tau_0=a)!}{2^{m(\tau_i=a,\tau_0=a)/2} (m(\tau_i=a,\tau_0=a)/2)!}\right)\left(\prod_{a\neq b\in V} \sqrt{m(\tau_0=a,\tau_i=b)!}\right).
\end{multline*}
From this, Equation \eqref{eprk} follows easily.
\end{proof}

\begin{lemma}[The analogue of Lemma \ref{eltolas}]
For any $q,r_1,r_2\in V^n$, we define the vector  $v=v_{q,r_1,r_2}\in\mathbb{R}^{V^{1+h}}$  by
\[v(t)=\frac{m_{q,r_2}(t_0,t_\Sigma)-m_{q,r_1}(t_0,t_\Sigma)}{|V|^{h-1}}\]
for every $t\in V^{1+h}$. Then we have

\[\pushQED{\qed}A^S(q,r_1)+v_{q,r_1,r_2}=A^S(q,r_2).\qedhere
\popQED\]
\end{lemma}

\begin{lemma}[The analogue of Lemma \ref{vanegesz}]\label{vanegesz2}
Assume that $n$ is large enough. For an $\alpha$-typical vector $q\in V^n$ and $r\in R^S(q,h)$, the affine subspace $A^S(q,r)$ contains an integral vector satisfying the parity constraints in \eqref{equ3k} of Lemma \ref{matching3}.
\end{lemma}

To prove Lemma \ref{vanegesz2} we need a few lemmas. The group $V$ has a decomposition $V=\bigoplus_{i=1}^\ell <v_i>$ such that $o_1|o_2|\cdots|o_\ell$, where $o_i$ is order of $v_i$. 

\begin{lemma}\label{vanszimmetrikus}
Let $q\in V^n$ be such that $m_q(v_i)>0$ for every $1\le i \le \ell$. Let $r\in V^n$ such that $<q\otimes r>\in I_2$.  Then there is a symmetric matrix $A$ over $\mathbb{Z}$ such that $r=Aq$ and all the diagonal entries of $A$ are even. 
\end{lemma}
\begin{proof}
 We express $q_k$ as $q_k=\sum_{i=1}^\ell q_k(i) v_i$, and similarly we express $r_k$ as \break $r_k=\sum_{i=1}^\ell r_k(i) v_i$, where $q_k(i),r_k(i)\in \mathbb{Z}$. The condition that $<q\otimes r>\in I_2$   is equivalent to the following. For $1\le i\le j\le \ell$, we have
\begin{equation}\label{tensorequiv}
\sum_{k=1}^n q_k(i)r_k(j)\equiv \sum_{k=1}^n q_k(j)r_k(i)\pmod{o_i}
\end{equation} 
and whenever  $o_i$ is even, we have
\begin{equation}\label{tensorequiv2}
\sum_{k=1}^n q_k(i)r_k(i)\text{ is even.}
\end{equation}  
Due to symmetries and the fact that $m_q(v_i)>0$ for every $1\le i\le \ell$, we may assume that $q_i=v_i$ for $1\le i\le \ell$. We define the symmetric matrix $A=(a_{ij})$ by
\[
a_{ij}=
\begin{cases}
r_i(j)&\text{for }\ell<i\le n\text{ and } 1\le j\le\ell,\\
r_j(i)&\text{for }1\le i\le \ell \text{ and } \ell< j\le n,\\
0& \text{for }\ell<i\le n \text{ and } \ell< j\le n,\\
r_i(j)+r_j(i)-\sum_{k=1}^n q_k(j) r_k(i) &\text{for }1\le i\le j\le\ell,\\
r_i(j)+r_j(i)-\sum_{k=1}^n q_k(i) r_k(j) &\text{for }1\le j< i\le\ell.
\end{cases}
\]

From Equation \eqref{tensorequiv} we obtain that for $1\le j<i\le\ell$, we have
\[a_{ij}\equiv r_i(j)+r_j(i)-\sum_{k=1}^n q_k(j) r_k(i) \pmod{o_j}.\]
In particular, $a_{ij} q_j=a_{ij} v_j=(r_i(j)+r_j(i))v_j-\sum_{k=1}^n q_k(j) r_k(i)v_j$ for every \break$1\le i,j\le \ell$.

Let $w=Aq$. We need to prove that $w_i=r_i$ for every $1\le i\le n$. It is easy to see for $i>\ell$. Now assume that $i\le \ell$. Then
\begin{align*}
w_i&=\sum_{h=1}^\ell \sum_{j=1}^n a_{ij}q_j(h)v_h=\sum_{h=1}^\ell \left(a_{ih}v_h+\sum_{j=\ell+1}^n r_j(i)q_j(h)v_h\right)\\&=
\sum_{h=1}^\ell \left(r_i(h)+r_h(i)-\sum_{k=1}^n q_k(h) r_k(i)+\sum_{j=\ell+1}^n r_j(i)q_j(h)\right)v_h\\&= \sum_{h=1}^\ell \left(r_i(h)+r_h(i)-\sum_{k=1}^\ell q_k(h) r_k(i) \right)v_h=\sum_{h=1}^\ell r_i(h) v_h=r_i.
\end{align*}
Now we modify $A$ slightly to achieve that all the diagonal entries are even. If $i>\ell$, then $a_{ii}=0$ which is even. If $1\le i\le \ell$ and $o_i$ is even, then 
$a_{ii}=2r_i(i)-\sum_{k=1}^n q_k(i) r_k(i)$, which is even using the condition \eqref{tensorequiv2} above. 
If $1\le i\le \ell$, $o_i$ is odd and $a_{ii}$ is odd, we replace $a_{ii}$ by $a_{ii}+o_i$, this way we can achieve that $a_{ii}$ is even, without changing $Aq$. To see this, observe that $o_i q_i=o_i v_i=0$.  
\end{proof}


For $q,w\in V^n$ and $c\in V$, we define
\[z_{q,w}(c)=\sum_{\substack{1\le i\le n\\q_i=c}} w_i.\]
Note that $<q\otimes w>=\sum_{c\in V} c\otimes z_{q,w}(c)$.

\begin{lemma}\label{lemmavm}
Let $q\in V^n$ such that $m_q(c)>10|V|^2$ for every $c\in V$, and let $z\in V^V$. Then there is an $m$-permutation  $w$ of $q$ such that $z_{q,w}=z$, if and only if
\begin{equation}\label{vm1}
\sum_{c\in V}z(c) =s(q)
\end{equation}
and
\begin{equation}\label{vm2}
\sum_{c\in V} c\otimes z(c)\in I_2.
\end{equation}
\end{lemma}
\begin{proof}
It is clear that the conditions are indeed necessary, so we only need to prove the other direction. Since $m_q(c)>0$ for all $c\in V$, we can find a $w_0$ such that $z_{q,w_0}=z$. (Of course $w_0$ is not necessarily  a matching permutation of $q$.) Condition \eqref{vm2} gives us that $<q\otimes w_0>\in I_2$. Using Lemma \ref{vanszimmetrikus}, it follows that there is a symmetric matrix $A=(a_{ij})$, such that $Aq=w_0$ and all the diagonal entries of $A$ are even. 
For $a,b\in V$ we define
\[m_0(a,b)=\sum_{\substack{1\le i,j\le n\\q_i=a,\quad q_j=b}} a_{ij}.\]
Since $A$ is symmetric and the diagonal entries are even, we have $m_0(a,b)=m_0(b,a)$ and $m(a,a)$ is even for every $a,b\in V$. 

Let $m=m_0$. Replace $m(a,b)$ by  $m(a,b)-2\ell |V|$, where $\ell$ is an integer chosen such that $0\le m(a,b)-\ell 2|V|<2|V|$.  Now for every $0\neq a\in V$, we do the following procedure. We find the unique integer $\ell$ such that for 
\[\Delta=m_q(a)-\sum_{b\in V} m(a,b)-\ell 2|V|,\]
we have $0\le\Delta<2|V|$. Now increase $m(a,a)$ by $\ell 2|V|$. (Note that $\ell$ is non-negative because of  the condition $m_q(a)>10|V|^2$.)  
 Increase both $m(a,0)$ and $m(0,a)$ by $\Delta$. Finally, let $\Delta_0=m_q(0)-\sum_{b\in V} m(0,b)$, and increase $m(0,0)$ by $\Delta_0$. (Once again $\Delta_0$ is non-negative because of the condition $m_q(a)>10|V|^2$.) 

This way we achieved that for every $a\in V$, we have $\sum_{b\in V} m(a,b)=m_q(a)$. It is clear that $m(a,b)$ is a non-negative integer and $m(a,b)=m(b,a)$ for every $a,b\in V$. Moreover, $m(a,a)$ is even for $0\neq a\in V$. It is also true for $a=0$, but this requires some explanation. Indeed, $m(0,0)$ can be expressed as
\begin{align*}
m(0,0)&=\sum_{a,b\in V}m(a,b)-2\sum_{\substack{\{a,b\}\\a\neq b\in V}} m(a,b)-\sum_{0\neq a\in V} m(a,a)\\&=n-2\sum_{\substack{\{a,b\}\\a\neq b\in V}}m(a,b)-\sum_{0\neq a\in V} m(a,a).  
\end{align*}
Here in the last row, every term is even, so $m(0,0)$ is even too. From these observations, it follows that there is an $m$-permutation $w$ of $q$ such that $m_{q,w}=m$. We will prove that $z_{q,w}=z$. Consider an $0\neq a\in V$. Observe that $m(a,b)\equiv m_0(a,b)$ modulo $|V|$ for $b\neq 0$. Thus, 

\begin{align*}
z_{q,w}(a)&=\sum_{\substack{1\le i\le n\\q_i=a}} w_i=\sum_{b\in V} m_{q,w}(a,b)b=\sum_{b\in V}m_0(a,b) b=
\sum_{b\in V} \sum_{\substack{1\le i,j\le n\\q_i=a,\quad q_j=b}} a_{ij}b\\&= \sum_{b\in V}\sum_{\substack{1\le i,j\le n\\q_i=a,\quad q_j=b}} a_{ij}q_j=
\sum_{\substack{1\le i\le n\\q_i=a}} \sum_{j=1}^n a_{ij}q_j=\sum_{\substack{1\le i\le n\\q_i=a}} w_0(i)=z_{q,w_0}(a)=z(a).
\end{align*}

Finally 
\begin{align*}
z_{q,w}(0)&=\sum_{a\in V} z_{q,w}(a)-\sum_{0\neq a\in V} z_{q,w}(a)=\sum_{i=1}^n q_i-\sum_{0\neq a\in V} z_{q,w}(a)\\&=
s(q)-\sum_{0\neq a\in V} z(a)=\sum_{a\in V} z(a)-\sum_{0\neq a\in V} z(a)=z(0),
\end{align*}
using condition \eqref{vm1}.
\end{proof}

The proof of Lemma \ref{vanegesz} also gives us the following statement.
\begin{lemma}\label{vanegesz3}
Let $q_1,q_2,\dots,q_h\in V^n$ and $r\in V^n$. Assume that $\sum_{i=1}^n s(q_i)=s(r)$. Then there is an integral vector $m$ indexed by $V^h$ such that\footnote{Unlike in the rest of this section, here  the components of a $t\in V^h$ are indexed from $1$ to $h$.} 
\[m(\tau_i=b)=m_{q_i}(b)\]
for every $i=1,2,\dots, h$ and $b\in V$, and
\[m(\tau_\Sigma=b)=m_{r}(b)\]
for every $b\in V$.\qed
\end{lemma} 

Now we are ready to prove Lemma \ref{vanegesz2}.

\begin{proof}
Fix an $\alpha$-typical $q$, and $r\in R^S(q,h)$. Let $W$ be the set of $z\in V^V$ satisfying the conditions \eqref{vm1} and  \eqref{vm2} of Lemma \ref{lemmavm}. Observe that $W$ is a coset of $V^V$. Moreover, $r\in R^S(q,h)$ implies that $z_{q,r}\in hW$. Thus, we can find $z_1,z_2,\dots,z_h\in W$ such that $z_{q,r}=\sum_{i=1}^h z_i$. 
If $n$ is large enough, then for an \break $\alpha$-typical $q$, we have $m_q(c)>10|V|^2$.  By using Lemma \ref{lemmavm}, for each $i\in\{1,2,\dots,h\}$ we can find a matching permutation $w_i$ of $q$ such that $z_{q,w_i}=z_i$. For $a\in V$, let \break $w_i^a\in V^{m_q(a)}$ be the vector obtained from $w_i$ by projecting to the coordinates in the set  $\{i|\quad q_i=a\}$. Similarly, $r^a$ is obtained from $r$ by projecting to the same set of coordinates. Observe that $\sum_{i=1}^h s(w^a_i)=\sum_{i=1}^h z_i(a)=z_{q,r}(a)=s(r^a)$. Thus, from Lemma \ref{vanegesz3}, we obtain an integral vector $m^a$ indexed by $V^h$ such that

\[m^a(\tau_i=b)=m_{w^a_i}(b)=m_{q,w_i}(a,b)\]
for every $i=1,2,\dots, h$ and $b\in V$, and
\[m^a(\tau_\Sigma=b)=m_{r^a}(b)=m_{q,r}(a,b)\]
for every $b\in V$.

Then the vector $m$ 
defined by
\[m((t_0,1_1,\dots,t_h))=m^{t_0}((t_1,\dots,t_h))\]
gives us an integral point in $A^S(q,r)$  satisfying the parity constraints in \eqref{equ3k} of Lemma \ref{matching3}.   
\end{proof}

\begin{lemma}[The analogue of Lemma \ref{lemma20}]
For an $\alpha$-typical $q\in V^n$, a \break $(q,\beta)$-typical $r\in R^S(q,h)$, $r_0=r_0(q)$ and $m\in \mathcal{M}^{S*}(q,r_0)$, we have that
\[\mathbb{P}((\Sigma(\bar{Q})=r_0)\wedge (m_{\bar{Q}}=m))\sim \mathbb{P}((\Sigma(\bar{Q})=r)\wedge (m_{\bar{Q}}=m+\hat{v}_{q,r_0,r}))\]
uniformly.
\end{lemma} 
\begin{proof}
For any $\alpha$-typical $q\in V^n$,  $(q,\beta)$-typical $r\in R^S(q,h)$ and $m\in \mathcal{M}^{S*}(q,r)$, we have
\begin{multline*}
\mathbb{P}(\Sigma(Q)=r\text{ and }m_Q=m)\sim 
f(q) \exp\left(\frac{1}{2n} B\left(m-\frac{1}{|V|^{h+1}}\mathbbm{1},m-\frac{1}{|V|^{h+1}}\mathbbm{1}\right)\right)
\end{multline*}
uniformly, where $f(q)$ is some function of $q$ and the bilinear form $B(x,y)$ is defined as
\begin{multline*}
B(x,y)=-|V|^{1+h}\sum_{t\in V^{1+h}} x(t) y(t)+\frac{|V|^2}{2}\sum_{i=1}^h \sum_{a,b\in V} x(\tau_0=a,\tau_i=b) y(\tau_0=a,\tau_i=b)\\+|V|^2\sum_{a,b\in V} x(\tau_0=a,\tau_\Sigma=b) y(\tau_0=a,\tau_\Sigma=b).
\end{multline*}
 
The statement follows from the fact that $v_{q,r_0,r}$ is in the radical of $B$.
\end{proof}

\begin{lemma}[The analogue of Lemma \ref{OTEM} part \eqref{OTEM4}]
The following holds
 \[\lim_{n\to\infty} \sup_{\substack{q\in V^n\quad \alpha-\text{typical}\\r\in R^S(q,h)\quad (q,\beta)-\text{typical}}} \mathbb{P}\left((\Sigma(\bar{Q})=r)\wedge(\bar{Q}\text{ is not }\gamma-\text{typical})\right)|V|^{n}=0.\]
\end{lemma}
\begin{proof}
Take any $\alpha$-typical $q\in V^n$ and $(q,\beta)$-typical $r\in R^S(q,h)$. We define
\[S(q,r)=\{r'\in V^n|\quad m_{q,r'}=m_{q,r}\}.\]
From symmetry, it follows that $\mathbb{P}\left((\Sigma(\bar{Q})=r')\wedge(\bar{Q}\text{ is not }\gamma-\text{typical})\right)$ is the same for every $r'\in S(q,r)$. Thus,
\[\mathbb{P}\left((\Sigma(\bar{Q})=r)\wedge(\bar{Q}\text{ is not }\gamma-\text{typical})\right)\le\frac{\mathbb{P}(\bar{Q}\text{ is not }\gamma-\text{typical})}{|S(q,r)|}.\]
Since there is $c>0$ such that $|S(q,r)|\ge |V^n| \exp(-c n^{2\beta-1})$, the statement follows as in the proof of Lemma \ref{tech3}. 
\end{proof}

This concludes the proof of Theorem \ref{mixingThmk}.

\end{proof}


The analogue of Theorem \ref{mixingThm2} is the following.

\begin{theorem}\label{mixingThm2k}
For any fixed finite abelian group $V$ and  $d\ge 3$,  we have
\[\lim_{n\to\infty} |V|^n \sup_{q\in V^n\quad \alpha-\text{typical}} d_\infty(C_n^{(d)}q,U^S_{q,d})=0.\]
\end{theorem}

This theorem follows immediately from Theorem \ref{mixingThmk} once we prove the following analogue of Lemma \ref{lemmamixingthm2}.

\begin{lemma}\label{elhalasztott}
Let $q\in V^n$ be $\alpha$-typical, $r\in V^n$, $h\ge 2$ and $Q$ is a random $(q,h)$-tuple. Then there is a polynomial $g$ and a constant $C$ (not depending on $q$ and $r$), such that
\[\mathbb{P}(\Sigma(Q)=r)\le g(n) |V|^{-n}\exp(Cn^{2\alpha-1}).\] 
\end{lemma}
This will be proved after Lemma \ref{Eesp}, because the proofs of these two lemmas share some ideas.

Once we have Theorem \ref{mixingThm2k}, we only need to control the non-typical vectors to obtain Theorem \ref{FoFo22}. This can be done almost the same way as in Section \ref{Secmom}. Here we list the necessary modifications. 

In the next few lemmas, our main tool will be the notion of \emph{Shannon entropy}. Given a random variable $X$ taking values in a finite set $\mathcal{R}$, its Shannon entropy $H(X)$ is defined as
\[H(X)=\sum_{r\in \mathcal{R}} -\mathbb{P}(X=r)\log\mathbb{P}(X=r).\]
In the rest of this discussion, we always assume that random variables have finite range, and all the random variables are defined on the same probability space. If $X_1,X_2,\dots,X_k$ is a sequence random variables, then their joint Shannon entropy $H(X_1,X_2,\dots,X_k)$ is defined as the Shannon entropy $H(X)$ of the (vector valued) random variable $X=(X_1,X_2,\dots,X_k)$. See \cite{information} for more information on Shannon entropy. 

A few basic properties of Shannon entropy are given in the next lemma.
\begin{lemma}\label{entl}
Let $X,Y,Z$ be three random variables. Then
\begin{equation}\label{szubad}
H(X,Y)\le H(X)+H(Y),
\end{equation}
and
\begin{equation}\label{szubmod}
H(X,Z)+H(Y,Z)\ge H(Z)+H(X,Y,Z).
\end{equation}
Let $X,Y$ be two random variables such that $Y$ is a function of $X$. Then
\[H(X,Y)=H(X).\]
\end{lemma}   
\begin{proof}
Note that the quantity $H(X,Z)+H(Y,Z)-H(Z)-H(X,Y,Z)$ is usually denoted by $I(X;Y|Z)$ and it is called conditional mutual information. It is well known that $I(X;Y|Z)\ge 0$. See \cite[(2.92)]{information}. This proves Inequality \eqref{szubmod}. We can obtain Inequality \eqref{szubad} as a special case of Inequality \eqref{szubmod}, if we we choose $Z$ to be constant. The last statement is straightforward from the definitions.
\end{proof}

Later  we will need the following lemma.

\begin{lemma}\label{shaererszeru}
For $d\ge 1$, let $Y_0,Y_1,\dots,Y_d$ be $d+1$ random variables. Then
\[H(Y_0,Y_1,\dots,Y_d)\le\sum_{i=1}^d H(Y_0,Y_i)-(d-1)H(Y_0).\] 
\end{lemma} 
\begin{proof}
The statement can be proved by induction. Indeed, from Inequality \eqref{szubmod}, we have
\[H(Y_0,Y_1,\dots,Y_d)+H(Y_0)\le H(Y_0,Y_1,\dots,Y_{d-1})+H(Y_0,Y_d).\]
Therefore,
\begin{align*}
H(Y_0,Y_1,\dots,Y_d)&\le H(Y_0,Y_1,\dots,Y_{d-1})+H(Y_0,Y_d)-H(Y_0)\\&\le \sum_{i=1}^d H(Y_0,Y_i)-(d-1)H(Y_0),
\end{align*}
where in the last step we used the induction hypothesis.
\end{proof}

In Section \ref{Secmom}, we used the fact that $|S(q)|\mathbb{P}(A_n^{(d)}q=r)=\mathbb{P}(r-A_n^{(d-1)}q\sim q)$. This equality is replaced by the following lemma.

\begin{lemma}\label{Eesp}
Let $q,r\in V^n$ and \[m\in \mathcal{M}^S(q,r)=\{m_Q|\quad Q\text{ is a $(q,1,d)$-tuple and }\Sigma(Q)=r \}.\] We define 
\[E(m)=|S(q)|\mathbb{P}(m_{\bar{Q}}=m\text{ and }\Sigma(\bar{Q})=r),\]
where $\bar{Q}$ is random $(q,1,d)$-tuple. 

Moreover, let $p(m)$ be the probability of the event that for a random $(q,1,d-1)$-tuple $\bar{Q}=(\bar{q}^{(0)},\bar{q}^{(1)},\dots,\bar{q}^{(d-1)})$, we have that $r-\Sigma(\bar{Q})$ is a matching permutation of $q$ and the $(q,1,d)$-tuple $Q'=(\bar{q}^{(0)},\bar{q}^{(1)},\dots,\bar{q}^{(d-1)},r-\Sigma(\bar{Q}))$ satisfies $m_{Q'}=m$.
Then there is a polynomial $f(n)$ (not depending on $q,r$ or $m$) such that
\[E(m)\le f(n) p(m)^{\frac{1}{d-1}}.\]
Furthermore, there is a polynomial $g(n)$ such that
\[|S(q)|\mathbb{P}(C_n^{(d)}q=r)\le g(n) \mathbb{P}(r-C_n^{(d-1)}q\sim q)^{\frac{1}{d-1}}.\] 

\end{lemma} 
\begin{proof}
Let $X=(X_0,X_1,X_2,\dots, X_d)\in V^{1+d}$ be a random variable, such that $\mathbb{P}(X=t)=\frac{m(t)}{n}$ for every $t\in V^{1+d}$. We define $X_\Sigma=\sum_{i=1}^d X_i$. 
Then

\[E(m)=c_1(m)\exp\left(n\left(H(X_0)+H(X)-H(X,X_\Sigma)-\frac{1}2\sum_{i=1}^d H(X_0,X_i)\right)\right),\]
and 
\[p(m)=c_2(m)\exp\left(n\left(H(X)-H(X_0,X_\Sigma)-\frac{1}{2}\sum_{i=1}^{d-1}H(X_0,X_i)\right)\right),\]
where $\frac{1}{b(n)}\le c_1(m),c_2(m)\le b(n)$ for some polynomial $b(n)$. 

Since $X_d=X_\Sigma-\sum_{i=1}^{d-1} X_i$ and $X_\Sigma=\sum_{i=1}^d X_i$,  applying the last statement  of Lemma \ref{entl} twice,  we get that
\begin{align}\label{mutdet}
H(X)=(X_0,X_1,\dots, X_{d})&=H(X_0,X_1,\dots, X_{d},X_\Sigma)\\&=H(X_0,X_1,\dots, X_{d-1},X_\Sigma).\nonumber
\end{align} 

Combining this with Lemma \ref{shaererszeru}, we get that 

\begin{align*}
H(X)&=H(X_0,...,X_{d-1},X_\Sigma)\\&\le \sum_{i=1}^{d-1} H(X_0,X_i)+H(X_0,X_\Sigma)-(d-1)H(X_0). 
\end{align*}

Or more generally, for every $i=1,2,\dots, d$, we have

\[H(X)\le \sum_{\substack{1\le j\le d\\j\neq i}} H(X_0,X_j)+H(X_0,X_\Sigma)-(d-1)H(X_0).\] 

Summing up these inequalities for $i=1,2,...,d-1$, we get that

\begin{multline}\label{Hegy}
(d -1)H(X)\\\le 
(d-2)\sum_{i=1}^{d-1} H(X_0,X_i)+(d-1)H(X_0,X_d)+(d-1)H(X_0,X_\Sigma)-(d-1)^2 H(X_0).
\end{multline} 

Note that $X_0,X_1,\dots, X_d$ all have the same distribution, so $H(X_0)=H(X_1)=\dots=H(X_d)$. Combining this with Equation \eqref{mutdet} and Inequality \eqref{szubad}, we have
\begin{align}\label{Hketto}
H(X)&=H(X_0,...,X_{d-1},X_\Sigma)\\&\le H(X_0,X_\Sigma)+\sum_{i=1}^{d-1} H(X_i)=H(X_0,X_\Sigma)+(d-1)H(X_0).\nonumber
\end{align}

Therefore,
\begin{align*}
H(X_0)+H(X)&-H(X_0,X_\Sigma)-\frac{1}2\sum_{i=1}^d H(X_0,X_i)\\&= H(X_0)+H(X)-H(X_0,X_\Sigma)-\frac{1}{2(d-1)}\sum_{i=1}^{d-1} H(X_0,X_i)\\&\qquad\qquad-\frac{1}2 \left(\frac{d-2}{d-1} \sum_{i=1}^{d-1} H(X_0,X_i)+H(X_0,X_d)\right)\\
& \le H(X_0)+H(X)-H(X_0,X_\Sigma)-\frac{1}{2(d-1)}\sum_{i=1}^{d-1} H(X_0,X_i)\\&\qquad\qquad-\frac{1}2 \left(H(X)+(d-1)H(X_0)-H(X_0,X_\Sigma)\right)\\&= \frac{1}{d-1}\left(H(X)-H(X_0,X_\Sigma)-\frac{1}2\sum_{i=1}^{d-1}H(X_0,X_i)\right)\\&\qquad\qquad+\frac{d-3}{2(d-1)}(H(X)-H(X_0,X_\Sigma))-\frac{(d-3)}{2}H(X_0)\\& \le \frac{1}{d-1}\left(H(X)-H(X_0,X_\Sigma)-\frac{1}2\sum_{i=1}^{d-1}H(X_0,X_i)\right),
\end{align*}
where at the first inequality, we used Inequality \eqref{Hegy}, and at the second inequality, we used Inequality \eqref{Hketto}. This gives the first statement. To get the second one, observe that
\begin{align*}
|S(q)|\mathbb{P}(C_n^{(d)}q=r)&=\sum_{m\in\mathcal{M}^S(q,r)}E(m)\le\sum_{m\in\mathcal{M}^S(q,r)}f(n)p(m)^{\frac{1}{d-1}}\\&\le |\mathcal{M}^S(q,r)|f(n)\mathbb{P}(r-C_n^{(d-1)}q\sim q)^{\frac{1}{p-1}}.
\end{align*}
\end{proof}
Now we prove Lemma \ref{elhalasztott}.
\begin{proof}
Clearly we may assume that $h=2$. 
The size of $\mathcal{M}^S(q,r)$ is polynomial in $n$, so it is enough to prove that for a fixed $m\in \mathcal{M}^S(q,r)$, we have a good upper bound on $\mathbb{P}(\Sigma(Q)=r\text{ and }m_Q=m)$. To show this, let $X=(X_0,X_1,X_2)\in V^{1+2}$ be a random variable, such that $\mathbb{P}(X=t)=\frac{m(t)}{n}$ for every $t\in V^{1+2}$, and let $X_\Sigma=X_1+X_2$. Then $\mathbb{P}(\Sigma(Q)=r\text{ and }m_Q=m)$ can be upper bounded by some polynomial multiple of
\begin{align*}
\exp&\left(n\left(H(X)-H(X_0,X_\Sigma)-\frac{1}2\left(H(X_0,X_1)+H(X_0,X_2)\right)\right)\right)\\&=
\exp\Big(n\big(-H(X_0)-\frac{1}{2}((H(X_0,X_1)+H(X_0,X_\Sigma)-H(X)-H(X_0))\\&\qquad\quad+(H(X_0,X_2)+H(X_0,X_\Sigma)-H(X)-H(X_0)))\big)\Big)\\&
\le \exp(-nH(X_0))\le |V|^{-n}\exp(Cn^{2\alpha-1}),
\end{align*}
using the fact that for $i\in\{1,2\}$, we have \[H(X_0,X_i)+H(X_0,X_\Sigma)\ge H(X_0)+H(X_0,X_i,X_\Sigma)=H(X_0)+H(X),\]
which is a combination of Inequality \eqref{szubmod} and the last statement of Lemma \ref{entl}.
\end{proof}

For any non-negative integral vector $m$ indexed by $V^{1+d}$ and for $i \in \{1,2,\dots,d\}$, we define
\[
E_0(m)=\frac{m(V^{1+d})!}{\prod_{c\in V} m(\tau_0=c)!} \frac{\prod_{a,b\in V} m(\tau_0=a,\tau_\Sigma=b)!}{\prod_{t\in V^{1+d}} m(t)!},
\]
and
\begin{multline*}
E_i(m)=\left(\frac{m(V^{1+d})!}{2^{m(V^{1+d})/2} (m(V^{1+d})/2)!}\right)^{-1}  \\\times \left(\prod_{a\in V} \frac{m(\tau_i=a,\tau_0=a)!}{2^{m(\tau_i=a,\tau_0=a)/2} (m(\tau_i=a,\tau_0=a)/2)!}\right)\left(\prod_{a\neq b\in V} \sqrt{m(\tau_0=a,\tau_i=b)!}\right).
\end{multline*}
Finally, let
\[E(m)=E_0(m)\prod_{i=1}^d E_i(m).\]

Here we need to define $(\ell+\frac{1}{2})!$ for an integer $\ell$. The simple definition \break $(\ell+\frac{1}{2})!=\ell!\sqrt{\ell+1}$ is good enough for our purposes.

Recall that for $q,r\in V^n$ and $m\in \mathcal{M}^S(q,r)$, we already defined $E(m)$ as
\[E(m)=|S(q)|\mathbb{P}(m_{\bar{Q}}=m\text{ and }\Sigma(\bar{Q})=r),\]
where $\bar{Q}$ is a random $(q,1,d)$-tuple.

Using Equation \eqref{eprk}, it is straightforward to verify that for a special $m$ like above, the two definitions coincide.

Given a $q\in V^n$, a $(q,1,d)$-tuple $Q$ or $m_Q$ itself will be called $W$-decent if for any $u\in W^{1+d}$ we have \[\frac{1+m_Q(\tau_0=u_0,\tau_\Sigma=u_\Sigma)}{1+m_Q(u)}\le \log^2 n.\] 

A non-negative integral vector $m$ indexed by $V^{1+d}$ will be called $W$-half-decent if for every $u\in W^{1+d}$, we have 
\[\frac{1+ m(\tau_0=u_0,\tau_\Sigma=u_\Sigma)}{1+m(u)}\le \log^4 n,\]
and for every $c\in W$, we have 
\[\left|m(\tau_0=c)-\frac{n}{|W|}\right|<2n^{\alpha},\]
 where $n=\sum_{t\in V^{1+d}} m(t)$.

\begin{lemma}[The analogue of Lemma \ref{csakdecent}]
For any coset $W\in\Cosets(V)$, we have
\begin{multline*}
\limsup_{n\to\infty} \sum_{q\in D_W^n} |S(q)|\mathbb{P}(\Sigma(\bar{Q})=r_q)\\=\limsup_{n\to\infty} \sum_{q\in D_W^n} |S(q)|\mathbb{P}(\Sigma(\bar{Q})=r_q\text{ and }\bar{Q}\text{ is }W-\text{decent}).
\end{multline*}
\end{lemma}
\begin{proof}

As in the proof of Lemma \ref{csakdecent}, it is enough to show that \[|S(q)|\mathbb{P}(\Sigma(\bar{Q})=r_q\text{ and }\bar{Q}\text{ is not  }W-\text{decent})<n^{-(|V|+1)}\] for every $(W,C\log n)$-typical vector $q\in V^n$ if $n$ is large enough. 

Consider a $(W,C\log n)$-typical vector $q\in V^n$, and let 
\begin{multline*}
\mathcal{M}_B=\\ \{m_Q|\quad Q\text{ is a not $W$-decent }(q,1,d)\text{-tuple, such that }\Sigma(Q)=r_q\}\subset\mathcal{M}^S(q,r_q).
\end{multline*}

Recall that for $m\in \mathcal{M}^S(q,r_q)$, we defined $p(m)$ as the probability of the event that for a random $(q,1,d-1)$-tuple $\bar{Q}=(\bar{q}^{(0)},\bar{q}^{(1)},\dots,\bar{q}^{(d-1)})$, we have that \break $r_q-\Sigma(\bar{Q})$ is a matching permutation of $q$ and the $(q,1,d)$-tuple \break $Q'=(\bar{q}^{(0)},\bar{q}^{(1)},\dots, \bar{q}^{(d-1)}, r_q-\Sigma(\bar{Q}))$ satisfies $m_{Q'}=m$.

Note that for $m\in \mathcal{M}_B$ the event above is contained in the  event that

\[
\text{there is a }t\in W^{1+(d-1)}\text{ and }c\in dW\text{ such that }\qquad\qquad\qquad\qquad\]
\begin{equation*}
\frac{1+|\{i|\quad r_q(i)=c\text{ and }q_i=t_0\}|}{1+|\{i|\quad r_q(i)=c\text{ and }\bar{Q}(i)=t\}|}> \log^2 n.   
\end{equation*}

Let $p'(q)$ be the probability of the latter event. As we just observed, $p(m)\le p'(q)$ for all $m\in \mathcal{M}_B$. Using Lemma \ref{Eesp} and Lemma \ref{bazumaM}, we obtain
\begin{align*}
|S(q)|\mathbb{P}(\Sigma(\bar{Q})=r_q\text{ and }\bar{Q}\text{ is not  }W-\text{decent})&=\sum_{n\in \mathcal{M}_B} E(m)\\&\le \sum_{n\in \mathcal{M}_B} f(n)p(m)^{\frac{1}{d-1}}\\&\le |\mathcal{M}_B|f(n)p'(q)^{\frac{1}{d-1}}<n^{-(|V|+1)}
\end{align*}
for large enough $n$.
\end{proof}

Let \[\mathcal{M}^S=\{m_Q\quad|\quad Q\text{ is a }(q,1,d)\text{-tuple for some }n\ge 0\text{ and }q\in V^n\}.\]

\begin{lemma}[The analogue of Lemma \ref{Hilbert}]
There are finitely many vectors $m_1,m_2,...,m_{\ell}\in \mathcal{M}^S$, such that
\[
\mathcal{M}^S=\{c_1 m_1+c_2m_2+\dots+ c_\ell m_\ell \quad |\quad c_1,c_2,\dots, c_{\ell}\text{ are non-negative integers}\}.
\]
\end{lemma}
\begin{proof}
We define \[\mathcal{R}=\left\{(m,g)\quad|\quad m\in\mathbb{R}^{V^{1+d}},\quad g\in \mathbb{R}^{\{1,2,\dots,d\}\times V}\right\}.\]
Consider the linear subspace $\mathcal{R}'$ of $\mathcal{R}$ consisting of pairs $(m,g)$ satisfying the following liner equations:
\[m(\tau_0=a\text{ and }\tau_i=b)=m(\tau_0=b\text{ and }\tau_i=a)\]
for all $a,b\in V$ and $i\in\{1,2,\dots, d\}$, moreover,
\[m(\tau_0=c\text{ and }\tau_i=c)=2g(i,c)\]
for all $c\in V$ and $i\in\{1,2,\dots,d\}$.

Let $\mathcal{M}_0$ be the set of non-negative integral points of $\mathcal{R}'$. Observe that $\mathcal{M}_0$ consists of the integral points of a rational polyhedral cone. From  \cite[Theorem~16.4]{schr}, we know that this cone is generated by an integral Hilbert basis, i. e., there are finitely many vectors $(m_1,g_1),(m_2,g_2),...,(m_{\ell},g_\ell)\in \mathcal{M}_0$, such that
\[
\mathcal{M}_0=\{c_1 \cdot(m_1,g_1)+\dots+ c_\ell\cdot (m_\ell,g_\ell)  |\quad c_1,c_2,\dots, c_{\ell}\text{ are non-negative integers}\}.
\]

Then the vectors $m_1,m_2,\dots,m_\ell\in \mathcal{M}^S$ have the required properties.

Note we only introduced the extra component $g$ to enforce the parity constraints in \eqref{equ3k}.
\end{proof}
As before, we may assume that the indices in the lemma above are chosen such that there is an $h$ such that the supports of $m_1,m_2,\dots,m_h$ are contained in $W^{1+d}$, and the supports of $m_{h+1},m_{h+2},...,m_\ell$ are not contained in $W^{1+d}$.

\begin{lemma}[The analogue of Lemma \ref{segedlem}]
Consider a non-negative integral $W$-half-decent vector $m_0\in\mathbb{R}^{V^{1+d}}$, such that  $\|m_0\|_{W^C}=m(t\not \in W^{1+d})=O(\log n)$, where $n=\sum_{t\in V^{1+d}} m(t)$. For $u\in V^{1+d}$, let $\chi_u\in\mathbb{R}^{V^{1+d}}$ be such that $\chi_u(u)=1$ and $\chi_u(t)=0$ for every $t\neq u\in V^{1+d}$.
\begin{itemize}
\item If $u\in W^{1+d}$, then $E(m_0+\chi_u)/E(m_0)=O\left( \log^4 n\right)$;
\item If $u_0\not\in W$, then $E(m_0+\chi_u)/E(m_0)=O\left(\frac{\log^{d+1} n}{n^{d/2-1}}\right)$;
\item If $u_0\in W$ and $u\not\in W^{1+d}$, then $E(m_0+\chi_u)/E(m_0)=O\left(\log^2 n\right).$
\end{itemize}
\end{lemma}
\begin{proof}

Let 
\begin{align*}
g&=\frac{1+ m_0(\tau_0=u_0,\tau_\Sigma=u_\Sigma)}{1+m_0(u)},\\ h&=\frac{n+1}{m(\tau_0=u_0)+1}, \qquad\qquad\text{ and}\\f_i&=\sqrt{\frac{1+ m_0(\tau_0=u_0,\tau_i=u_i)} {n+1}}.
\end{align*}

\begin{lemma} 
\[E(m_0+\chi_u)/E(m_0)=O(g\cdot h\cdot\prod_{i=1}^d f_i).\]
\end{lemma}
\begin{proof}
It is straightforward to check that $E_0(m_0+\chi_u)/E_0(m_0)=g\cdot h$. Let \break $i\in \{1,2,\dots,d\}$. First assume that $u_i\neq u_0$, then
\[E_i(m_0+\chi_u)/E_i(m_0)=\frac{\sqrt{2}}{n+1} \cdot\frac{\left(\frac{n+1}{2}\right)!}{\left(\frac{n}{2}\right)!}\cdot \sqrt{m_0(\tau_i=u_i,\tau_0=u_0)+1}.\]
Recall that for any integer $\ell$ we defined $(\ell+\frac{1}{2})!$ as $(\ell+\frac{1}{2})!=\ell!\sqrt{\ell+1}$. Thus, if $n$ is even, then \[\frac{\left(\frac{n+1}{2}\right)!}{\left(\frac{n}{2}\right)!}=\sqrt{\frac{n}{2}+1}=O(\sqrt{n+1}),\] 
and if $n$ is odd, then
\[\frac{\left(\frac{n+1}{2}\right)!}{\left(\frac{n}{2}\right)!}=\sqrt{\frac{n+1}{2}}=O(\sqrt{n+1}).\]
Therefore, $E_i(m_0+\chi_u)/E_i(m_0)=O(f_i)$. In the case $u_i=u_0=c$,  we have 
\[E_0(m_0+\chi_u)/E_0(m_0)=\frac{\sqrt{2}}{n+1} \cdot\frac{\left(\frac{n+1}{2}\right)!}{\left(\frac{n}{2}\right)!}\cdot \frac{m_0(\tau_i=c,\tau_0=c)+1}{\sqrt{2}}
\cdot \frac{\left(\frac{m_0(\tau_i=c,\tau_0=c)}{2}\right)!}{\left(\frac{m_0(\tau_i=c,\tau_0=c)+1}{2}\right)!}.
\]
A similar argument as above gives that $E_i(m_0+\chi_u)/E_i(m_0)=O(f_i)$ also holds in this case.
 The statement follows from the fact that 
\[E(m_0+\chi_u)/E(m_0)=\prod_{i=0}^d E_i(m_0+\chi_u)/E_i(m_0).\]
\end{proof}

If $u\in W^{1+d}$, then since $m_0$ is $W$-half-decent, we have $g\le\log^4 n$, $h=O(1)$ and clearly $f_i\le 1$, thus the statement follows.

If $u_0\not\in W$, then $g=O(\log n)$, $h=O(n)$, $f_i=O(\frac{\log n}{\sqrt{n}})$, and the statement follows.

If $u_0\in W$ and $u\not\in W^{1+d}$, then we consider two cases: 
\begin{enumerate}
\item If $u_\Sigma\in dW$, then $g=O(n)$, $h=O(1)$, moreover there are at least two indices $i$ such that $u_i\not\in W$. For such an $i$, we have $f_i=O(\frac{\log n}{\sqrt{n}})$,  otherwise we have $f_i\le 1$, from these the statement follows. 
\item
If $u_\Sigma\not\in dW$, then $g=O(\log n)$, $h=O(1)$ and $f_i\le 1$ for every $i$. The statement follows.
\end{enumerate}
\end{proof}

The previous lemma has the following consequence.

\begin{lemma}[The analogue of Lemma \ref{DeltaMbecsElo}]\label{DBecsS}
There are $D,\delta>0$, such that for any $i\in \{h+1,h+2,\dots,\ell\}$ and any  non-negative integral $W$-half-decent vector $m_0\in\mathbb{R}^{V^{1+d}}$, such that $\|m_0\|_{W^C}=O(\log n)$,  we have
 \[E(m_0+m_i)/E(m_0)=O\left(\left(n^{-\delta}\log^D n\right)^{\|m_i\|_{W^C}}\right).\]
\end{lemma}

\begin{proof}
Take any $i\in\{h+1,h+2,\dots,\ell\}$. Since $m_i$ is not supported on $W^{1+d}$,\break we have a $u\not\in W^{1+d}$ such that $m_i(u)\ge 1$. If $u_0\not \in W$, then \break $m_i(\tau_0\not\in W)\ge m_i(\tau_0=u_0)\ge 1$. If $u_0\in W$, then there is a $j$ such that $u_j\not\in W$, thus
\[m_i(\tau_0\not\in W)\ge m_i(\tau_0=u_j,\tau_j=u_0)=m_i(\tau_0=u_0,\tau_j=u_j)\ge m_i(u)\ge 1.\]
In both cases, we obtained that $m_i(\tau_0\not\in W)\ge 1$. Note that for $d\ge 3$, we have $d/2-1>0$. From the previous statements and Lemma \ref{DBecsS}, it follows that for a large enough $D$ and a small enough $\delta>0$, we have

 \[E(m_0+m_i)/E(m_0)=O\left(\left(\log^D n\right)^{\|m_i\|_{W^C}}n^{-(d/2-1)}\right)=O\left(\left(n^{-\delta}\log^D n\right)^{\|m_i\|_{W^C}}\right).\]

\end{proof}

With these modifications above, we proved Theorem \ref{FoFo22}.

As an easy consequence of Theorem \ref{FoFo22} we obtain following analogue of Corollary \ref{correduced}.  
The random $(n-1)\times (n-1)$ matrix $C_n'$ is obtained from $C_n$ by deleting its last row and last column. Recall $q\in V^{n-1}$ the subgroup generated by $q_1,q_2,\dots,q_{n-1}$ is denoted by $G_q$. Let $U_q^S$ be a uniform random element of the set \break $\{w\in G_q^{n-1}|\quad <q\otimes w>\in I_2\}$. 

\begin{corr}\label{correduced22}
We have
\[\pushQED{\qed} 
\lim_{n\to\infty}  \sum_{q\in V^{n-1}} d_\infty(C_n'q,U_{q}^S)=0.\qedhere
\popQED\]
\end{corr}

Note that for $q\in V^{n-1}$ such that $G_q=V$, if $r\in V^{n-1}$ and $<q\otimes r>\in I_2$ then $\mathbb{P}(U_q^S=r)=|V|^{-(n-1)}2^{\rank2(V)} |\wedge^2 V|$. Therefore, Theorem \ref{momentumokUD} can be proved using the following observation.
\begin{lemma}
If $d$ is even, then $<q\otimes dq>\in I_2$ for every $q\in V^{n-1}$. If $d$ is odd, then $<q\otimes dq>\in I_2$ if and only if $s(q)$ is an element of the subgroup $V'=\{2v|v\in V\}$. The subgroup $V'$ has index $2^{\rank2(V)}$ in $V$. \qed
\end{lemma} 

For odd $d$, Theorem \ref{CohenlenstraUD} follows from Theorem \ref{momentumokUD} and Theorem \ref{momentproblemspec} part (2).

\section{The $2$-Sylow subgroup in the case of even $d$}\label{parosparos}


Assume that $d$ is even. Let $\Delta_n$ be the reduced Laplacian of $H_n$, and $\Gamma_n$ be the corresponding sandpile group. Theorem \ref{momentumokUD} provides us the limit of the surjective $V$-moments of $\Gamma_n$. However, these moments grow too fast, so Theorem \ref{woodmoment} can not be applied to get the existence of a limit distribution. We can overcome this difficulty by using that $\Gamma_n$ has a special property given in the next lemma.
\begin{lemma}\label{paratlanrang}
The group $\Gamma_n\otimes \mathbb{Z}/2\mathbb{Z}$ has odd rank.
\end{lemma}
Given any integral matrix $M$, let $\overline{M}$ be its mod 2 reduction. That is, $\overline{M}$ is a matrix over the $2$ element field, where an entry is $1$ if and only if the corresponding entry of $M$ is odd.

\begin{prop}
Let $M$ be a integral $m\times m$ matrix. Then
\[\rankk(\cok(M)\otimes \mathbb{Z}/2\mathbb{Z})=\dim\ker\overline{M}=m-\rankk(\overline{M}).\]
\end{prop} 
\begin{proof}
It is straightforward to verify the statement if $M$ is diagonal. If $M$ is not diagonal, then $M$ can be written as $M=ADB$, where $D$ is diagonal, and \break $A,B\in \text{GL}_m(\mathbb{Z})$. This is the so-called Smith normal form. The statement follows from the fact that $\dim\ker\overline{M}=\dim\ker\overline{ADB}=\dim\ker\overline{A}\cdot\overline{D}\cdot\overline{B}=\dim\ker\overline{D}$, and $\cok{M}=\cok{ADB}=\cok{D}$.
\end{proof}

\begin{proof}(Lemma \ref{paratlanrang})
Observe that $\overline{\Delta_n}$ is a symmetric matrix, where all the diagonal entries are $0$. Such a matrix alway has even rank. See for example \cite[Theorem 3]{macw}. Recall that $\Delta_n$ is an $(n-1)\times (n-1)$ matrix, where $n$ is even. Thus, the statement follows from the previous proposition.   
\end{proof}

In the first part of this section, we prove a modified version of Theorem \ref{woodmoment}, which allows us to make use of the fact that $\Gamma_n\otimes \mathbb{Z}/2\mathbb{Z}$ has odd rank. For most of the proof we can follow the original argument of Wood \cite{wood} almost word by word with only small modifications. A few proofs are deferred to the Appendix, since they are almost identical to the proofs of Wood \cite{wood}.

 We start by giving a few definitions. A partition $\lambda$ of length $m$ is a  sequence $\lambda_1\ge\lambda_2\ge\dots\ge \lambda_m\ge 1$ of positive integers. It will be a convenient notation to also  define $\lambda_i=0$ for $i>m$.   The transpose partition $\lambda'$  of $\lambda$ is defined by setting $\lambda'_j$ to be the number of $\lambda_i$ that are at least $j$. Thus, the length of $\lambda'$ is  $\lambda_1$. Recall that any finite abelian $p$-group $G$ is isomorphic to
 \[\bigoplus_{i=1}^m \mathbb{Z}/p^{\lambda_i}\mathbb{Z}\]
 for some partition $\lambda$ of length $m$. We call $\lambda$ the type of the group $G$.  In fact, this provides a bijection between the set of isomorphism classes of finite abelian $p$-groups and the set of partitions.

\if \rver 1

\begin{lemma} \label{L:Hacts}\hfill
\begin{enumerate}
\item

Given a positive integer $m$, and $b\in \mathbb{Z}^m$ such that $b_1$ is odd, $b_1\ge b_2\ge \dots \ge b_m$,
we have an entire analytic function in the $m$ variables $z_1,\dots,z_m$
$$
H_{m,2,b}(z)=\sum_{\substack{d_1,\dots,d_m \geq 0\\d_2+\cdots + d_m \leq b_1 }} a_{d_1,\dots,d_m} z_1^{d_1} \cdots z_m^{d_m}
$$
and a constant $E$ such that
$$
a_{d_1,\dots,d_m}\leq E 2^{-b_1d_1 - d_1(d_1+1)}.
$$
Further, 
if $f$ is a partition of length $\le m$ such that  $f > b$ (in the lexicographic ordering), $f_1$ is odd, then  ${H}_{m,2,b}(2^{f_1}, 2^{f_1+f_2}, \dots, 2^{f_1+\dots+f_m})=0$.  If $f=b$, then 
${H}_{m,2,b}(2^{f_1}, 2^{f_1+f_2}, \dots, 2^{f_1+\dots+f_m})\ne 0.$

\item

Given a positive integer $m$, a prime $p>2$,\footnote{In fact, this statement is also true for $p=2$, but we will not use this.} and $b\in \Z^m$ with $b_1\geq b_2\geq \dots \geq b_m$,
we have an entire analytic function in the $m$ variables $z_1,\dots,z_m$
$$
H_{m,p,b}(z)=\sum_{\substack{d_1,\dots,d_m \geq 0\\d_2+\cdots + d_m \leq b_1 }} a_{d_1,\dots,d_m} z_1^{d_1} \cdots z_m^{d_m}
$$
and a constant $E$ such that
$$
a_{d_1,\dots,d_m}\leq E p^{-b_1d_1 - \frac{d_1(d_1+1)}{2}}.
$$
Further, 
if $f$ is a partition of length $\leq m$ and  $f > b$ (in the lexicographic ordering), then  $H_{m,p,b}(p^{f_1}, p^{f_1+f_2}, \dots, p^{f_1+\dots+f_m})=0$.  If $f=b$, then 
$H_{m,p,b}(p^{f_1}, p^{f_1+f_2}, \dots, p^{f_1+\dots+f_m})\ne0.$
\end{enumerate}
\end{lemma}
\begin{proof}
See the Appendix for the proof.
\end{proof}

In the original proof of Wood \cite{wood}, the prime $2$ was not handled separately. That is, the functions given in part (2) of Lemma \ref{L:Hacts} were used for all primes. Let us restrict our attention to random groups $G$ where $G\otimes \mathbb{Z}/2\mathbb{Z}$ has odd rank. Then, for the prime $2$, we can use the functions given in part (1) of Lemma \ref{L:Hacts} instead of the ones given in part (2), and still proceed with the proof, as we show in the next lemmas. Note that part (1) provides better bounds for the coefficients. This allows us to handle faster growing moments.

\begin{theorem}\label{T:Momdet}
Let $2=p_1,\dots,p_s$ be distinct primes.  Let $m_1,\dots,m_s\geq 1 $ be integers.

Let $M_j$ be the set of partitions $\lambda$ at most $m_j$ parts.  
Let $M=\prod_{j=1}^{s} M_j$.  For $\mu\in M$, we write $\mu^j$ for its $j$th entry, which is a partition consisting of non-negative integers $\mu^j_i$ with $\mu^j_1\geq\mu^j_2\geq\dots \mu^j_{m_j}$.
Let
\[M_0=\{\mu\in M\quad|\quad \mu^1_1\text{ is odd}\}.\]
 Suppose we have non-negative reals $x_\mu, y_\mu$, for each tuple of partitions $\mu\in M_0$.
Further suppose that we have non-negative reals $C_\lambda$ for each $\lambda\in M$ such that
$$C_\lambda\leq 2^{\lambda_1^1}\prod_{j=1}^{s} {F^{m_j}p_j^{\sum_i \frac{\lambda^{j}_i(\lambda^{j}_i-1)}{2}}},$$
where $F>0$ is an absolute constant.
Suppose that 
for all $\lambda\in M$,
\begin{equation}\label{E:infsys}
\sum_{\mu\in M_0} x_\mu \prod_{j=1}^{s}  p_j^{\sum_i \lambda^{j}_i \mu^{j}_i} =\sum_{\mu\in M_0} y_\mu \prod_{j=1}^{s} p_j^{\sum_i \lambda^{j}_i \mu^{j}_i} =C_\lambda.
\end{equation}
Then for all $\mu\in M_0$, we have that $x_\mu=y_\mu$.
\end{theorem}
\begin{proof}
See the Appendix for the proof.
\end{proof}

\begin{lemma}\label{L:BoundMom}
There is a constant $F$, such that for any finite abelian $p$-group $G$ of type $\lambda$, we have
\[\sum_{G_1\text{ subgroup of G}}|\wedge^2 G_1|\le F^{\lambda_1}p^{\sum_i \frac{\lambda_i'(\lambda_i'-1)}{2}}.\]
Moreover, if $G$ finite abelian $2$-group $G$ of type $\lambda$, we have
\[\sum_{G_1\text{ subgroup of G}}2^{\rank2(G_1)}|\wedge^2 G_1|\le F^{\lambda_1}2^{\lambda_1'+\sum_i \frac{\lambda_i'(\lambda_i'-1)}{2}}.\]
\end{lemma}
\begin{proof}
The first statement is the same as \cite[Lemma 7.5]{wood}.\footnote{In the latest arxiv version of this paper this is Lemma 7.4} The second statement follows from first by using the elementary fact that for any subgroup $G_1$ of $G$, we have $\rank2(G_1)\le \rank2(G)=\lambda_1'$. 
\end{proof}

\begin{lemma}(\cite[Lemma 7.1]{wood})\label{HomSize}
Let $G_\mu$ and $G_\lambda$ be two finite abelian $p$-groups of type $\mu$ and $\lambda$. Then
\[|\Hom(G_\mu,G_\lambda)|=p^{\sum_i \mu_i'\lambda_i'}.\]
\end{lemma}

\begin{theorem}\label{T:MomDetDetail}
Let $X_n$ be a sequence of random variables taking values in finitely generated abelian groups.
Let $a$ be an even positive integer and $A$ be the set of (isomorphism classes of) abelian groups with exponent dividing $a$. Assume that $\rankk(X_n\otimes \mathbb{Z}/2\mathbb{Z})$ is odd with probability $1$ for every $n$.
Suppose that for every $G\in A$, we have
$$
\lim_{n\ra \infty} \mathbb{E}| \Sur(X_n, G)| = 2^{\rank2(G)}|\wedge^2 G|. 
$$
Then for every $H\in A$,
 the limit
$
\lim_{n\ra\infty} \mathbb{P}(X_n\tensor \Z/a\Z \isom H)
$
exists, and for all $G\in A$, we have
$$
\sum_{H\in A} \lim_{n\ra\infty} \P(X_n\tensor \Z/a\Z \isom H) |\Sur(H,G)|=2^{\rank2(G)}|\wedge^2 G|.
$$

 Suppose $Y_n$ is a sequence of random variables taking values in finitely generated abelian groups
 such that $\rankk(Y_n\otimes \mathbb{Z}/2\mathbb{Z})$ is odd with probability $1$ for every $n$, and for every $G\in A$, we have
$$
\lim_{n\ra \infty} \mathbb{E}| \Sur(Y_n, G)| = 2^{\rank2(G)}|\wedge^2 G|. 
$$
Then, we have that for every every $H\in A$
$$
\lim_{n\ra\infty} \P(X_n\tensor \Z/a\Z \isom H) =\lim_{n\ra\infty} \P(Y_n\tensor \Z/a\Z \isom H).
$$
\end{theorem}
\begin{proof}
See the Appendix for the proof.
\end{proof}

\else


For a measure $\nu$ on the finitely generated abelian groups and a finite abelian group $V$, we define $\Sur(\nu,V)=\sum_{G}\nu(G)|\Sur(G,V)|$.

\begin{lemma}\label{wool}

\begin{enumerate}[(i)]
\item \label{wooi}
Let $k$ be a positive integer. Then there is a unique probability measure $\nu_k$ on $FA_{odd}(2,k)$, such that $\Sur(\nu_k,V)=2^{\rank2(V)}|\wedge^2 V|$ for any $V\in FA(2,k)$. Moreover, if $X_n$ is a sequence of random finitely generated abelian groups of odd rank, such that for any $V\in FA(2,k)$ we have 
\[\lim_{n\to\infty}\mathbb{E}|\Sur(X_n,V)|=2^{\rank2(V)}|\wedge^2 V|,\]
then for every $V\in FA_{odd}(2,k)$
\[\lim_{n\to\infty} P(X_n\otimes \mathbb{Z}/2^k\mathbb{Z}\isom V)=\nu_k(V).\]

\item \label{wooii} There is a unique probability measure $\nu$ on the set of finite abelian $2$-groups of odd rank, such that $\Sur(\nu,V)=2^{\rank2(V)}|\wedge^2 V|$ for any finite abelian $2$-group $V$. Moreover, if $X_n$ is a sequence of random finite abelian $2$-groups of odd rank\footnote{Or slightly more generally $X_n$ has odd rank and it is a direct sum of a finite abelian 2-group and a  finitely generated free abelian group.}, such that for any finite abelian $2$-group $V$ we have 
\[\lim_{n\to\infty}\mathbb{E}|\Sur(X_n,V)|=2^{\rank2(V)}|\wedge^2 V|,\]
then for every finite abelian $2$-group $V$ of odd rank we have
\[\lim_{n\to\infty} \mathbb{P}(X_n\isom V)=\nu(V).\]
\end{enumerate}
\end{lemma}
\begin{proof}

 The statements of (\ref{wooi}) can be obtained by slightly modifying the argument of \break Wood \cite{wood}[Theorem 8.3] by making use of the fact that the ranks of the groups $X_n$ are odd. We only point out the details   needed to be changed. Lemma 8.1. should be replaced by the following lemma.
\begin{lemma}\label{lemma8p1}
Given a positive integer $m$, and $b\in \mathbb{Z}^m$ such that $b_1$ is odd, $b_1\ge b_2\ge \dots \ge b_m$,
we have an entire analytic function in the $m$ variables $z_1,\dots,z_m$
$$
\bar{H}_{m,2,b}(z)=\sum_{\substack{d_1,\dots,d_m \geq 0\\d_2+\cdots + d_m \leq b_1 }} a_{d_1,\dots,d_m} z_1^{d_1} \cdots z_m^{d_m}
$$
and a constant $E$ such that
$$
a_{d_1,\dots,d_m}\leq E 2^{-b_1d_1 - d_1(d_1+1)}.
$$
Further, 
if $f$ is a partition of length $\le m$ such that  $f > b$ (in the lexicographic ordering) and $f_1$ is odd, then  $\bar{H}_{m,2,b}(2^{f_1}, 2^{f_1+f_2}, \dots, 2^{f_1+\dots+f_m})=0$.  If $f=b$, then 
$\bar{H}_{m,2,b}(2^{f_1}, 2^{f_1+f_2}, \dots, 2^{f_1+\dots+f_m})\ne 0.$
\end{lemma}
\begin{proof}
The proof is the same as the proof of  \cite{wood}[Lemma 8.1]. But instead of $G(z_1)$ we use 
\[\bar{G}_1=\prod_{\substack{j>b_1\\j\text{ is odd}}} \left(1-\frac{z_1}{2^j}\right)=\sum_{d_1\ge 0}c_{d_1} z_1^{d_1}.\]
Observe that $\bar{G}_1(4z)=(1-\frac{z}{2^{b_1}})\bar{G}_1(z)$. So $4^n c_n=c_n-2^{-b_1}c_{n-1}$, or equivalently $c_n=\frac{-2^{-b_1-2n}}{1-4^{-n}}c_{n-1}$. Since $c_0=1$, by induction we obtain that \[c_n=\frac{(-1)^n 2^{-n b_1-n(n+1)}}{\prod_{i=1}^n (1-4^{-i})}.\]
So $|c_n|\le 2^{-n b_1-n(n+1)} \prod_{i=1}^{\infty} (1-4^{-i})^{-1}$.
\end{proof}  

Lemma 7.4. of Wood \cite{wood} states that there is a constant $F$ such that for any abelian 2-group $G$ of type $\lambda$ we have
\[\sum_{G_1\text{ subgroup of }G} |\wedge^2 G_1|\le F^{\lambda_1} 2^{\sum_i\frac{\lambda_i'(\lambda_i'-1)}{2}}.\]
As a simple corollary of this we obtain that
\begin{equation}
\sum_{G_1\text{ subgroup of }G} |\wedge^2 G_1|2^{\rank2(G_1)}\le F^{\lambda_1} 2^{\rank2(G)} 2^{\sum_i\frac{\lambda_i'(\lambda_i'-1)}{2}}=F^{\lambda_1}  2^{\lambda_1'+\sum_i\frac{\lambda_i'(\lambda_i'-1)}{2}}.
\end{equation}

Instead of Theorem 8.2. we use the following lemma.
\begin{lemma}
Let $m\ge 1$ be an integer.
Let $M_0$ be the set of partitions $\lambda$ with at most $m$ parts. 
Let $M$ be the set of partitions $\lambda\in M_0$ such that $\lambda_1$ is odd.  

 Suppose we have non-negative reals $x_\mu, y_\mu$, for each  partition $\mu\in M$.
Further suppose that we have non-negative reals $C_\lambda$ for each $\lambda\in M_0$ such that
$$C_\lambda\le F^m  2^{\lambda_1+\sum_i\frac{\lambda_i(\lambda_i-1)}{2}},$$
where $F>0$ is an absolute constant.
Suppose that 
for all $\lambda\in M_0$,
\begin{equation}
\sum_{\mu\in M} x_\mu   2^{\sum_i \lambda_i \mu_i} =\sum_{\mu\in M} y_\mu 2^{\sum_i \lambda_i \mu_i} =C_\lambda.
\end{equation}
Then for all $\mu$, we have that $x_\mu=y_\mu$.
\end{lemma}
\begin{proof}
We can again follow the proof of Wood. However,  when we define $A_\lambda$ we use the function $\bar{H}_{m,2,\mu}(z)=\sum_d a_d z_1^{d_1}z_2^{d_2}\dots z_m^{d_m}$ provided by Lemma \ref{lemma8p1}, instead of $H_{m,2,\mu}$. That is for every $\lambda\in M_0$ we define
\[A_\lambda=a_{\lambda_1-\lambda_2,\lambda_2-\lambda_3,\dots, \lambda_{m}}.\]   
To proceed with the proof we need to prove that $\sum_{\lambda\in M_0} A_\lambda C_\lambda$ is absolute convergent. We have
\begin{multline*}
\sum_{\lambda\in M_0} |A_\lambda C_\lambda| \le \sum_{\substack{d_1,\dots,d_m \ge 0\\
d_2+\dots+d_m \le \mu_1
}} |a_{d_1,d_2,\dots, d_m}| F^{m}2^{\sum_i d_i +\sum_i \frac{\sum_{k=i}^{m}d_k (\sum_{k=i}^{m}d_k-1)}{2}} \\
\le \sum_{\substack{d_1,\dots,d_m \geq 0\\
d_2+\dots+d_m \le \mu_1
}} E 2^{-b_1 d_1 -{d_1(d_1+1)}} F^{m}2^{\sum_i d_i +\sum_i \frac{\sum_{k=i}^{m}d_k (\sum_{k=i}^{m}d_k-1)}{2} }  .
\end{multline*}
For each choice of $d_2,\dots d_m$, the remaining sum over $d_1$ is a constant times\break
$
\sum_{d_1\geq 0} 2^{d_1(-b_1-\frac{1}{2}+d_2+\dots+d_m)-\frac{d_1^2}{2}},
$
which converges, so it follows that $\sum_{\lambda\in M} A_\lambda C_\lambda$ converges absolutely.
\end{proof}
The rest of the proof follows by repeating the arguments of Wood \cite{wood}.

To prove (\ref{wooii}) first we need to following lemma.
\begin{lemma}\label{tight}
For a random finite $2$-group $X$ we have
\[\mathbb{P}(X\not \in FA(2,k))\le \frac{\mathbb{E}|\Sur(X,\mathbb{Z}/2^{k+1}\mathbb{Z})|}{2^k}.\]
\end{lemma}
\begin{proof}
Observe that if a finite abelian $2$-group has exponent at least $2^{k+1}$, then it has at least $2^k$ surjective homomorphism to $\mathbb{Z}/2^{k+1}\mathbb{Z}$. Thus the statement follows from
\[\mathbb{E}|\Sur(X,\mathbb{Z}/2^{k+1}\mathbb{Z})|\ge \mathbb{P}(x\not \in AF(2,k))2^k.\]
\end{proof}
Now we prove the uniqueness of the measure $\nu$. Let $\nu$ be a measure satisfying the properties of (\ref{wooii}). Let $X$ be a random group with distribution $\nu$. Let $V\in FA_{odd}(2,k)$. Take any $m>k$. Then for any $W\in FA(2,m)$ we have
\[\mathbb{E}|\Sur(X\otimes \mathbb{Z}/2^m \mathbb{Z},W)|=\mathbb{E}|\Sur(X,W)|=|W|2^{\rank2(W)}.\]
Using the statement of (\ref{wooi}) we get that $X\otimes \mathbb{Z}/2^m \mathbb{Z}$ has distribution $\nu_m$. Thus
\begin{equation}\label{wooiieq}
\nu(V)=\mathbb{P}(X\isom V)=\mathbb{P}(X\otimes \mathbb{Z}/2^m \mathbb{Z}\isom V)=\nu_m(V). 
\end{equation}

This shows that the only possible measure is the one defined as follows. For $V\in FA_{odd}(2,k)$ we set $\nu(V)=\nu_m(V)$, where $m>k$. A similar argument as above shows that this does not depend on the choice of $m$ as long as $m>k$. An alternative way to express $\nu(V)$ is $\nu(V)=\lim_{m\to\infty} \nu_m(V)$. We need to prove that for any $W\in FA(2,k)$ we have $\Sur(\nu,W)=|W|2^{\rank2(W)}$. Let $\bar{\nu}$ be the push forward of the measure $\nu$ by the map $X\mapsto X\otimes \mathbb{Z}/2^k \mathbb{Z}$. It is enough to prove that $\bar{\nu}=\nu_k$. If $V$ has exponent smaller than $2^k$, then $\bar{\nu}(V)=\nu(V)=\nu_k(V)$. If $V\in FA_{odd}(2,k)\backslash FA_{odd}(2,k-1)$ then
\[\bar{\nu}(V)=\sum_{U\otimes \mathbb{Z}/2^k \mathbb{Z}\isom V} \nu(U)=\sum_{U\otimes \mathbb{Z}/2^k \mathbb{Z}\isom V} \lim_{m\to\infty}\nu_m(U)\le \lim_{m\to\infty}  \sum_{U\otimes \mathbb{Z}/2^k \mathbb{Z}\isom V} \nu_m(U)=\nu_k(V),\]
using Fatou's lemma and the fact that if $X_m$ has distribution $\nu_m$ for $m>k$, then $X_m\otimes \mathbb{Z}/2^k \mathbb{Z}$ has distribution $\nu_k$. Using the latter fact and Lemma \ref{tight} we obtain that
\begin{multline*}\nu_k(V)=\mathbb{P}(X_m\otimes \mathbb{Z}/2^k \mathbb{Z}\isom V)=\\ \sum_{\substack{U\otimes \mathbb{Z}/2^k \mathbb{Z}\isom V\\ U\in \cup_{i=k}^{m-1} FA_{odd}(2,i)}} \nu_m(U)+\sum_{\substack{U\otimes \mathbb{Z}/2^k \mathbb{Z}\isom V\\ U\in FA_{odd}(2,m)}} \nu_m(U)\le\\ \sum_{\substack{U\otimes \mathbb{Z}/2^k \mathbb{Z}\isom V\\ U\in \cup_{i=k}^{m-1} FA_{odd}(2,i)}} \nu(U) +\mathbb{P}(X_m\not\in FA(2,m-1))\le\\ \bar{\nu}(V)+\frac{\mathbb{E}|\Sur(X,\mathbb{Z}/2^{m}\mathbb{Z})|}{2^{m-1}}= \bar{\nu}(V)+2^{-(m-2)}.
\end{multline*}
Tending to infinity with $m$, we obtain that $\nu_k(V)\le \bar{\nu}(V)$. So indeed $\nu_k=\bar{\nu}$. 

The last statement of (\ref{wooii}) follows from (\ref{wooi}) and (\ref{wooiieq}). 

\end{proof}

In Lemma \ref{wool} above we concentrated only on the prime $2$ for simplicity, but the using the same argument we can handle finitely many primes simultaneously by following the argument of Wood \cite{wood}, that is, we have the following theorem. 
\fi

In the rest of the section we find a sequence of random groups, such that they have same limiting surjective moments  as the sequence of sandpile groups of $H_n$. The nice algebraic properties of these groups allow us to give an explicit formula for their limiting distribution. Then the previous theorem can be used to conclude that the sandpile group of $H_n$ has the same limiting distribution.

We start by showing that Lemma \ref{vanszimmetrikus} is true under slightly weaker conditions. 
\begin{lemma}\label{vanszimmetrikus2}
Assume that $n\ge 2|V|$. Let $q\in V^n$ be such that $G_q=V$. Let $r\in V^n$ such that $<q\otimes r>\in I_2$.  Then there is a symmetric matrix $A$ over $\mathbb{Z}$ such that $r=Aq$ and all the diagonal entries of $A$ are even. 
\end{lemma}
\begin{proof}
We start by the following lemma. As in Lemma \ref{vanszimmetrikus}, let $V=\bigoplus_{i=1}^\ell <v_i>$.
\begin{lemma}
There is an invertible integral matrix $B$, such that $B^{-1}$ is integral, and $q'=Bq$ satisfies that $m_{q'}(v_i)>0$ for every $1\le i \le \ell$.
\end{lemma}
\begin{proof}
Using the condition $n\ge 2|V|$ and $G_q=V$, we can choose $n-\ell$ components of $q$ such that they generate $V$. Due to symmetry we may assume that $q_{\ell+1},q_{\ell+2},\dots, q_n$ generates $V$. Let us define $q'=(v_1,v_2,\dots,v_\ell,q_{\ell+1},q_{\ell+2},\dots, q_n)$. We define the integral matrix $B=(b_{ij})$ by 
\begin{equation*}
b_{ij}=
\begin{cases}
1&\text{for } 1\le i=j\le n,\\
0&\text{for } 1 \le j<i\le n,\\
0& \text{for }\ell<i<j\le n,\\
0& \text{for }1\le i<j\le \ell .
\end{cases}
\end{equation*}
We still have not defined $b_{ij}$ for $1\le i\le \ell$ and $\ell<j\le n$. Since $q_{\ell+1},q_{\ell+2},\dots, q_n$ generates $V$ we can choose these entries such that $Bq=q'$. Since $B$ is an upper triangular integral matrix such that each diagonal entry is $1$, it is invertible and the inverse is an integral matrix.
\end{proof}
Let $B$ the matrix provided by the lemma above. Set $q'=Bq$ and $r'=\left(B^{-1}\right)^T r$. Observe that \[<q'\otimes r'>=<Bq\otimes \left(B^{-1}\right)^T r>=<B^{-1}Bq\otimes r>=<q\otimes r>\in I_2.\] Applying Lemma \ref{vanszimmetrikus}, we obtain a symmetric integral matrix $A'$ with even diagonal entries such that $r'=A'q'$. Consider $A=B^TA'B$. Then $A$ is a a symmetric integral matrix with even diagonal entries. Moreover, \[Aq=B^TA'Bq=B^T A' q'=B^Tr'=B^T \left(B^{-1}\right)^T r=r.\]
\end{proof}

\begin{lemma}\label{moment2k}
Let $V$ be a finite  abelian $2$-group. Assume that $2^k$ is divisible by the exponent of $V$. Let $A_n$ be uniformly chosen from the set of symmetric matrices in  $M_n(\mathbb{Z}/2^k\mathbb{Z})$, such that all the diagonal entries are even. 
Then we have
\[\lim_{n\to\infty} \mathbb{E}|\{q\in V^n|\quad G_q=V,\quad A_nq=0\}|=2^{\rank2(V)}|\wedge^2 V|.\] 
\end{lemma}
\begin{proof}
Take  any $q\in V^n$ such that $G_q=V$. 
Let $N_n$ be the set of symmetric matrices with even diagonal entries in $M_n(\mathbb{Z}/2^k\mathbb{Z})$. The distribution of $A_nq$ is the uniform distribution on the image of the $N_n\to V^n$  homomorphism  
${C\mapsto Cq}$.\break  From Lemma~\ref{vanszimmetrikus2} one can see that if $n$ is large enough then this image is \break $\{r\in V^n|<q\otimes r>\in I_2\}$, which has size $|V|^n\left(2^{\rank2(V)}|\wedge^2 V|\right)^{-1}$. It is clear that $0$ is always contained in the image, thus $\mathbb{P}(A_nq=0)=|V|^{-n}2^{\rank2(V)}|\wedge^2 V|$. Thus
\begin{multline*}
\lim_{n\to\infty} \mathbb{E}|\{q\in V^n|\quad G_q=V,\quad A_nq=0\}|=\\ \lim_{n \to\infty}\mathbb{E}|\{q\in V^n|\quad G_q=V\}|\frac{2^{\rank2(V)}|\wedge^2 V|}{|V^n|}=2^{\rank2(V)}|\wedge^2 V|.
\end{multline*} 
\end{proof}

Let $\mathbb{Z}_2$ be the ring of $2$-adic integers. Recall the fact that $\mathbb{Z}_2$ is the inverse limit of $\mathbb{Z}/2^k \mathbb{Z}$. Thus combining the lemma above with the analogue of Proposition \ref{prop1}, we get the following.


\begin{lemma}\label{masik2}
Let $\Symm_0(n)$ be the set of $n\times n$ symmetric matrices over $\mathbb{Z}_2$, such that all diagonal entries are even. Let $Q_n$ be a Haar-uniform element of $\Symm_0(n)$. For any finite abelian $2$-group $V$, we have
\[\lim_{n\to\infty}\mathbb{E}|\Sur(\cok(Q_n),V)|=2^{\rank2(V)}|\wedge^2 V|.\]
Moreover, if $\overline{Q}_n\in M_n(\mathbb{Z}/2\mathbb{Z})$ is obtained by reducing each entry of $Q_n$ modulo $2$, then $\overline{Q}_n$ is a symmetric matrix with $0$ as its
diagonal entries. Consequently, $\rankk(\cok(Q_n))\equiv n$ modulo $2$.  \qed
\end{lemma}



The next lemma gives an explicit formula for the limiting distribution of $\cok(Q_n)$. The author is grateful to Melanie Wood who proved this result for him.

\begin{lemma}\label{explicitnu}(Wood \cite{personal})
For any finite abelian $2$-group $G$ of odd rank, we have
\begin{multline*}
\nu(G)=\lim_{\substack{n\to\infty\\n\text{ is odd}}} \mathbb{P}(\cok(Q_n)\isom G)=\\
2^{\rankk(G)}\frac{|\{\phi:G \times G \to \mathbb{C}^* \SBP \}|}{|G||\Aut(G)|} \prod_{j=0}^\infty (1-2^{-2j-1}).
\end{multline*}
\end{lemma}

\begin{proof}
Assume that $G=\bigoplus_{i=1}^k(\mathbb{Z}/2^{e_i}\mathbb{Z})^{n_i}$ where $e_1>e_2>\dots>e_k>0$.

We consider $\mathbb{Z}_2^n$ as a $\mathbb{Z}_2$ module. Let $L_n(G)$ be the set of submodules $M$ of $\mathbb{Z}_2^n$ such that $\mathbb{Z}_2^n/M$ is isomorphic to $G$.

\[\mathbb{P}(\cok(Q_n)\isom G)=\P(\RowS(Q_n)\in L_n(G))=\sum_{M\in L_n(G)} \P(\RowS(Q_n)=M).\]

Let $\mu_n$ be the Haar probability measure on $\Symm_0(n)$. Fix $M\in L_n(G)$. We are interested in the probability \[\mathbb{P}(\RowS(Q_n)=M)=\mu_n(\{S\in \Symm_0(n)|\RowS(S)=M\}).\] 

Fix any (not necessary symmetric) $n\times n$ matrix $N$ over $\mathbb{Z}_p$ such that \break ${\RowS(N)=M}$. Observe that
\[\{S\in \Symm_0(n)|\RowS(S)=M\} =\{CN|\quad CN\in\Symm_0(n), C\in GL_n(\mathbb{Z}_2)\}.\]
Since $\mathbb{Z}_p$ is a principal ideal domain $N$ has a Smith normal form, that is, we can find $A,B\in GL_n(\mathbb{Z}_2)$ such that $D=ANB$ is a diagonal matrix. 
Since each nonzero element of $\mathbb{Z}_2$ can written as $2^du$, where $d$ is a nonnegative integer, $u$ is a unit in $\mathbb{Z}_2$, we may assume each entry of $D$ is of the form $2^d$ for some $d$. 
But since $\mathbb{Z}_2^n/\RowS(D)\isom \mathbb{Z}_2^n/\RowS(N)\isom G$, we know exactly what is $D$. Let $ n_{k+1}=n-\sum_{i=1}^k n_i$, and $e_{k+1}=0$. From now on it will be convenient to view $n\times n$ matrices as $(k+1)\times (k+1)$ block matrices, where the block at the position $(i,j)$ is an $n_i\times n_j$ matrix. 
Then $D$ is a block matrix $(D_{ij})_{i,j=1}^{k+1}$ where all the off-diagonal blocks are zero and $D_{ii}=2^{e_i}I$.   

Observe that map $S\mapsto B^T S B$ is an automorphism of the abelian group $\Symm_0(n)$. Thus, it pushes forward $\mu_n$ to $\mu_n$, which gives us
\begin{align*}
\mu_n(\{CN|&\quad CN\in\Symm_0(n), C\in GL_n(\mathbb{Z}_2)\})\\&= \mu_n(\{B^TCNB|\quad B^TCNB \in\Symm_0(n), C\in GL_n(\mathbb{Z}_2)\})\\&=
\mu_n(\{B^TCA^{-1}ANB|\quad B^TCA^{-1}ANB \in\Symm_0(n), C\in GL_n(\mathbb{Z}_2)\})\\&= \mu_n(\{B^TCA^{-1}D|\quad B^TCA^{-1}D \in\Symm_0(n), C\in GL_n(\mathbb{Z}_2)\})\\&=
\mu_n(\{FD|\quad FD \in\Symm_0(n), F\in GL_n(\mathbb{Z}_2)\}).
\end{align*}

We consider $F=(F_{ij})_{i,j=1}^{k+1}$ as $(k+1)\times (k+1)$ block matrix as it was described above. Then $FD\in \Symm_0(n)$ if and only if for every $i<j$, we have
\begin{equation}\label{oszt}
F_{ij} = 2^{e_i-e_j}F_{ji}^T 
\end{equation} 
and the diagonal entries of $F_{k+1,k+1}$ are even. Assuming that $F$ has these properties, when does $F$ belong to   $GL_n(\mathbb{Z}_2)$? Observe that  $F\in GL_n(\mathbb{Z}_2)$ if and only if the mod 2 reduction $\overline{F}$ of $F$ is invertible, but Equation \eqref{oszt} tells us $\overline{F}$ is a block lower triangular matrix, so $F\in GL_n(\mathbb{Z}_2)$  if and only if $F_{ii}\in GL_{n_i}(\mathbb{Z}_2)$ for each $i$.

From this it follows that $\{FD|\quad FD \in\Symm_0(n), F\in GL_n(\mathbb{Z}_2)\}$ consists of all block matrices $H\in \Symm_0(n)$, such that
\begin{enumerate}
\item For $1\le i,j\le k+1$ all entries of the block $H_{ij}$ is divisible by $2^{\max(e_i,e_j)}$.

\item For $1\le i \le k+1$ the mod $2$ reduction of the matrix $2^{-e_i}H_{ii}$ is an invertible symmetric matrix over $\mathbb{F}_2$. Moreover, if $i=k+1$, then all its diagonal entries are zero.
\end{enumerate}

Let $p_m$ be the probability that a uniform random symmetric $m\times m$ matrix over $\mathbb{F}_2$ is invertible, and  let $p_m'$ be the probability that a uniform random symmetric $m\times m$ matrix over $\mathbb{F}_2$ is invertible and all its diagonal entries are zero.

\begin{align*}
\mathbb{P}(\RowS(Q_n)=M)&=\mu_n(\{FD|\quad FD \in\Symm_0(n), F\in GL_n(\mathbb{Z}_2)\})\\&=
2^{n}p_{n_{k+1}}'\prod_{i=1}^k p_{n_i} 2^{e_i\left(n_i(n-\sum_{j=1}^i n_j)+{{n_i+1}\choose{2}}\right)} .
\end{align*}
In particular, this does not depend on the choice of $M\in L_n(G)$. Thus, we obtain that
\[\mathbb{P}(\cok(Q_n)\isom G)=|L_n(G)|2^{n}p_{n_{k+1}}'\prod_{i=1}^k p_{n_i} 2^{e_i\left(n_i(n-\sum_{j=1}^i n_j)+{{n_i+1}\choose{2}}\right)}.\]

Now let $Q_n' $ be a Haar-uniform $n\times n$ symmetric matrix over $\mathbb{Z}_2$. A very similar calculation as above gives that
\[\mathbb{P}(\cok(Q_n')\isom G)=|L_n(G)|p_{n_{k+1}}\prod_{i=1}^k p_{n_i} 2^{e_i\left(n_i(n-\sum_{j=1}^i n_j)+{{n_i+1}\choose{2}}\right)}.\]  
Therefore,
\begin{align}\label{rank2sz}
\frac{\mathbb{P}(\cok(Q_n)\isom G)}{\mathbb{P}(\cok(Q_n')\isom G)}&=2^n\frac{p_{n_{k+1}}'}{p_{n_{k+1}}}=2^{n-n_{k+1}}\frac{2^{n_{k+1}}p_{n_{k+1}}'}{p_{n_{k+1}}}\\&=2^{\rankk(G)}\frac{2^{n_{k+1}}p_{n_{k+1}}'}{p_{n_{k+1}}}=2^{\rankk(G)}.\nonumber
\end{align}
The last equality follows from the results of MacWilliams \cite{macw}. Note that here we needed to use that $n$ and $\rankk(G)$ are both odd, therefore $n_{k+1}$ is even. As we already mentioned in the Introduction in line \eqref{szimformula} by the result of \cite{clp14}, we have
\begin{multline*}
\lim_{n\to\infty} \mathbb{P}(\cok(Q_n')\isom G)\\=\frac{|\{\phi:G \times G \to \mathbb{C}^* \SBP \}|}{|G||\Aut(G)|} \prod_{j=0}^\infty (1-2^{-2j-1}).
\end{multline*}
Combining this with line \eqref{rank2sz} above, we get the statement.
\end{proof}


Now we can prove the remaining part of Theorem \ref{CohenlenstraUD}
\begin{proof}(Theorem \ref{CohenlenstraUD} for even $d$)

Let $p_i^{k_i}$ be the exponent of $G_i$. 

Let $Q_{n,1}$ be  a Haar-uniform element of the the set of $(2n-1)\times (2n-1)$ symmetric matrices over $\mathbb{Z}_2$, where all the 
diagonal entries are even. For $i>1$, let $Q_{n,i}$ be  a Haar-uniform element of the the set of $(2n-1)\times (2n-1)$ symmetric matrices over $\mathbb{Z}_{p_i}$.  All the choices are made independently. Let $\bar{Q}_{n,i}\in M_{2n-1}(\mathbb{Z}/p_i^{k_i+1} \mathbb{Z})$ be the mod $p_i^{k_i+1}$ reduction of $Q_{n,i}$.

Let $a=\prod_{i=1}^s p_i^{k_i+1}$. Let $X_n$ be the sandpile group $\Gamma_{2n}$ of $H_{2n}$. Let \break $Y_n=\bigoplus_{i=1}^s \cok(\bar{Q}_{n,i})$. Let  $V$ be a finite abelian group with exponent dividing $a$. Then, from Theorem \ref{momentumokUD}, we have
\[\lim_{m\to\infty} \mathbb{E}|\Sur(X_n,V)|=2^{\rank2(V)}|\wedge^2 V|.\]
Let $V_i$ be the $p_i$-Sylow subgroup of $V$. From Lemma \ref{moment2k}, we have 
\[\lim_{n\to\infty} \mathbb{E}|\Sur(\cok(\bar{Q}_{n,1}),V_1)|=2^{\rank2(V_1)}|\wedge^2 V_1|.\]
For $i>1$, from \cite[Theorem 11]{clp14}, we have
\[\lim_{n\to\infty} \mathbb{E}|\Sur(\cok(\bar{Q}_{n,1}),V_1)|=|\wedge^2 V_i|.\]
It is also clear that 
\[ |\Sur(Y_n,V)|=\prod_{i=1}^s |\Sur(\cok(\bar{Q}_{n,i}),V_i)|.\]
Thus, from the independence of $Q_{n,i}$, we get that
\begin{align*}
\lim_{n\to\infty}\mathbb{E}|\Sur(Y_n,V)|&=\prod_{i=1}^s \lim_{n\to\infty} \mathbb{E}|\Sur(\cok(\bar{Q}_{n,i}),V_i)|\\ &=2^{\rank2(V_1)}\prod_{i=1}^{s}|\wedge^2 V_i|=2^{\rank2(V)}|\wedge^2 V|.
\end{align*}

From Lemma \ref{explicitnu} and \cite[Theorem 2]{clp14}, we have

\begin{multline*}
\lim_{n\to\infty}\mathbb{P}( Y_n\otimes \mathbb{Z}/a\mathbb{Z}\isom \bigoplus_{i=1}^s G_i)=\lim_{n\to\infty} \prod_{i=1}^s \mathbb{P}(\cok(Q_{n,i})\isom G_i)=\\2^{\rankk(G_1)}\prod_{i=1}^s \left(\frac{|\{\phi:G_i \times G_i \to \mathbb{C}^* \SBP \}|}{|G_i||\Aut(G_i)|} \prod_{j=0}^\infty (1-p_i^{-2j-1})\right).
\end{multline*}

Note that $\bigoplus_{i=1}^s \Gamma_{n,i}\isom\bigoplus_{i=1}^s G_i$ if and only if $X_n\otimes \mathbb{Z}/a\mathbb{Z}\isom \bigoplus_{i=1}^s G_i$. Note that both $\rank2(X_n\otimes \mathbb{Z}/2\mathbb{Z})$ and $\rank2(Y_n\otimes \mathbb{Z}/2\mathbb{Z})$ are odd. Therefore, Theorem \ref{T:MomDetDetail} can be applied to finish the proof. 
\end{proof}

\section{The sublinear growth of rank}\label{sublin}
In this section we prove Theorem \ref{rankThm}. Let $\Gamma_n$ be the sandpile group of $H_n$. We start by a simple lemma. Recall that $\rankk_p(\tors(\Gamma_n))$ is the rank of the $p$-Sylow subgroup of $\tors(\Gamma_n)$.
\begin{lemma}\label{nagyprim}
There is a constant $c_d$ such that $|\tors(\Gamma_n)|< c_d^n$. Consequently, for any prime $p$, we have
\[\rankk_p(\tors(\Gamma_n))\le \frac{n\log c_d}{\log p}.\]
\end{lemma}
\begin{proof}
Let $v_1,v_2,...,v_k=n$ be a subset of the vertices of $H_n$, such that each connected component of $H_n$ contains exactly one of them. (With high probability $k=1$.) Let $\Delta_0$ be the matrix obtained from the Laplacian by deleting the rows and columns corresponding to the vertices $v_1,v_2,\dots,v_k$. Observe that \break $\tors(\Gamma_n)=|\det \Delta_0|$. Each row of $\Delta_0$ has Euclidean norm at most $c_d=\sqrt{2d^2}$. Thus,  $\tors(\Gamma_n)=|\det \Delta_0|\le c_d^{n-k}<c_d^n$, from Hadamard's inequality \cite{hadam}. The proof of the second statement is straightforward from this.
\end{proof}

The lemma above will be used for large primes, for small primes we will use the next lemma.

\begin{lemma}\label{kicsiprim}
For every prime $p$, there is a constant $C_p$ such that for any $n$ and $\varepsilon>0$, we have
\[\mathbb{P}(\rankk(\Gamma_n\otimes \mathbb{Z}/p\mathbb{Z})\ge \varepsilon n)\le C_p p^{-\varepsilon n}.\]
\end{lemma}
\begin{proof}
It is an easy consequence of Corollary \ref{correduced22} and Proposition \ref{prop1} that \[\lim_{n\to\infty} \mathbb{E} |\Hom(\Gamma_n\otimes \mathbb{Z}/p\mathbb{Z},\mathbb{Z}/p\mathbb{Z})|\] exists. This implies that there is a constant $C_p$ such that  \[\mathbb{E} |\Hom(\Gamma_n\otimes \mathbb{Z}/p\mathbb{Z},\mathbb{Z}/p\mathbb{Z})|\le C_p\] for any $n$. Note that $|\Gamma_n\otimes \mathbb{Z}/p\mathbb{Z}|=|\Hom(\Gamma_n\otimes \mathbb{Z}/p\mathbb{Z},\mathbb{Z}/p\mathbb{Z})|$. Thus, from Markov's inequality
\begin{align*}
\mathbb{P}(\rankk(\Gamma_n\otimes \mathbb{Z}/p\mathbb{Z})\ge \varepsilon n)&=\mathbb{P}(|\Gamma_n\otimes \mathbb{Z}/p\mathbb{Z}|\ge p^{\varepsilon n})\le
 p^{-\varepsilon n} \mathbb{E}|\Gamma_n\otimes \mathbb{Z}/p\mathbb{Z}|\\&= p^{-\varepsilon n} \mathbb{E}|\Hom(\Gamma_n\otimes \mathbb{Z}/p\mathbb{Z},\mathbb{Z}/p\mathbb{Z})|\le  C_p p^{-\varepsilon n}.
\end{align*} 
\end{proof}
Now we are ready to prove Theorem \ref{rankThm}. Take any $\varepsilon>0$. Set $K=\exp(\varepsilon^{-1}\log c_d)$. Let $\{p_1,p_2,\dots,p_s\}$ be the set of primes that are at most $K$. Using Lemma \ref{kicsiprim}, we get that
\[\mathbb{P}(\rankk(\Gamma_n\otimes \mathbb{Z}/p_i\mathbb{Z})\ge \varepsilon n\text{ for some }i\in\{1,2,\dots,s\})\le \sum_{i=1}^s C_{p_i} p_i^{-\varepsilon n}.\] 
Since $\sum_{n=1}^{\infty} \sum_{i=1}^s C_{p_i} p_i^{-\varepsilon n}$ is convergent, the Borel-Cantelli lemma gives us the following. With probability $1$ there is an $N$ such that for every $n>N$ and \break $i=1,2,\dots,s$, we have  $\rankk(\Gamma_n\otimes \mathbb{Z}/p_i\mathbb{Z})< \varepsilon n$. By the choice of $K$ and Lemma \ref{nagyprim}, for a prime $p>K$, we have   $\rankk_p(\tors(\Gamma_n))\le \varepsilon n$. Write $\Gamma_n$ as $\Gamma_n=\mathbb{Z}^f \times \tors(\Gamma_n)$. Then for $n>N$, we have
\begin{align*}
\rankk(\Gamma_n)&=f+\max_{p\text{ is a prime}} \rankk_p(\tors(\Gamma_n))\\&\le  \rankk(\Gamma_n\otimes \mathbb{Z}/2\mathbb{Z}) +\max_{p\text{ is a prime}} \rankk_p(\tors(\Gamma_n))\le \varepsilon n+\varepsilon n.
\end{align*}
Tending to $0$ with $\varepsilon$, we get the statement.

\section{Bounding the probabilities of non-typical events}\label{Bounding}

At several points of the paper we need to bound the probability of that something is not-typical. These estimates are all based on the following lemma.

\begin{lemma}\label{bazuma}
Given $0\le a,b\le n$, let $A$ and $B$ be a uniform independent random subset of $\{1,2,\dots,n\}$ such that $|A|=a$ and $|B|=b$. Then for any $k>0$, we have
\[\mathbb{P}\left(\left||A\cap B|-\frac{ab}n\right|\ge k\right)\le
 2\exp\left(-\frac{2k^2}{a}\right)\le
2\exp\left(-\frac{2k^2}{n}\right).\] 
\end{lemma}  
\begin{proof}
Note that $A\cap B$ has the same distribution as $\sum_{i=1}^a X_i$, where $X_1,X_2\dots,X_a$ is a random sample drawn without replacement from an $n$ element multiset, where $1$ has multiplicity $b$ and $0$ has multiplicity $n-b$. Then the statement follows from \cite[Proposition 1.2]{concent}. 
\end{proof}

Applying this iteratively we get the following lemma.
\begin{lemma}\label{itbazuma}
Given $0\le a_1,a_2,...,a_d\le n$, let $A_1,A_2,...,A_d$ be uniform independent random subsets of $\{1,2,\dots,n\}$ such that $|A_i|=a_i$ for $i=1,2,\dots, d$. Then we have
\begin{align*}
\mathbb{P}\left(\left||A_1\cap\dots\cap A_d|-n\prod_{i=1}^d \frac{a_i}{n}\right|\ge (d-1) k\right)&\le 2(d-1)\exp\left(-\frac{2k^2}{a_1}\right)\\& \le 2(d-1)\exp\left(-\frac{2k^2}{n}\right).
\end{align*}
\end{lemma}
\begin{proof}
The proof is by induction. For $d=2$, it is true as Lemma \ref{bazuma} shows. Now we prove for $d$. By induction  
\[\mathbb{P}\left(\left||A_1\cap\dots A_{d-1}|-n\prod_{i=1}^{d-1} \frac{a_i}{n}\right|\ge (d-2) k\right)\le 2(d-2)\exp\left(-\frac{2k^2}{a_1}\right).\]
Using Lemma \ref{bazuma} for $A_1\cap\dots A_{d-1}$ and $A_d$ and the fact that $|A_1\cap\dots A_{d-1}|\le a_1$, we have
\[\mathbb{P}\left(\left||A_1\cap\dots A_{d}|-\frac{|A_1\cap\dots\cap A_{d-1}|a_d}{n}\right|\ge  k\right)\le 2\exp\left(-\frac{2k^2}{a_1}\right).\]
Thus, with probability at least $1-2(d-1)\exp\left(-\frac{2k^2}{a_1}\right)$, we have that
\[\left||A_1\cap\dots A_{d}|-\frac{|A_1\cap\dots\cap A_{d-1}|a_d}{n}\right|\le  k\]
and for
\[\Delta=|A_1\cap\dots A_{d-1}|-n\prod_{i=1}^{d-1} \frac{a_i}{n},\]
the inequality $|\Delta|\le (d-2)k$ holds. Therefore,
\begin{align*}
\left||A_1\cap\dots\cap A_d|-n\prod_{i=1}^d \frac{a_i}{n}\right|&=\left||A_1\cap\dots\cap A_d|-\frac{a_d(|A_1\cap \dots \cap A_{d-1}|-\Delta)}{n}\right|\\&\le
\left||A_1\cap\dots\cap A_d|-\frac{a_d|A_1\cap \dots \cap A_{d-1}|}{n}\right|+\frac{a_d|\Delta|}n\\&\le k+(d-2)k\le (d-1)k.\qedhere
\end{align*}
\end{proof}

Next we give the analogue of Lemma \ref{bazuma} for uniform random perfect matchings.
\begin{lemma}\label{bazumaM}
Assume that $n$ is even. Let $A$ and $B$ be two fixed subsets of \break $\{1,2,\dots,n\}$,  let  $|A|=a$ and $|B|=b$. Let $M$ be uniform random perfect matching on the set $\{1,2,\dots,n\}$. Let $X$ be the number of elements in $A$ that are paired with an element in $B$ in the matching $M$. Then for any $k>0$, we have
\[\mathbb{P}\left(\left|X-\frac{ab}n\right|\ge 4k\right)\le
 6\exp\left(-\frac{2k^2}{a}\right)\le
6\exp\left(-\frac{2k^2}{n}\right).\] 
\end{lemma}  
\begin{proof}
Observe that the uniform random matching $M$ can be generated as follows. First we partition the set $\{1,2,\dots,n\}$ into two disjoint subsets $H_1$ and $H_2$ of size $\frac{n}{2}$ uniformly at random. Then we consider a uniform random perfect matching between $H_1$ and $H_2$. For $i\in \{1,2\}$, let $a_i=|A\cap H_i|$, and let $b_i=|B\cap H_i|$. Let $X_i$ be the number of element in $A\cap H_i$ that are paired with an element in $B$. From Lemma \ref{bazuma}, we have
\begin{align*}
\mathbb{P}\left(\left|a_1-\frac{a}2\right|\ge k\right)&\le
 2\exp\left(-\frac{2k^2}{a}\right),\\
\mathbb{P}\left(\left|X_1-\frac{2a_1b_2}{n}\right|\ge k\right)&\le
 2\exp\left(-\frac{2k^2}{a_1}\right),\\
 \mathbb{P}\left(\left|X_2-\frac{2a_2b_1}{n}\right|\ge k\right)&\le
 2\exp\left(-\frac{2k^2}{a_2}\right).
\end{align*}
It follows from the union bound that with probability at least $1-6\exp\left(-\frac{2k^2}{a}\right)$, we have that
\[\left|a_1-\frac{a}2\right|< k,\quad  \left|X_1-\frac{2a_1b_2}{n}\right|<k\text{ and }\left|X_2-\frac{2a_2b_1}{n}\right|<k.\]
On this event
\begin{align*}
\left|X-\frac{ab}n\right|&=\left|\left(X_1-\frac{ab_2}{n}\right)+\left(X_2-\frac{ab_1}{n}\right)\right|\\
                         &\le \left|X_1-\frac{ab_2}{n}\right|+\left|X_2-\frac{ab_1}{n}\right|\\
                         &\le \left|X_1-\frac{2a_1b_2}{n}\right|+\left|\frac{2a_1b_2}{n}-\frac{ab_2}{n} \right|+\left|X_2-\frac{a_2b_1}{2n}\right|+\left|\frac{2a_2b_1}{n}-\frac{ab_1}{n} \right|\\
                         &<2k+\frac{2b_1}{n}\left|a_2-\frac{a}{2}\right|+\frac{2b_2}{n}\left|a_1-\frac{a}{2}\right|<4k.
\end{align*}
\end{proof}
Applying this iteratively,  we can get  a lemma similar to Lemma \ref{itbazuma}.


\bibliographystyle{amsplain}

\appendix

\section*{Appendix. Proofs omitted from Section \ref{parosparos}}

\begin{proof}[\textbf{The proof of Lemma \ref{L:Hacts}}]
We only prove the first statement. The second statement is the same as  \cite[Lemma 8.1]{wood}, and it can be proved essentially the same way.

We define analytic functions
$$
G(z_1)=\prod_{\substack{j> b_1\\j\text{ is odd}}}(1-\frac{z_1}{2^j}) =\sum_{d_1\geq 0} c_{d_1} z^{d_1}
$$
and
\begin{align*}
H(z_2,&\dots,z_m)=\\&
\prod_{j= b_1+b_2+1}^{2b_1}(1-\frac{z_2}{2^j}) \prod_{j= b_1+b_2+b_3+1}^{b_1+2b_2}(1-\frac{z_3}{2^j}) \cdots
\prod_{j= b_1\dots+b_m+1}^{b_1+\dots+b_{m-2}+2b_{m-1}}(1-\frac{z_m}{2^j})=\\&
\sum_{d_2,\dots,d_m \geq 0} e_{d_1,\dots,d_m} z_2^{d_2} \cdots z_m^{d_m}.
\end{align*}
In each of the $z_i$ separately, for $2\leq i \leq m$, we have that $H$ is a polynomial of degree $b_{i-1}-b_i$.
We then have an entire, analytic function
in  $m$ variables
$$
H_{m,2,b}(z)=G(z_1)H(z_2,\dots,z_m)=\sum_{\substack{d_1,\dots,d_m \geq 0\\
d_2+\dots+d_m\leq b_1
}} a_{d_1,\dots,d_m} z_1^{d_1} \cdots z_m^{d_m}.
$$

We now estimate the size of the $a_d$.  
We see that $a_d=c_{d_1} e_{d_2,\dots d_m}$.
We have that 
$
G(4z)=(1-\frac{z}{2^{b_1}})G(z).
$
So
$
c_n 4^n=c_n-2^{-b_1}c_{n-1}.
$
Thus 
$
c_n= - \frac{2^{-b_1}c_{n-1}}{4^n-1}, 
$
and by induction,
$
c_n=(-1)^n \frac{2^{-b_1n}}{\prod_{i=1}^n (4^i-1)}.
$
So
$
|c_n|\leq 2^{-b_1n - n(n+1)}  \prod_{i\geq 1} (1-4^{-i})^{-1} .
$
Thus, 
$$
a_d\leq \frac{1}{\prod_{i\geq 1} (1-4^{-i})} 2^{-b_1d_1 - d_1(d_1+1)} \max_{d_2,\dots,d_m}e_{d_2,\dots d_m}.
$$

Now we check the final statements of the lemma.
If $f>b$, suppose $f_i=b_i$ for $i\leq t$ and $f_{t+1} > b_{t+1}$ for some $0\leq t \leq m-1$.
Then, in particular $f_1+\dots+ f_i = b_1+\dots+ b_i$ for $i \leq t$, and 
$f_1+\dots+ f_{t+1} \geq b_1+\dots+ b_{t+1}+1$.
However, (when $t\geq 1$) since $f_{t+1}\leq f_t=b_t,$ we have $f_1+\dots+ f_{t+1} \leq  b_1+\dots+ b_{t-1}+2b_{t}.$
Since $H$ vanishes whenever $z_{t+1}=p^{k}$ for integers $k$ with $b_1+\dots+ b_{t+1}+1 \leq k
\leq b_1+\dots+ b_{t-1}+2b_{t},$ we obtain the desired vanishing.

For the last statement, we first note that since the  product in the definition of $G$ is absolutely
convergent, we have that $z_1=p^{b_1}$ is not a root of $G$.  Then we observe all the other finitely many factors in $H$ are non-zero in this case as well.
\end{proof}

\begin{proof}[\textbf{The proof of Lemma \ref{T:Momdet}}]
 We will induct on the size of $\mu$ in the lexicographic total ordering (we take the lexicographic ordering for partitions and then the lexicographic ordering on top of that for tuples of partitions).  Suppose we have $x_\pi=y_\pi$ for every $\pi<\nu$.

We use Lemma~\ref{L:Hacts} to find $H_{m_j,p_j,\nu^j}(z)=\sum_d a(j)_d z_1^{d_1}\dots z_{m_j}^{d_{m_j}}.$ Note the definition of  $H_{m_j,p_j,\nu^j}(z)$ is different for $j=1$ and $j>1$. Namely, for $j=1$ we use the first part of Lemma~\ref{L:Hacts}, and for $j>1$ we use the second part.

For $\lambda\in M$, we define
$$
A_\lambda=\prod_{j=1}^{s} a(j)_{\lambda^j_1-\lambda^j_2,\lambda^j_2-\lambda^j_3,\dots, \lambda^j_{m_j}}.
$$ 
We wish to show that the sum
$
\sum_{\lambda\in M} A_\lambda C_\lambda
$
converges absolutely.
We have
\begin{align*}
 \sum_{\lambda\in M} |A_\lambda C_\lambda| 
&\leq\sum_{\lambda\in M} 2^{\lambda_1^1}\prod_{j=1}^{s} \left| a(j)_{\lambda^j_1-\lambda^j_2,\lambda^j_2-\lambda^j_3,\dots, \lambda^j_{m_j}}{F^{m_j}p_j^{\sum_i \frac{\lambda^{j}_i(\lambda^{j}_i-1)}{2}}}  \right|\\
&= \left(\sum_{\lambda\in M_1}\left|a(1)_{\lambda_1-\lambda_2,\lambda_2-\lambda_3,\dots, \lambda_{m_1}}{F^{m_1}2^{\lambda^1+\sum_i \frac{\lambda_i(\lambda_i-1)}{2}}}  \right| \right)\\
&\quad\cdot \prod_{j=2}^{s} \sum_{\lambda\in M_j}  \left|a(j)_{\lambda_1-\lambda_2,\lambda_2-\lambda_3,\dots, \lambda_{m_j}}{F^{m_j}p_j^{\sum_i \frac{\lambda_i(\lambda_i-1)}{2}}}  \right|.
\end{align*}
First we investigate the first term in the product above. We drop the index $1$, and let $b=\nu^1$.  We apply the first part of Lemma~\ref{L:Hacts}  to obtain
\begin{multline*}
\sum_{\substack{d_1,\dots,d_m \geq 0\\
d_2+\dots+d_m \leq b_1
}} |a_{d_1,d_2,\dots, d_m}| F^{m}2^{\sum_i d_i+\sum_i \frac{\sum_{k=i}^{m}d_k (\sum_{k=i}^{m}d_k-1)}{2}}  
\leq \\ \sum_{\substack{d_1,\dots,d_m \geq 0\\
d_2+\dots+d_m \leq b_1
}} E 2^{-b_1 d_1 -{d_1(d_1+1)}} F^{m}2^{\sum_i d_i+\sum_i \frac{\sum_{k=i}^{m}d_k (\sum_{k=i}^{m}d_k-1)}{2} }  .
\end{multline*}
For each choice of $d_2,\dots d_m$, the remaining sum over $d_1$ is a constant times\break
$
\sum_{d_1\geq 0} 2^{d_1(-b_1-\frac{1}{2}+d_2+\dots+d_m)-\frac{d_1^2}{2}},
$
which converges.

We now investigate the inner sum in the second term.  We drop the $j$ index, and let $b=\nu^j$.  We apply the second part of Lemma~\ref{L:Hacts} to obtain
\begin{multline*}
\sum_{\substack{d_1,\dots,d_m \geq 0\\
d_2+\dots+d_m \leq b_1
}} |a(j)_{d_1,d_2,\dots, d_m}| F^{m}p^{\sum_i \frac{\sum_{k=i}^{m}d_k (\sum_{k=i}^{m}d_k-1)}{2}}  
\leq \\ \sum_{\substack{d_1,\dots,d_m \geq 0\\
d_2+\dots+d_m \leq b_1
}} E p^{-b_1 d_1 -\frac{d_1(d_1+1)}{2}} F^{m}p^{\sum_i \frac{\sum_{k=i}^{m}d_k (\sum_{k=i}^{m}d_k-1)}{2} }  .
%
\end{multline*}
For each choice of $d_2,\dots d_m$, the remaining sum over $d_1$ is a constant times\break 
$
\sum_{d_1\geq 0} p^{d_1(-b_1-1+d_2+\dots+d_m)},
$
which converges, so it follows that $\sum_{\lambda\in M} A_\lambda C_\lambda$ converges absolutely.

Suppose we have $x_\mu$ for $\mu \in M_0$ all non-negative, such that for all $\lambda\in M$,
$$
\sum_{\mu\in M_0} x_\mu \prod_{j=1}^{s}  p_j^{\sum_i \lambda^{j}_i \mu^{j}_i}=C_\lambda.
$$
So we have that
$$
\sum_{\lambda\in M} \sum_{\mu\in M_0} A_\lambda  x_\mu \prod_{j=1}^{s}  p_j^{\sum_i \lambda^{j}_i \mu^{j}_i}
$$
converges absolutely.
Thus,
\begin{align*}
 \sum_{\lambda\in M} A_\lambda C_\lambda &=\sum_{\lambda\in M} \sum_{\mu\in M_0} A_\lambda  x_\mu \prod_{j=1}^{s}  p_j^{\sum_i \lambda^{j}_i \mu^{j}_i}\\
&=\sum_{\mu\in M_0} x_\mu \sum_{\lambda\in M}  A_\lambda  \prod_{j=1}^{s}  p_j^{\sum_i \lambda^{j}_i \mu^{j}_i}\\
&=\sum_{\mu\in M_0} x_\mu  \prod_{j=1}^{s}  \sum_{\lambda\in M_j}  
a(j)_{\lambda_1-\lambda_2,\lambda_2-\lambda_3,\dots, \lambda_{m_j}}
  p_j^{\sum_i \lambda_i \mu^{j}_i}.
\end{align*}
Now we consider the inner sum.  Again we drop the $j$ indices.  
We have
\begin{align*}
\sum_{\lambda\in M_j}  &
a(j)_{\lambda_1-\lambda_2,\lambda_2-\lambda_3,\dots, \lambda_{m}}
  p^{\sum_i \lambda_i \mu_i}\\&=
  \sum_{d_1,\dots,d_{m}\geq 0}  a(j)_{d_1,\dots,d_m}  (p^{\mu_1})^{d_1}(p^{\mu_1+\mu_2})^{d_2} \cdots (p^{\mu_1+\dots+\mu_m})^{d_m}
  \\
&= H_{m,p,\nu}(p^{\mu_1}, p^{\mu_1+\mu_2}, \dots, p^{\mu_1+\dots+\mu_m}).
\end{align*}
If $\mu,\nu\in M_0$ and $\mu > \nu$ (in the lexicographic total ordering), then 
some $\mu^j>\nu^j$ and so for $m=m_j$ and $p=p_j$,
by Lemma~\ref{L:Hacts}, 
 $H_{m,p,\nu^j}(p^{\mu_1}, p^{\mu_1+\mu_2}, \dots, p^{\mu_1+\dots+\mu_m})=0$.  Furthermore, if $\mu=\nu$, then for each (implicit) $j$, we have\break  
$H_{m,p,\nu}(p^{\mu_1}, p^{\mu_1+\mu_2}, \dots, p^{\mu_1+\dots+\mu_m})\ne0.$
So for some non-zero $u$,
$$
 \sum_{\lambda\in M} A_\lambda C_\lambda  =  x_\nu u 
+ \sum_{\mu\in M_0, \mu < \nu } x_\mu \sum_{\lambda\in M}  A_\lambda  \prod_{j=1}^{s}  p_j^{\sum_i \lambda^{j}_i \mu^{j}_i}.
$$
So since by assumption $x_\mu$ with $\mu<\nu$ we determined by the $C_\lambda$, we conclude that
$x_\nu$ is determined as well.
\end{proof}

\begin{proof}[\textbf{The proof of Lemma \ref{T:MomDetDetail}}]
First, we will suppose that the limits
$
\lim_{n\ra\infty} \P(X_n\tensor \Z/a\Z \isom H)
$
exist, and from that show that
$$
\sum_{H\in A} \lim_{n\ra\infty} \P(X_n\tensor \Z/a\Z \isom H) |\Sur(H,G)|=2^{\rank2(G)} |\wedge^2 G|.
$$



For each $G \in A$, we claim we can find an abelian group $G'\in A$ such that
$$
\sum_{H\in A} \frac{|\Hom(H,G)|}{|\Hom(H,G')|}
$$
converges.  We can factor over the primes $p$ dividing $a$, and reduce to the problem when $a=p^e$.  
Then if $G$ has type $\lambda$, we take $G'$ 
of type $\pi$ with $\pi_i'=2\lambda_i'+1$ for $1\leq i \leq e$.
Then using Lemma \ref{HomSize} we see that
$$
\sum_{H\in A} \frac{|\Hom(H,G)|}{|\Hom(H,G')|}=\sum_{c_1\geq  \dots \geq c_e \geq 0} p^{\sum_{i=1}^e c_i(\lambda_i'-2\lambda_i'-1)}
=\sum_{c_1\geq  \dots \geq c_e \geq 0} p^{\sum_{i=1}^e c_i(-\lambda_i'-1)}
$$
converges.

We have
\begin{align*}
\sum_{B\in A}\P(X_n\tensor \Z/a\Z \isom B) |\Hom(B,G')|&=\mathbb{E}|\Hom(X_n,G')|\\
&=\sum_{H<G'} \mathbb{E}|\Sur(X_n,H)|,
\end{align*}
and by supposition, each of the finite summands on the right-hand side has a finite
limit as $n\to \infty$ (and in particular is bounded above for all n). Thus, there is some
constant $D_G$ such that for all n we have

$$
\P(X_n\tensor \Z/a\Z \isom H) |\Hom(H,G')| \leq \sum_{H\in A} \P(X_n\tensor \Z/a\Z \isom H) |\Hom(H,G')| \leq D_G.
$$
Thus, for all $n$,
$$
\P(X_n\tensor \Z/a\Z \isom H) |\Hom(H,G)|\leq 
D_G|\Hom(H,G)| \cdot|\Hom(H,G')|^{-1} .
$$
Since
$
\sum_{H\in A } D_G |\Hom(H,G)|\cdot |\Hom(H,G')|^{-1} 
$
converges, by the Lebesgue Dominated Convergence Theorem, we have
\begin{multline*}
\sum_{H\in A } \lim_{n\ra\infty } \P(X_n\tensor \Z/a\Z \isom H) |\Hom(H,G)|\\
= \lim_{n\ra\infty } \sum_{H\in A } \P(X_n\tensor \Z/a\Z \isom H) |\Hom(H,G)|.
\end{multline*}
As this holds for every $G\in A$, we also have (by a finite number of additions and subtractions)
\begin{align*}
\sum_{H\in A } &\lim_{n\ra\infty } \P(X_n\tensor \Z/a\Z \isom H) |\Sur(H,G)|\\
&= \lim_{n\ra\infty } \sum_{H\in A } \P(X_n\tensor \Z/a\Z \isom H) |\Sur(H,G)|\\
&=  2^{\rank2(G)}|\wedge^2 G|.
\end{align*}

Next, we show that if for every $G\in A$, 
\begin{align*}
\sum_{H\in A } &\lim_{n\ra\infty } \P(X_n\tensor \Z/a\Z \isom H) |\Sur(H,G)|\\
&= \sum_{H\in A } \lim_{n\ra\infty } \P(Y_n\tensor \Z/a\Z \isom H) |\Sur(H,G)|\\
&=  2^{\rank2(G)}|\wedge^2 G|,
\end{align*}
then we have for every $H\in A$ that
\[
\lim_{n\ra\infty } \P(X_n\tensor \Z/a\Z \isom H)=\lim_{n\ra\infty } \P(Y_n\tensor \Z/a\Z \isom H).
\]
For each $G$, by a finite number of additions, we have 
\begin{align*}
\sum_{H\in A } &\lim_{n\ra\infty } \P(X_n\tensor \Z/a\Z \isom H) |\Hom(H,G)|\\
&= \sum_{H\in A } \lim_{n\ra\infty } \P(Y_n\tensor \Z/a\Z \isom H) |\Hom(H,G)|\\&= \sum_{G_1 \textrm{ subgroup of } G }  2^{\rank2(G_1)} |\wedge^2 G_1|.
\end{align*}

Now we will explain how to apply Theorem~\ref{T:Momdet} to conclude that\break 
$
\lim_{n\ra\infty } \P(X_n\tensor \Z/a\Z \isom H)=\lim_{n\ra\infty } \P(Y_n\tensor \Z/a\Z \isom H).
$
We factor $a=\prod_{j=1}^s p_j^{m_j}$. Since $a$ is even we may assume that $p_1=2$. The partition $\lambda^j\in M_j$ is the transpose of the type of the Sylow $p_j$-subgroup of $H$, which gives a bijection between $M$ and $A$.

Let $A_0=\{G\in A\quad|\quad\rank2(G)\text{ is odd}\}$. By restricting the bijection above we get a bijection between $M_0$ and $A_0$, where $M_0$ was defined in Lemma \ref{T:Momdet}.

We have that for $G\in A$ with corresponding $\lambda\in M$,
$$
C_\lambda=\sum_{G_1 \textrm{ subgroup of } G } 2^{\rank2(G_1)} |\wedge^2 G_1|\leq 2^{\lambda_1^1} \prod_{j=1}^s {F^{m_j}p_j^{\sum_i \frac{\lambda^{j}_i(\lambda^{j}_i-1)}{2}}}
$$
by Lemma~\ref{L:BoundMom}.  
For $H,G\in A$ with corresponding $\mu,\lambda\in M$, we have
$
| \Hom(H,G)| = \prod_{j=1}^s p_j^{\sum_i \lambda_i^j \mu_i^j }
$ by Lemma \ref{HomSize}.
 So for $H\in A_0$ with corresponding $\mu\in M_0$, we let
\[
x_\mu= \lim_{n\ra\infty } \P(X_n\tensor \Z/a\Z \isom H)
\]
and similarly for $y_\mu$ and we can apply Theorem~\ref{T:Momdet}.

Now, we suppose for the sake of contradiction that the limit\break 
$
\lim_{n\ra\infty} \P(X_n\tensor \Z/a\Z \isom H)
$
does not exist for at least some $H\in A_0$.  Then we can use a diagonal argument to find a subsequence
of $X_n$ where the limits do exist for all $H\in A_0$, and then another subsequence where the limits
do also exist for all $H \in A_0$, but at least one is different.
But since in each subsequence the limits
$
\lim_{n\ra\infty} \P(X_{i_n}\tensor \Z/a\Z \isom H)
$
exist,  we can use the argument above to conclude that these limits have to be the same for both subsequences, a contradiction.
\end{proof}

\end{document}